\documentclass[hidelinks,onefignum,onetabnum]{siamart251216}
\usepackage[T1]{fontenc}
\usepackage{amsfonts,amssymb,mathtools,booktabs,array,longtable,graphicx,enumitem}
\newsiamremark{remark}{Remark}
\newsiamremark{example}{Example}
\theoremstyle{plain}
\theoremheaderfont{\normalfont\scshape}
\theorembodyfont{\normalfont}
\theoremseparator{.}
\theoremsymbol{}
\newtheorem{problemstatement}{Problem}
\newenvironment{problem}[1]{
  \renewcommand{\theproblemstatement}{#1}
  \begin{problemstatement}
  \interlinepenalty=10000\relax
  \postdisplaypenalty=10000\relax
  \leavevmode\phantomsection\label{prob:#1}\ignorespaces
}{\end{problemstatement}}
\newcommand{\R}{\mathbb R}
\newcommand{\C}{\mathbb C}

\newcommand{\E}{\mathbb E}
\newcommand{\F}{\mathbb F}

\newcommand{\eps}{\varepsilon}
\newcommand{\norm}[1]{\lVert #1\rVert}

\newcommand{\opnorm}[1]{\lVert #1\rVert_2}

\DeclareMathOperator{\tr}{tr}
\DeclareMathOperator{\diag}{diag}
\DeclareMathOperator{\rank}{rank}

\DeclareMathOperator{\adj}{adj}

\newcommand{\one}{\mathbf 1}

\headers{Numerical Linear Algebra, Explained}{M. J. Colbrook and G. Stepaniants}
\title{Numerical Linear Algebra, Explained:\\Twenty Problem Resolutions}
\author{Matthew J. Colbrook\thanks{Department of Applied Mathematics and Theoretical Physics,
  University of Cambridge, Wilberforce Road, Cambridge CB3~0WA, United Kingdom
  (\email{m.colbrook@damtp.cam.ac.uk}, \email{gs852@cam.ac.uk}).}
\and George Stepaniants\footnotemark[1]\enspace\thanks{Also at the Courant Institute,
  New York University, New York, USA.}}
\begin{document}
\maketitle
\begin{abstract}
We resolve twenty problems from the ``Open Problems in Numerical Linear Algebra'' repository.
This work is part of the \emph{Numerical Linear Algebra, Explained} (\url{https://nla-explained.com/})
project, which collects explanations and human digestions of numerical linear algebra proofs.
(1,~2)~For products drawn from a family of $d\times d$ matrices, we prove local $1/d$-H\"older
continuity of their exponential growth rate (the joint spectral radius) and a uniform
$O(n^{d-1})$ bound on $n$-factor product norms at radius one.
(3)~For any $\gamma\geq0$, some pair of real matrices has largest length-$n$ product
norm of order $n^\gamma$ at every length.
(4,~5)~After $r$ random pivots, Cholesky and LU have sharp factors $2^r$ and $4^r$ relative
to the best rank-$r$ error, in expected trace norm and expected squared Frobenius norm,
respectively.
(6)~Even the simple concave function $f(x)=\min(x,1)$ can spoil a best approximation in
spectral norm and amplify the relative nuclear-norm excess without bound.
(7)~For a continuous, nonnegative, nondecreasing concave function, the sharp amplification
factor for the relative excess in Frobenius error is one under semidefinite order and two in general.
(8,~9)~In high dimensions with many more rows than columns, independent rows uniform on the
unit sphere asymptotically maximize the normalized smallest singular value under worst-case deletion of a fixed fraction of the rows.
(10)~If the prescribed tensor-train ranks do not allow an exact fit, a different choice
of first matrix truncation makes the usual squared-error bound strict.
(11)~Even with smaller border rank, the first saving for rank relative to a smooth
projective curve can occur beyond any prescribed tensor power.
(12)~The derivative of a nonnegative matrix's characteristic polynomial can have roots
that no nonnegative matrix realizes as its spectrum.
(13)~The Rowland--Wu identity gives a polynomial equation for an entry of a positive
rectangular matrix after Sinkhorn scaling to uniform row and column sums.
(14)~When matrix entries are expressed as trace products, the space of minimum-size positive
semidefinite factorizations can be disconnected modulo congruence.
(15)~A random cubic graph typically supports stable Kuramoto equilibria with different
phases even though neighboring phases attract.
(16,~17)~For tridiagonal systems with interval entries, finding the exact determinant
range or smallest box containing all solutions is NP-hard even with guaranteed nonsingularity.
(18,~19,~20)~Restricting the $q$-permanent to permutations that preserve a subset can
increase the sum, and two inequalities involving geometric means fail, one for norms
of sums and one for products of leading singular values in the weighted case.
\end{abstract}

\vspace{-3mm}

\section{Introduction}\label{sec:introduction}
We present twenty solutions to problems in numerical linear algebra, with formal verifications in Lean~\cite{palomar-2026-11-21-000005-v1}:

\vspace{1mm}

\noindent
\begin{center}
\textcolor[rgb]{0,0,1}{\url{https://github.com/MColbrook/NLA-formal-verification-20-problems}}
\end{center}

\vspace{1mm}

\noindent
The problems in this paper are drawn from Townsend's \emph{Open Problems in Numerical Linear Algebra} repository~\cite{TownsendRepository}.
Questions in this repository cover matrix products, tensor rank, factorization, synchronization, and matrix inequalities. At the time of writing, solutions have
been submitted to 131 of the 217 problems. This paper covers only a small part of the
collection, and we invite the reader to contribute explanations of other solutions.
Proofs to problems from this repository are collected at \url{https://nla-explained.com/}~\cite{NLAExplained}, which invites mathematicians at different stages to submit papers and videos with clearly explained solutions to problems in numerical linear algebra. 
As AI moves us into an era of proof abundance, there is now an even greater emphasis and need for clear explanations to mathematical problems.

\begin{table}[t]
\caption{The twenty questions and their resolutions.}\label{tab:roadmap}
\centering
\begin{tabular}{@{}ll>{\raggedright\arraybackslash}p{9cm}@{}}
\toprule
$\S$ & Problems & Result \\
\midrule
\ref{mfu-section} & \ref{prob:MF-05}, \ref{prob:MF-07} & Local $1/d$-H\"older continuity of the joint spectral radius and uniform bounds on product growth. \\
\ref{mfg-section} & \ref{prob:MF-12} & A pair of real matrices can realize any nonnegative exponent of product growth. \\
\ref{rap-section} & \ref{prob:RA-02}, \ref{prob:RA-03} & After $r$ random pivots, the sharp error factors are $2^r$ for Cholesky and $4^r$ for LU. \\
\ref{rac-section} & \ref{prob:RA-08} & A concave function can spoil spectral-norm optimality and amplify the relative nuclear-norm excess without bound. \\
\ref{raf-section} & \ref{prob:RA-09} & Sharp bounds on how much a concave function can increase relative Frobenius error excess. \\
\ref{sec:rows} & \ref{prob:IE-21}, \ref{prob:IE-22} & Random unit-norm rows attain the best asymptotic constant under adversarial row deletion. \\
\ref{sec:tr04} & \ref{prob:TR-04} & A tensor-train algorithm that makes the usual error bound strict whenever the optimal error is positive. \\
\ref{sec:tr27} & \ref{prob:TR-27} & For rank relative to a curve, savings can first occur beyond any prescribed tensor power. \\
\ref{sec:is03} & \ref{prob:IS-03} & The derivative roots of a characteristic polynomial need not have a nonnegative realization. \\
\ref{sec:nm04} & \ref{prob:NM-04} & A proof of the Rowland--Wu identity for rectangular Sinkhorn scaling. \\
\ref{sec:pf02} & \ref{prob:PF-02} & Disconnected quotients of minimal real positive semidefinite factorizations. \\
\ref{sec:md} & \ref{prob:MD-06} & Stable nonsynchronized equilibria with high probability on random cubic graphs. \\
\ref{sec:intervals} & \ref{prob:IV-02}, \ref{prob:IV-04} & NP-hardness of exact determinant ranges and solution hulls for regular tridiagonal interval matrices. \\
\ref{sec:mi19} & \ref{prob:MI-19} & A counterexample to the subset inequality for the $q$-permanent. \\
\ref{sec:mi21} & \ref{prob:MI-21} & A two-dimensional counterexample to a geometric-mean norm inequality. \\
\ref{sec:mi22} & \ref{prob:MI-22} & A counterexample to weighted log-majorization certified by error bounds for matrix roots. \\
\bottomrule
\end{tabular}\vspace{-5mm}
\end{table}

\subsection{Organization and notation}
We organize problems in such a way, that problems which use the same proof or construction share a section.
Table~\ref{tab:roadmap} summarizes the problem
statements and the section numbers link to their solutions.

We write $A^T$ for the transpose and $A^*$ for the conjugate transpose.
The norms $\norm{A}_2$, $\norm{A}_{\mathrm F}=(\tr(A^*A))^{1/2}$, and
$\norm{A}_*=\sum_j\sigma_j(A)$ are the spectral, Frobenius, and nuclear
norms, with the last also known as the trace norm. We order singular
values $\sigma_j(A)$ and the eigenvalues of positive semidefinite matrices
from largest to smallest, with multiplicity. When we seek a best rank-$k$
approximation, we minimize the stated error over matrices of rank at most
$k$. For a vector norm $p$, the induced norm $\norm{A}_p$ is the supremum
of $p(Ax)$ over vectors with $p(x)=1$. In the L\"owner order on Hermitian
matrices, $A\succeq B$ means that $A-B$ is positive semidefinite, while
$A\succ B$ requires positive definiteness. We call a matrix nonnegative
when all its entries are nonnegative.

\section*{Use of AI}
The solutions and proofs to these problems are generated using AI with several versions of ChatGPT 5 and 6, as documented in the repository. The authors take full responsibility for the writing of solutions in this paper.

\section{Uniform growth and the joint spectral radius (MF-05, MF-07)}
\label{mfu-section}
\begingroup
\newcommand{\mfuopnorm}[2]{\left\lVert #1\right\rVert_{#2}}
\newcommand{\mfuecc}{\operatorname{ecc}_2}

For one matrix $A$, the spectral radius gives the exponential growth
factor of the powers $A^n$. At spectral radius one there can still be
polynomial growth. For example, the norm of the $n$th power of a Jordan
block of order $d$ with eigenvalue one is of order $n^{d-1}$.

Now take a bounded family $\mathcal M$ of $d\times d$ matrices. We ask
how large the norm of $A_n\cdots A_1$ can get when each factor $A_i$ is
chosen separately from $\mathcal M$. With $\mathcal M=\{A\}$, every
factor is $A$ and the product is just $A^n$.

Let $\mathcal M\subset\F^{d\times d}$ be nonempty and bounded, where
$d\ge1$ and $\F$ is $\R$ or $\C$. We write $a_n(\mathcal M)$ for the
supremum of the norms of length-$n$ products. Its exponential growth
factor $\widehat\rho(\mathcal M)$ is the joint spectral radius,
\[
 a_n(\mathcal M)=\sup_{A_1,\ldots,A_n\in\mathcal M}
                   \opnorm{A_n\cdots A_1},\qquad a_0(\mathcal M)=1,
 \qquad
 \widehat\rho(\mathcal M)=\lim_{n\to\infty}a_n(\mathcal M)^{1/n}.
\]
The limit exists by submultiplicativity and is independent of the
choice of matrix norm. For compact families the supremum in $a_n$ is
attained, and for $\mathcal M=\{A\}$ the limit equals the ordinary
spectral radius $\rho(A)$.

For two nonempty compact families $\mathcal M$ and $\mathcal N$, their Hausdorff distance in the spectral
matrix norm is
\[
 d_H(\mathcal M,\mathcal N)
 =\max\left\{\sup_{A\in\mathcal M}\inf_{B\in\mathcal N}\opnorm{A-B},
             \sup_{B\in\mathcal N}\inf_{A\in\mathcal M}\opnorm{A-B}\right\}.
\]
Epperlein and Wirth posed the two questions in parts (L1) and (L3)
of Conjecture~3~\cite{EpperleinWirth2025}.

\begin{problem}{MF-05}
For every $d\ge1$ and every nonempty compact family $\mathcal M_0$ of
complex $d\times d$ matrices, do there exist $r,C>0$ such that
\[
 |\widehat\rho(\mathcal M)-\widehat\rho(\mathcal N)|
 \le C d_H(\mathcal M,\mathcal N)^{1/d}
\]
for all nonempty compact families $\mathcal M,\mathcal N$ of complex
$d\times d$ matrices with $d_H(\mathcal M,\mathcal M_0)<r$ and
$d_H(\mathcal N,\mathcal M_0)<r$?
\end{problem}

\begin{problem}{MF-07}
Does there exist a constant $\Theta_d>0$ depending only on $d\ge1$ such
that every nonempty compact family $\mathcal M$ of complex $d\times d$
matrices with $\widehat\rho(\mathcal M)=1$ satisfies
\[
 a_n(\mathcal M)\le\Theta_d(Ln)^{d-1}\quad(n\ge1),
 \qquad L=\max_{A\in\mathcal M}\opnorm A?
\]
\end{problem}

Wirth~\cite{Wirth2002} proved local Lipschitz continuity for irreducible compact families,
using extremal norms. Such families have no common nontrivial invariant
subspace. Epperlein and Wirth~\cite[Theorem~23]{EpperleinWirth2025}
proved pointwise H\"older continuity with exponent $1/(d^2+d)$ for all
compact families, and with exponent $1/(d+\varepsilon)$ at finite
families for every $\varepsilon>0$. In dimension two they obtained a uniform
product bound, as well as local $1/6$-H\"older continuity at families with positive
joint spectral radius.

In a local continuity estimate, one constant has to work for any two families
in the same neighborhood. We show that once we fix the dimension and a bound
on the matrix norms, the same constant works for all such families. This includes
real and complex matrices, reducible families, and families with infinitely many matrices.

\begin{theorem}[Uniform product bound]\label{mfu-thm-growth}
Let $\mathcal M\subset\F^{d\times d}$ be nonempty and bounded, and let $L>0$ be any upper bound for $\sup_{A\in\mathcal M}\opnorm{A}$. Set $r=\widehat\rho(\mathcal M)$ and define the constant
\[
 \Theta_d=d\left(\frac{2ed^2}{d-1}\right)^{d-1}\quad\text{for }d\geq 2,\qquad \Theta_1=1.
\]
If $r>0$, then for every $n\geq1$,
\begin{equation}\label{mfu-eq-general-growth}
 a_n(\mathcal M)\leq\Theta_d(Ln)^{d-1}r^{n-d+1}.
\end{equation}
If $r=0$, every product of length at least $d$ is zero.
In particular, when $r=1$,
\begin{equation}\label{mfu-eq-radius-one-growth}
 a_n(\mathcal M)\leq\Theta_d(Ln)^{d-1}\qquad(n\geq1).
\end{equation}
\end{theorem}

\begin{remark}
To see that the exponent $d-1$ of $n$ is sharp, take a single Jordan block
$J=I+N$ of order $d$ with eigenvalue one. The top-right entry of $J^n$ is
$\binom n{d-1}$, asymptotic to $n^{d-1}/(d-1)!$.
\end{remark}

\begin{theorem}[Uniform H\"older estimate]\label{mfu-thm-holder}
Let $\mathcal M,\mathcal N\subset\F^{d\times d}$ be nonempty compact families, both contained in the spectral-norm ball of radius $L>0$. Then
\begin{equation}\label{mfu-eq-holder}
 |\widehat\rho(\mathcal M)-\widehat\rho(\mathcal N)|
 \leq d(2d+1)L^{1-1/d}d_H(\mathcal M,\mathcal N)^{1/d}.
\end{equation}
For $d=1$, we can replace the coefficient $d(2d+1)$ with $1$.
\end{theorem}

\begin{remark}
A nilpotent shift shows that the H\"older exponent $1/d$ is sharp. Let $N$ be the nilpotent Jordan shift of order $d$, and let $E_{d1}$ have a one in position $(d,1)$ and zeros elsewhere. For $0<\delta\leq1$, the matrix $N+\delta E_{d1}$ has characteristic
polynomial $z^d-\delta$ and spectral norm at most one. Its eigenvalues
give the perturbed radius, and the perturbation norm $\delta$ gives
the Hausdorff distance. In particular,
\[
 \widehat\rho(\{N\})=0,\qquad
 \widehat\rho(\{N+\delta E_{d1}\})=\delta^{1/d},\qquad
 d_H(\{N\},\{N+\delta E_{d1}\})=\delta.
\]
A H\"older estimate of exponent $\alpha$ at $\{N\}$ would require
$\delta^{1/d}/\delta^\alpha$ to stay bounded as $\delta\downarrow0$.
This fails for $\alpha>1/d$, so in particular there is no local Lipschitz
bound at $\{N\}$ when $d\geq2$.
\end{remark}

\subsection{Choosing coordinates for an approximate extremal norm}
Let $a$ be larger than the joint spectral radius. We can bound the norm of each
matrix in the family by $a$ by measuring vectors in a norm built from all products.
The difficulty is comparing this norm with the Euclidean norm. Choose a basis
from the unit ball that maximizes the absolute determinant and take its
singular value decomposition. We can put all the large scale factors into
the diagonal matrix $D$ and work with a norm that lies within a factor $d$
of the Euclidean norm in the new coordinates.
The connection between this Auerbach-basis construction and the joint
spectral radius is discussed in \cite{EpperleinWirthProjection2025}.

\begin{lemma}[A well-rounded approximate extremal norm]\label{mfu-lem-rounded}
Let $\mathcal M$ and $L$ be as in Theorem~\ref{mfu-thm-growth}. For every $a>\widehat\rho(\mathcal M)$, there exist a unitary matrix $Q$ (orthogonal over $\R$), a diagonal matrix
\[
 D=\diag(\sigma_1,\ldots,\sigma_d),\qquad
 \sigma_1\geq\cdots\geq\sigma_d>0,
\]
and a norm $w$ on $\F^d$ such that, for every $z\in\F^d$ and $A\in\mathcal M$,
\begin{equation}\label{mfu-eq-rounding}
 \opnorm{z}\leq w(z)\leq d\opnorm{z},
 \qquad \mfuopnorm{B_A}{w}\leq a,
 \qquad B_A=D^{-1}Q^*AQD,
\end{equation}
and
\begin{equation}\label{mfu-eq-entry-old}
 |(B_A)_{ij}|\leq L\frac{\sigma_j}{\sigma_i}
 \quad(A\in\mathcal M).
\end{equation}
\end{lemma}
\begin{proof}
Define
\[
 v(x)=\sup_{k\geq0}\sup_{A_1,\ldots,A_k\in\mathcal M}
               a^{-k}\opnorm{A_k\cdots A_1x}.
\]
Choose $b$ strictly between $\widehat\rho(\mathcal M)$ and $a$.
The root limit gives $a_k(\mathcal M)\leq b^k$ for all sufficiently
large $k$, so $a^{-k}a_k(\mathcal M)$ is bounded. Thus $v$ is a finite
norm equivalent to the Euclidean norm. The $k=0$ term gives $v(x)\geq\opnorm{x}$, and appending a generator gives $v(Ax)\leq av(x)$.

Choose a matrix $T$ whose columns lie in the $v$-unit ball and maximize
$|\det T|$. The ball is compact and has nonempty interior, so a maximizer
exists and has nonzero determinant. Each column has $v$-norm one, since otherwise
scaling that column up to norm one would increase $|\det T|$.
For $x=Tz$ in the unit ball, replacing column $j$ of $T$ by $x$
multiplies the determinant by $z_j$, so maximality implies $|z_j|\leq1$.
Scaling to the unit ball gives the lower bound for any $z$.
The upper bound follows by writing $Tz$ in terms of the columns of $T$ and applying
the triangle inequality. Thus
\[
 \lVert z\rVert_\infty\leq v(Tz)\leq\lVert z\rVert_1.
\]

Take a singular value decomposition $T=QDR^*$, with the singular values in decreasing order, and define $w(z)=\sqrt d\,v(QDz)$. Since $QD=TR$ and $R$ is unitary,
\[
 \frac{\opnorm{z}}{\sqrt d}\leq v(TRz)\leq\sqrt d\opnorm{z}.
\]
This proves the first part of \eqref{mfu-eq-rounding}; the operator-norm bound follows from $v(Ax)\leq av(x)$. Finally, each entry of $Q^*AQ$ has modulus at most $L$, proving \eqref{mfu-eq-entry-old}.
\end{proof}

\subsection{Truncating large changes of scale}
The change of coordinates can be poorly conditioned when successive diagonal entries
of $D$ have a large ratio. Capping each ratio at $s\geq1$ bounds the condition number
by $s^{d-1}$. Make a block boundary at each ratio larger than $s$. By
\eqref{mfu-eq-entry-old}, the entries above these boundaries are small, so the matrices
are close to block lower triangular ones. The following lemma controls what happens
to their products when we rescale the blocks.

\begin{lemma}[Block diagonal similarity]\label{mfu-lem-damping}
Suppose $X$ is block lower triangular for a fixed ordered block decomposition, and
\[
 H=\diag(h_1I_1,\ldots,h_bI_b),\qquad
 h_1\geq\cdots\geq h_b>0.
\]
Here $I_j$ is the identity on block $j$, for $1\leq j\leq b$. Then
\begin{equation}\label{mfu-eq-damping}
 \opnorm{HXH^{-1}}\leq\opnorm{X}.
\end{equation}
\end{lemma}
\begin{proof}
For each cut $\ell=1,\ldots,b-1$, set $t_\ell=h_{\ell+1}/h_\ell\in(0,1]$ and let $S_\ell$ be $+I$ on the first $\ell$ blocks and $-I$ on the remaining blocks. The map
\[
 \Phi_\ell(Y)=\frac{1+t_\ell}{2}Y
              +\frac{1-t_\ell}{2}S_\ell YS_\ell
\]
is a convex combination of unitary conjugations, so $\opnorm{\Phi_\ell(Y)}\leq\opnorm{Y}$. It multiplies blocks crossing the cut by $t_\ell$ and leaves all other blocks unchanged.

Composing these maps multiplies a lower block $(i,j)$, $i\geq j$, by
$\prod_{\ell=j}^{i-1}t_\ell=h_i/h_j$. Since the upper blocks of $X$
are zero, the composite maps $X$ to $HXH^{-1}$, proving
\eqref{mfu-eq-damping}. For a single block, the assertion is equality.
\end{proof}

\begin{proposition}\label{mfu-prop-comparison}
Under the hypotheses of Theorem~\ref{mfu-thm-growth}, set $r=\widehat\rho(\mathcal M)$. For every $s\geq1$ and every integer $n\geq0$,
\begin{equation}\label{mfu-eq-comparison}
 a_n(\mathcal M)\leq d s^{d-1}
                \left(r+\frac{2d^2L}{s}\right)^n.
\end{equation}
\end{proposition}
\begin{proof}
Fix $a>r$ and obtain $Q,D,w$ from Lemma~\ref{mfu-lem-rounded}. Place a block
boundary between coordinates $i$ and $i+1$ whenever $\sigma_i/\sigma_{i+1}>s$.
Subtract the strictly upper block part $E_A$ of $B_A$ to obtain the
block lower triangular matrix $C_A=B_A-E_A$.

Each potentially nonzero entry of $E_A$ crosses a gap greater than $s$. By \eqref{mfu-eq-entry-old}, its modulus is at most $L/s$. A $d\times d$ matrix with entries bounded by $L/s$ has spectral norm at most $dL/s$. Hence
\[
 \mfuopnorm{E_A}{w}\leq d\opnorm{E_A}\leq\frac{d^2L}{s},
 \qquad \mfuopnorm{C_A}{w}\leq a+\frac{d^2L}{s}=:u.
\]
For a product of $t\geq0$ factors $C_A$, the operator norm induced by $w$
is at most $u^t$. Converting back to the spectral norm with
\eqref{mfu-eq-rounding} gives
\begin{equation}\label{mfu-eq-C-old}
 \opnorm{C_{A_t}\cdots C_{A_1}}\leq d u^t.
\end{equation}
For $t=0$, the product is the identity and satisfies the same bound.

Define $D'=\diag(\sigma'_1,\ldots,\sigma'_d)$ by $\sigma'_d=1$ and
\[
 \frac{\sigma'_i}{\sigma'_{i+1}}
 =\min\left\{\frac{\sigma_i}{\sigma_{i+1}},s\right\}.
\]
Then
\begin{equation}\label{mfu-eq-condition}
 \opnorm{D'}\opnorm{D'^{-1}}\leq s^{d-1}.
\end{equation}
The diagonal matrix $H=D'^{-1}D$ is scalar on each block. Its block scalars are nonincreasing: within a block, consecutive diagonal entries of $H$ have ratio one, while across a wide gap their ratio is $(\sigma_i/\sigma_{i+1})/s>1$.

Set
\[
 C'_A=HC_AH^{-1},\quad E'_A=HE_AH^{-1},\quad
 B'_A=C'_A+E'_A=D'^{-1}Q^*AQD'.
\]
The product of any $t$ of the $C_A$ is block lower triangular. Lemma~\ref{mfu-lem-damping} and \eqref{mfu-eq-C-old} therefore give
\begin{equation}\label{mfu-eq-C-new}
 \opnorm{C'_{A_t}\cdots C'_{A_1}}\leq d u^t
 \qquad(t\geq0).
\end{equation}
Every nonzero entry of $E'_A$ crosses a wide gap, which we replace with $s$ in $D'$. All other consecutive ratios in $D'$ are at least one. Directly from $B'_A=D'^{-1}Q^*AQD'$, this gives
\begin{equation}\label{mfu-eq-E-new}
 |(E'_A)_{ij}|\leq L/s,\qquad
 \opnorm{E'_A}\leq dL/s=:v.
\end{equation}

Expand a length-$n$ product of $B'_A=C'_A+E'_A$, grouping terms by the
positions occupied by $E'$. For $k$ such positions,
the remaining factors form $k+1$ consecutive $C'$-products, some possibly
empty, with total length $n-k$. Estimate each of these products by
\eqref{mfu-eq-C-new}, and each $E'$ by \eqref{mfu-eq-E-new}. The resulting
bound is $d^{k+1}u^{n-k}v^k$. There are $\binom nk$ ways to choose the
positions, so the binomial theorem gives
$$
 \opnorm{B'_{A_n}\cdots B'_{A_1}}
 \leq \sum_{k=0}^n\binom nk d^{k+1}u^{n-k}v^k=d(u+dv)^n
 =d\left(a+\frac{2d^2L}{s}\right)^n.
$$
Undoing the similarity costs at most the factor $s^{d-1}$ from
\eqref{mfu-eq-condition}, giving
\[
 a_n(\mathcal M)\leq d s^{d-1}
                \left(a+\frac{2d^2L}{s}\right)^n.
\]
The bound involves only $d,s,L$ and $a$, so letting $a\downarrow r$ gives \eqref{mfu-eq-comparison}.
\end{proof}

\begin{proof}[Proof of Theorem~\ref{mfu-thm-growth}]
For $d=1$, $r=\sup_{A\in\mathcal M}|A|$ and scalar multiplication proves the assertions. Let $d\geq2$ and set $m=d-1$. For $r>0$, keep the accumulated error bounded by choosing
\[
 s=\frac{2d^2Ln}{mr}\geq1
\]
in Proposition~\ref{mfu-prop-comparison}. Since $r\leq L$ and $(1+m/n)^n\leq e^m$,
$$
 a_n(\mathcal M)
 \leq d\left(\frac{2d^2Ln}{mr}\right)^m
                  r^n(1+m/n)^n
 \leq d\left(\frac{2ed^2}{m}\right)^m
                  (Ln)^m r^{n-m},
$$
which is \eqref{mfu-eq-general-growth}. If $r=0$, the same proposition gives
\[
 a_n(\mathcal M)\leq d(2d^2L)^n s^{d-1-n}\qquad(s\geq1).
\]
For $n\geq d$, letting $s\to\infty$ gives $a_n(\mathcal M)=0$.
\end{proof}

\subsection{Bounds for perturbed families}
For a vector norm $p$,
define its eccentricity relative to the Euclidean norm by

\[
 \mfuecc(p)=\frac{\sup_{\opnorm{x}=1}p(x)}{\inf_{\opnorm{x}=1}p(x)}.
\]
\begin{corollary}\label{mfu-cor-ecc}
Under the hypotheses of Theorem~\ref{mfu-thm-growth}, set $r=\widehat\rho(\mathcal M)$. For every $0<\varepsilon\leq2d^2L$ there is a norm $p_\varepsilon$ satisfying
\begin{equation}\label{mfu-eq-ecc}
 \sup_{A\in\mathcal M}\mfuopnorm A{p_\varepsilon}
 \leq r+\varepsilon,
 \qquad
 \mfuecc(p_\varepsilon)\leq d\left(\frac{2d^2L}{\varepsilon}\right)^{d-1}.
\end{equation}
\end{corollary}
\begin{proof}
Set $s=2d^2L/\varepsilon\geq1$ and $h=r+\varepsilon$. Proposition~\ref{mfu-prop-comparison} gives $a_k(\mathcal M)\leq d s^{d-1}h^k$ for every $k\geq0$. Consequently,
\[
 p_\varepsilon(x)=\sup_{k\geq0}\sup_{A_1,\ldots,A_k\in\mathcal M}
                     h^{-k}\opnorm{A_k\cdots A_1x}
\]
is a finite norm satisfying
\[
 \opnorm{x}\leq p_\varepsilon(x)\leq d s^{d-1}\opnorm{x},
 \qquad p_\varepsilon(Ax)\leq h p_\varepsilon(x).
\]
\end{proof}

\begin{proof}[Proof of Theorem~\ref{mfu-thm-holder}]
Set $\delta=d_H(\mathcal M,\mathcal N)$. If $\delta=0$, the two compact
families coincide. For $0<\delta\leq L$, take $s\geq1$ and set
$\varepsilon=2d^2L/s\leq2d^2L$.

Each $B\in\mathcal N$ lies within $\delta$ in spectral norm of some
$A\in\mathcal M$. The eccentricity bound in Corollary~\ref{mfu-cor-ecc}
gives $\mfuopnorm{B-A}{p_\varepsilon}\leq d s^{d-1}\delta$. Using the same
corollary to bound $\mfuopnorm{A}{p_\varepsilon}$, we get
\begin{align*}
 \widehat\rho(\mathcal N)
 \leq\sup_{B\in\mathcal N}\mfuopnorm B{p_\varepsilon}
 \leq\widehat\rho(\mathcal M)+2d^2L/s+d s^{d-1}\delta.
\end{align*}
Choose $s=(L/\delta)^{1/d}\geq1$ so that $L/s=\delta s^{d-1}$. The bound
then becomes
\[
 \widehat\rho(\mathcal N)\leq
 \widehat\rho(\mathcal M)+d(2d+1)L^{1-1/d}\delta^{1/d}.
\]
Interchanging the families proves \eqref{mfu-eq-holder}. If $\delta\geq L$, both radii belong to $[0,L]$, which gives the bound. For $d=1$, $\widehat\rho(\mathcal M)=\sup_{a\in\mathcal M}|a|$, a 1-Lipschitz function for Hausdorff distance.
\end{proof}

\endgroup

\section{Prescribed marginal growth of matrix products (MF-12)}
\label{mfg-section}
\begingroup
For real square matrices $A$ and $B$ of the same order, let
$a_n(\{A,B\})$ be the largest spectral norm of a product formed by
choosing either matrix at each of $n$ positions. At joint spectral
radius one, $a_n(\{A,B\})^{1/n}$ tends to one, and the maximum product
norm can still grow like a power of $n$. Varney and Morris~\cite{VarneyMorris2025}
call the sequence of maxima at radius one the marginal instability rate
sequence. They asked whether each power $n^\gamma$ with $\gamma\ge0$ can
occur as the growth rate for some finite family of matrices.

Protasov and Jungers~\cite{ProtasovJungers2015} constructed pairs of
$3\times3$ matrices for which $a_n$ is bounded above by a fixed
multiple of $n^{1/3}$ at every length and below by a positive multiple
of $n^{1/3}$ along an infinite subsequence. Varney and Morris
extended these bounds to exponents in $[1/3,1/2)$~\cite[Theorem~5]{VarneyMorris2025}.
Their Corollary~6.1 gives a pair of $6\times6$ matrices with both
bounds of order $n^{1/3}$ at every length.

\begin{problem}{MF-12}
For every real $\gamma\ge0$, does there exist a finite nonempty family
$\mathcal M$ of real square matrices of some common size with
$\widehat\rho(\mathcal M)=1$ and constants $0<c\le C<\infty$ such that
\[
 c n^\gamma\le a_n(\mathcal M)\le C n^\gamma\quad(n\ge1)?
\]
This is Question~2 of Varney and Morris~\cite{VarneyMorris2025}.
\end{problem}

\begin{theorem}\label{mfg-thm-main}
For every $\gamma\geq0$, there are two distinct real square matrices forming a family $\mathcal M_\gamma$, and constants $0<c_\gamma\leq C_\gamma<\infty$, such that
\[
 \widehat\rho(\mathcal M_\gamma)=1,\qquad
 c_\gamma n^\gamma\leq a_n(\mathcal M_\gamma)
                  \leq C_\gamma n^\gamma\quad(n\geq1).
\]
For $0<\gamma<1$, the matrices have order six. For noninteger $\gamma>0$, order $6(\lfloor\gamma\rfloor+1)$ suffices. For integer $\gamma=m\geq0$, order $m+1$ suffices. If $\gamma$ is rational, the matrices can have dyadic-rational entries.
\end{theorem}

\subsection{The six-dimensional construction and an upper bound}
We first construct a pair with growth exponent between zero and one. Choose
\[
 0<\lambda\leq\frac14,\qquad \lambda<\mu<1,\qquad
 \alpha=1-\frac{\log\mu}{\log\lambda}\in(0,1).
\]
Using the equivalent relation $\mu=\lambda^{1-\alpha}$, define
\[
 J_t=t\begin{pmatrix}1&1\\0&1\end{pmatrix},\qquad
 A=\diag(1,J_\lambda,J_\mu,1),
\]
so that the two Jordan blocks in $A$ have powers that decay at different
rates. For the second generator, define $V\in\R^{6\times2}$ and
$U\in\R^{2\times6}$ by
\[
 V=\begin{pmatrix}
 1&0\\0&0\\1&0\\0&0\\0&1\\0&1
 \end{pmatrix},\qquad
 U=\begin{pmatrix}
 1&-1&0&1&0&0\\
 0&0&0&0&0&1
 \end{pmatrix}.
\]
The second generator is the fixed matrix
\begin{equation}\label{mfg-eq-P}
 P=VU=\begin{pmatrix}
 1&-1&0&1&0&0\\
 0&0&0&0&0&0\\
 1&-1&0&1&0&0\\
 0&0&0&0&0&0\\
 0&0&0&0&0&1\\
 0&0&0&0&0&1
 \end{pmatrix}.
\end{equation}
Multiplication gives $UV=I_2$ and hence $P^2=P$. Thus $P$ is a
rank-two projection, though not an orthogonal one.

A vector in the range of $P$ has the form $Vy$ for $y\in\R^2$. After
applying $A^q$ and projecting again, we get $V(UA^qV)y$. Thus, in the
coordinates $y$, the effect of taking $q$ steps with $A$ and then projecting
is given by $UA^qV$.

\begin{lemma}[Compressed powers]\label{mfg-lem-compressed}
For every integer $q\geq0$,
\begin{equation}\label{mfg-eq-T}
 T_q:=UA^qV=\begin{pmatrix}1-q\lambda^q&q\mu^q\\0&1\end{pmatrix}.
\end{equation}
With $\ell_q=q\lambda^q$ and $b_q=q\mu^q$, these entries satisfy
\[
 \ell_0=b_0=0,\qquad 0<\ell_q\leq\tfrac14\quad(q\geq1),
 \qquad b_q=q^\alpha\ell_q^{1-\alpha}\quad(q\geq1).
\]
\end{lemma}
\begin{proof}
The Jordan identity
\[
 J_t^q=t^q\begin{pmatrix}1&q\\0&1\end{pmatrix}
\]
holds for every $q\geq0$. Substitution in $UA^qV$ gives \eqref{mfg-eq-T}, including $T_0=I_2$. Since $q4^{-q}\leq1/4$ for $q\geq1$ and $\lambda\leq1/4$, the estimate for $\ell_q$ follows. Finally,
\[
 q^\alpha(q\lambda^q)^{1-\alpha}
 =q\lambda^{q(1-\alpha)}=q\mu^q.
\]
\end{proof}

Consider a word of length $n$ containing $s\geq1$ copies of $P$. It has a representation
\begin{align}
 W=A^{q_s}P A^{q_{s-1}}P\cdots P A^{q_0}
  =A^{q_s}V\, T_{q_{s-1}}\cdots T_{q_1}\,U A^{q_0},
 \qquad \sum_{i=0}^s q_i+s=n.\label{mfg-eq-factor}
\end{align}
The middle product is empty when $s=1$.

\begin{lemma}[Bound for products of compressed powers]\label{mfg-lem-budget}
For any finite list of nonnegative integers $q_1,\ldots,q_k$, the product has the form
\[
 T_{q_k}\cdots T_{q_1}=\begin{pmatrix}a&z\\0&1\end{pmatrix}.
\]
Then $0\leq a\leq1$ and
\[
 0\leq z\leq\left(\sum_{i=1}^k q_i\right)^\alpha.
\]
For the empty list, the product is $I_2$ and $z=0$.
\end{lemma}
\begin{proof}
Set
\[
 w_i=\prod_{j=i+1}^k(1-\ell_{q_j}).
\]
By triangular multiplication,
\[
 a=\prod_{i=1}^k(1-\ell_{q_i}),\qquad
 z=\sum_{i=1}^k b_{q_i}w_i.
\]
The weight $w_i$ records the later contractions of the contribution
from gap $i$; it lies in $[0,1]$. Taking differences of successive
products of the factors $1-\ell_{q_i}$ gives
\begin{equation}\label{mfg-eq-telescoping}
 \sum_{i=1}^k\ell_{q_i}w_i
 =1-\prod_{i=1}^k(1-\ell_{q_i})\leq1.
\end{equation}
Since the terms with $q_i=0$ vanish, applying H\"older's inequality to the
remaining terms with conjugate exponents $1/\alpha$ and $1/(1-\alpha)$ gives
\begin{align*}
 z=\sum_i(q_iw_i)^\alpha(\ell_{q_i}w_i)^{1-\alpha}
  \leq\left(\sum_iq_iw_i\right)^\alpha
        \left(\sum_i\ell_{q_i}w_i\right)^{1-\alpha}
   \leq\left(\sum_iq_i\right)^\alpha.
\end{align*}
\end{proof}

By Lemma~\ref{mfg-lem-budget}, the middle product in \eqref{mfg-eq-factor} has spectral norm at most $1+n^\alpha$: decompose it as $\diag(a,1)+zE_{12}$, where $E_{12}$ is the $2\times2$ matrix with its only nonzero entry equal to one in position $(1,2)$, and apply the triangle inequality.

To bound the powers of $A$ in the end factors of \eqref{mfg-eq-factor},
use $\opnorm{J_t^r}\leq(r+1)t^r\leq\sum_{j=0}^rt^j\leq(1-t)^{-1}$ for
$t=\lambda,\mu$. As $\lambda<\mu<1$, the block diagonal form of $A$ gives
\begin{equation}\label{mfg-eq-boundedA}
 \opnorm{A^r}\leq H,\qquad H=(1-\mu)^{-1}\quad(r\geq0).
\end{equation}
Also
\[
 \opnorm{V}=\sqrt2,\qquad \opnorm{U}=\sqrt3.
\]

Consequently,
\[
 \opnorm{W}\leq\sqrt6 H^2(1+n^\alpha)
             \leq2\sqrt6 H^2n^\alpha.
\]
The word $A^n$, which has no $P$, satisfies the same final bound by \eqref{mfg-eq-boundedA}. Thus
\begin{equation}\label{mfg-eq-upper}
 a_n(\{A,P\})\leq C_{\lambda,\alpha}n^\alpha,
 \qquad C_{\lambda,\alpha}=\frac{2\sqrt6}{(1-\mu)^2}.
\end{equation}

\subsection{A lower bound at every length}
For the lower bound, repeat a run of $A$'s of the same length, applying $P$ after
each run. Suppose first that $n\geq1/\lambda$ and choose
\[
 q=\left\lfloor\log_{1/\lambda}n\right\rfloor\geq1,
 \qquad k=\left\lfloor\frac{n}{q+1}\right\rfloor,
 \qquad r=n-k(q+1).
\]
Then
\begin{equation}\label{mfg-eq-wordlower}
 W_n=A^r(PA^q)^k
\end{equation}
is an allowed word of exactly length $n$.

To estimate its norm, use $(PA^q)^kV=VT_q^k$, where
\[
 T_q^k=\begin{pmatrix}
 (1-\ell_q)^k&\dfrac{b_q}{\ell_q}\bigl(1-(1-\ell_q)^k\bigr)\\[3pt]
 0&1
 \end{pmatrix}.
\]
For every $(x,y)\in\R^2$, the first coordinate of $A^rV(x,y)^T$ equals $x$. For $e_2=(0,1)^T$, the vector $Ve_2$ has norm $\sqrt2$. Applying $W_n$ to this vector therefore gives
\begin{equation}\label{mfg-eq-lowerRayleigh}
 \opnorm{W_n}\geq\frac1{\sqrt2}\frac{b_q}{\ell_q}
                 \bigl(1-(1-\ell_q)^k\bigr).
\end{equation}
From the definition of $q$,
\[
 \lambda^{-q}\leq n<\lambda^{-(q+1)},\qquad
 n\geq\lambda^{-q}\geq4^q\geq2(q+1).
\]
It follows that $k\geq n/[2(q+1)]$ and $\lambda^q\geq1/n$, so
\[
 k\ell_q\geq\frac{q}{2(q+1)}\geq\frac14.
\]
Using $1-u\leq e^{-u}$ for $u\geq0$, we obtain
\[
 1-(1-\ell_q)^k\geq1-e^{-1/4}.
\]
Furthermore,
\[
 \frac{b_q}{\ell_q}=(\mu/\lambda)^q
  =(\lambda^{-q})^\alpha\geq(\lambda n)^\alpha.
\]
Substituting in \eqref{mfg-eq-lowerRayleigh} proves
\begin{equation}\label{mfg-eq-lower}
 a_n(\{A,P\})\geq c_{\lambda,\alpha}n^\alpha,
 \qquad c_{\lambda,\alpha}
       =\frac{1-e^{-1/4}}{\sqrt2}\lambda^\alpha.
\end{equation}
For the remaining lengths $1\leq n<1/\lambda$, the eigenvalue one of
$A$ gives $\opnorm{A^n}\geq1>c_{\lambda,\alpha}n^\alpha$, so
\eqref{mfg-eq-lower} also holds in this range.

The upper bound \eqref{mfg-eq-upper} implies $\widehat\rho(\{A,P\})\leq1$. Since $A$ has eigenvalue one, $a_n(\{A,P\})\geq1$ for every $n$, and the reverse inequality follows. For any prescribed $\alpha\in(0,1)$, setting $\lambda=1/4$ and $\mu=4^{-(1-\alpha)}$ gives the required pair.

\subsection{All exponents and rational entries}
Tensoring with a Jordan block supplies the integer part of the exponent. We first bound its powers for every $n\geq1$.
\begin{lemma}\label{mfg-lem-Jordan}
Let $m\geq0$ be an integer and let $J=I_{m+1}+N$ be the Jordan block with eigenvalue
one, where $N$ has ones on the first superdiagonal and zeros elsewhere. Then
$\opnorm{J^n}\asymp n^m$ for $n\geq1$, with positive comparison constants independent of $n$.
If $m\geq1$, explicit bounds are
\[
 m^{-m}n^m\leq\opnorm{J^n}\leq(m+1)n^m.
\]
For $m=0$, $\opnorm{J^n}=1$.
\end{lemma}
\begin{proof}
The identity $J^n=\sum_{j=0}^m\binom njN^j$, with $\binom nj=0$ for $j>n$, gives the upper bound. When $n\geq m\geq1$, the top-right entry is $\binom nm$, and
\[
 \binom nm=\prod_{j=0}^{m-1}\frac{n-j}{m-j}
           \geq(n/m)^m.
\]
For $1\leq n<m$, use $\opnorm{J^n}\geq1\geq(n/m)^m$.
\end{proof}

For noninteger $\gamma>0$, write $\gamma=m+\alpha$ with $m\geq0$ integer and $0<\alpha<1$. As in the closure argument of Varney and Morris~\cite[Proposition~3.1]{VarneyMorris2025}, tensor both generators with the Jordan block $J$ of order $m+1$:
\[
 \widehat A=A\otimes J,\qquad \widehat P=P\otimes J.
\]
Whatever choices we make in a length-$n$ product, the second tensor factor
is always $J^n$. Taking singular value decompositions shows that the spectral
norm of a tensor product is the product of the two norms. Thus
\[
 a_n(\{\widehat A,\widehat P\})
  =a_n(\{A,P\})\opnorm{J^n}\asymp n^{\alpha+m}.
\]
Taking $n$th roots gives joint spectral radius one.

For integer $\gamma=m\geq0$, use the pair $\{J,0_{(m+1)\times(m+1)}\}$. It has the same maximal growth as $J$ alone and consists of two distinct matrices. This proves the existence and dimension assertions of Theorem~\ref{mfg-thm-main}.

\begin{corollary}\label{mfg-cor-rational}
Two matrices with dyadic-rational entries realize every nonnegative rational exponent.
\end{corollary}
\begin{proof}
For integer $\gamma$ the pair above is already integral. Otherwise write
$\gamma=m+\alpha$, where $\alpha=a/b$ with integers $1\leq a<b$, and take
\[
 \lambda=2^{-b},\qquad \mu=2^{-(b-a)}.
\]
We have $\lambda\leq1/4$ since $b\geq2$, and $\mu=\lambda^{1-\alpha}$.
The denominators of $\lambda$ and $\mu$ both divide $2^b$, so $2^b A$
has integer entries. Since $P$ and $J$ are integral, so are
$2^b\widehat A=(2^b A)\otimes J$ and $\widehat P=P\otimes J$.
\end{proof}

\begin{corollary}\label{mfg-cor-irrational}
The largest length-$n$ product norm of a rational matrix pair can grow like $n^\alpha$
for irrational $\alpha$. The exponents attainable by rational pairs are dense in $[0,\infty)$.
\end{corollary}
\begin{proof}
Take $\lambda=1/4$ and $\mu=1/3$. Both matrices are rational, and the exponent is
\[
 \alpha=1-\frac{\log3}{\log4}\in(0,1).
\]
If $\log3/\log4=p/q$ for positive integers $p,q$, then $3^q=4^p$, contrary to unique factorization. Density follows from Corollary~\ref{mfg-cor-rational}, since the nonnegative rationals are dense in $[0,\infty)$.
\end{proof}

\endgroup

\section{Sharp worst-case losses in random pivoting (RA-02, RA-03)}\label{rap-section}
\begingroup
\newcommand{\rapind}{\mathbf 1}
In randomly pivoted Cholesky and LU, each pivot adds one rank-one term
to the approximation of $A$, giving rank at most $r$ after $r$ pivots.
Each step starts with the residual, obtained by subtracting the
current approximation from $A$. Whenever the residual is nonzero,
Cholesky samples a nonzero diagonal entry, while LU samples a nonzero
entry anywhere in the residual. For a pivot at $(i,j)$, the rank-one
term is the outer product of column $j$ and row $i$ of the residual,
divided by the pivot entry.

We compare the expected error after $r$ pivots with the best rank-$r$
error. The ratio is at most $2^r$ for Cholesky in trace norm and at
most $4^r$ for LU in squared Frobenius norm. We
construct a single family of positive definite matrices for which both
ratios approach their respective bounds.

For Cholesky, let $A\in\C^{n\times n}$ be Hermitian positive
semidefinite and take an integer $0\leq r\leq n$. Keeping the $r$
largest eigenvalues and their eigenvectors gives a best rank-$r$
approximation in the trace norm. Its error is the sum of the
remaining eigenvalues,
\[
 \tau_r(A)=\sum_{j>r}\lambda_j(A),
 \qquad \lambda_1(A)\ge\cdots\ge\lambda_n(A)\ge0.
\]
In both algorithms, once the residual reaches zero, all subsequent
residuals are zero. In exact arithmetic, randomly pivoted Cholesky
(RPCholesky) starts from $R_0=A$. Given a nonzero residual $R_s$
after $s$ steps, it chooses the diagonal pivot $(j,j)$ with
probability $(R_s)_{jj}/\tr R_s$ and sets
\begin{equation}\label{rap-eq:chol-update}
 R_{s+1}=R_s-\frac{(R_s)_{:j}(R_s)_{j:}}{(R_s)_{jj}}.
\end{equation}
Here a colon in an index denotes all rows or all columns. The
residual stores the Schur complement on the unpivoted indices,
with zeros on the pivoted rows and columns. It remains positive
semidefinite, so the error in the trace norm is $\tr R_s$. For
RPCholesky the known bound is \cite[Lemma 5.5]{ChenEpperlyTroppWebber2025}
\begin{equation}\label{rap-eq:known-chol}
 \E\tr R_r\le 2^r\tau_r(A).
\end{equation}

Randomly pivoted LU (RPLU) starts from $B_0=A$ for a general rectangular
complex matrix $A\in\C^{m\times n}$. Given $B_s\ne0$, it chooses $(i,j)$ with probability
$|(B_s)_{ij}|^2/\norm{B_s}_{\mathrm F}^2$ and sets
\begin{equation}\label{rap-eq:lu-update}
 B_{s+1}=B_s-\frac{(B_s)_{:j}(B_s)_{i:}}{(B_s)_{ij}}.
\end{equation}
As in Cholesky, the residual has zeros on the pivoted rows and columns.
Let $A_r^{\mathrm{best}}$ be a best rank-at-most-$r$ approximation of
$A$ in the Frobenius norm. The known bound is
\cite[Theorem 3]{GillesWilber2026}
\begin{equation}\label{rap-eq:known-lu}
 \E\norm{B_r}_{\mathrm F}^2\le4^r\norm{A-A_r^{\mathrm{best}}}_{\mathrm F}^2.
\end{equation}

For Epperly's conjecture~\cite[Conjecture~11.2, p.~182]{Epperly2025Thesis},
we use the expected-residual formulation in the repository~\cite{TownsendRepository}.

\begin{problem}{RA-02}
Do there exist constants $C>0$ and $p\ge0$ such that the RPCholesky
residual after $r$ pivots satisfies
\[
 \E\tr R_r\le Cr^p\tau_r(A)
\]
for every Hermitian positive semidefinite $A\in\C^{n\times n}$ and
$1\le r\le n$, with $C,p$ independent of $n,r,A$?
\end{problem}

\begin{problem}{RA-03}
For every $A\in\C^{m\times n}$ and $1\le r\le\min(m,n)$, does the
RPLU residual after $r$ pivots satisfy
\[
 \E\norm{B_r}_{\mathrm F}^2
 \le 2^r\norm{A-A_r^{\mathrm{best}}}_{\mathrm F}^2
\]
as conjectured by Gilles and Wilber~\cite[discussion after Theorem~3]{GillesWilber2026}?
\end{problem}

To make both bounds sharp, fix $r\ge1$ and set $n=r+1$. For $0<t<1$, define
\begin{equation}\label{rap-eq:family}
 \begin{split}
 L(t)_{ij}&=
 \begin{cases}
  1,&i=j,\\
  t^j/(i+j),&i>j,\\
  0,&i<j,
 \end{cases}\qquad 1\le i,j\le n,\\
 \eps&=t^{n^2},\qquad
 D=\diag(1,\eps,\ldots,\eps^{n-1}),\qquad
 A_r(t)=L(t)D L(t)^T.
 \end{split}
\end{equation}
Each column of $L(t)$ contributes a rank-one term to $A_r(t)$, weighted by
the corresponding power of $\eps$ in $D$. For small $t$, the later terms
are much smaller than the first. Since $L(t)$ is invertible and the diagonal of
$D$ is positive, $A_r(t)$ is positive definite. All the columns of $L(t)$
are nonnegative and the first is strictly positive, so its contribution alone
makes every entry of $A_r(t)$ positive. The entries are rational when $t$ is rational.

\begin{theorem}[Sharp lower-bound family]\label{rap-thm:main}
For a fixed $r\ge1$, apply $r$ steps of the two algorithms to $A=A_r(t)$ in
\eqref{rap-eq:family}. Then
\begin{equation}\label{rap-eq:main-limits}
 \lim_{t\downarrow0}\frac{\E\tr R_r}{\tau_r(A)}=2^r,
 \qquad
 \lim_{t\downarrow0}\frac{\E\norm{B_r}_{\mathrm F}^2}
 {\norm{A-A_r^{\mathrm{best}}}_{\mathrm F}^2}=4^r,
\end{equation}
where $A_r^{\mathrm{best}}$ denotes a best rank-$r$ approximation of $A$.
Every square submatrix of $A_r(t)$ is nonsingular for all sufficiently small
$t>0$, with the threshold allowed to depend on $r$. Thus every ordered pivot
history appearing in the proof has positive probability.
\end{theorem}

Combining with the upper bounds in \eqref{rap-eq:known-chol} and \eqref{rap-eq:known-lu}, we obtain the following.

\begin{corollary}[Worst-case constants]\label{rap-cor:sup}
For each $r\ge1$, the supremum of the RPCholesky ratio over all
positive-semidefinite inputs with nonzero rank-$r$ tail is $2^r$. The supremum
of the RPLU squared-error ratio over all matrix inputs with nonzero rank-$r$
error is $4^r$. Restricting the lower-bound examples to real positive-definite,
entrywise-positive matrices does not change either supremum.
\end{corollary}

\begin{corollary}[Consequences for the proposed improvements]\label{rap-cor:negative}
There are no fixed constants $C>0$ and $p\ge0$ such that
$\E\tr R_r\le Cr^p\tau_r(A)$ for every $r$ and every positive-semidefinite $A$.
For RPLU, no factor $q(r)2^r$, with $q$ a fixed nonnegative polynomial,
can replace $4^r$ uniformly.
\end{corollary}
\begin{proof}
For the first statement, choose $r$ so that $2^{r-1}>Cr^p$, and then use
Theorem~\ref{rap-thm:main} with $t$ sufficiently small. For RPLU,
a zero proposed factor is immediately impossible. Otherwise,
$4^r/(2^r q(r))=2^r/q(r)$ is unbounded; again choose $r$ first and then $t$.
The rank is fixed before choosing $t$, so the fixed-rank limits in
Theorem~\ref{rap-thm:main} suffice.
\end{proof}

\subsection{Minors from the leading columns}

To estimate a minor of order $s$ of $LDL^T$, we use Cauchy--Binet and show
that the term from the first $s$ columns of $L$ is the leading one.
We let $t$ tend to zero with $n=r+1$ fixed, so constants may depend on
$n$. Write $f\sim g$ when $f/g\to1$, and use the Euclidean norm for
vector estimates. In $A[I,J]$ the indices are in increasing order
unless we specify another order.
For $S\subseteq\{1,\ldots,n\}$ of size $s$, use $1\!:\!s$ for the first
$s$ columns and set
\[
 M_s(S)=L[S,1\!:\!s],\qquad
 d_s(S)=\sum_{j\in\{1,\ldots,s\}\setminus S}j.
\]
When $s=0$, the column list is empty, and we take the determinant
of the empty matrix to be one.

\begin{lemma}[Minors in the first $s$ columns]\label{rap-lem:prefix}
For every such $S$,
\begin{equation}\label{rap-eq:prefix-asymptotic}
 \det M_s(S)=\kappa_s(S)t^{d_s(S)}(1+O(t)),
 \qquad \kappa_s(S)\ne0,
\end{equation}
where $\kappa_s(S)$ is independent of $t$. For $s\ge1$, $M_s(S)$ is invertible for sufficiently small $t$, and
$\norm{M_s(S)^{-1}}_2=O(t^{-s})$.
\end{lemma}
\begin{proof}
Let $P=S\cap\{1,\ldots,s\}$, $I=S\setminus\{1,\ldots,s\}$, and
$J=\{1,\ldots,s\}\setminus S$. Order the rows as $(P,I)$ and the columns
as $(P,J)$, where $|I|=|J|$. In the columns indexed by $P$, the diagonal
ones are the only entries that survive as $t\downarrow0$. For each
$j\in J$, the diagonal entry is absent, and dividing column $j$ by
$t^j$ leaves $1/(i+j)$ for $i>j$ and zero for $i<j$. In particular,
every row $i\in I$ has $i>s\ge j$, so the lower-right block is already
constant after scaling. The limit is
\[
 \begin{pmatrix}I&F\\0&C\end{pmatrix},
 \qquad C=\left[\frac1{i+j}\right]_{i\in I,\,j\in J},
\]
where $F$ consists of the scaled entries on rows $P$ and columns $J$.
For $J$ nonempty, list the elements of $I$ and $J$ in increasing order
as $x_a$ and $y_b$. The Cauchy determinant formula gives
\[
 \det\left[\frac1{x_a+y_b}\right]_{a,b=1}^{m}
 =\frac{\prod_{a<a'}(x_{a'}-x_a)\prod_{b<b'}(y_{b'}-y_b)}
 {\prod_{a,b}(x_a+y_b)}\ne0.
\]
The limiting matrix is also invertible when $J$ is empty, in which case
only the identity block remains. Undoing the column scaling multiplies
the determinant by $\prod_{j\in J}t^j=t^{d_s(S)}$, and restoring the
original row and column order can change its sign, giving
\eqref{rap-eq:prefix-asymptotic}. The inverse of the scaled matrix stays
bounded as $t\downarrow0$. Undoing the scaling increases its norm by
at most a factor $t^{-s}$ because each scaled column has index at most $s$.
\end{proof}

\begin{lemma}[Dominance of the leading spectral scales]\label{rap-lem:CB}
For row and column sets $I,J$ of cardinality $s\le n-1$,
\begin{equation}\label{rap-eq:CB-asymptotic}
 \det A[I,J]
 \sim \eps^{s(s-1)/2}\det M_s(I)\det M_s(J).
\end{equation}
For $s=n$, the corresponding equality is exact, with $M_n(\{1,\ldots,n\})=L$ and
$\det L=1$. Consequently, every square minor of $A$ is nonzero for all
sufficiently small positive $t$.
\end{lemma}
\begin{proof}
Cauchy--Binet gives
\[
 \det A[I,J]=\sum_{\substack{T\subseteq\{1,\ldots,n\}\\|T|=s}}
 \det L[I,T]\det L[J,T]\,
 \eps^{\sum_{j\in T}(j-1)}.
\]
The smallest power of $\eps$ occurs uniquely at $T=\{1,\ldots,s\}$ and is
$s(s-1)/2$. All other powers are at least one larger. Boundedness of $L(t)$ bounds every minor in this finite sum, so
\[
 \det A[I,J]=\eps^{s(s-1)/2}
 \left(\det M_s(I)\det M_s(J)+O(\eps)\right).
\]
Lemma~\ref{rap-lem:prefix} shows that the leading product has nonzero coefficient
and $t$-degree at most $s(s+1)$. Since $\eps=t^{n^2}$ and
$n^2>s(s+1)$ for $s\le n-1$, the remainder divided by that product tends to
zero. If $s=n$, the sum has only one term. There are finitely many minors in
a fixed dimension, so their eventual nonvanishing holds simultaneously.
\end{proof}

\subsection{The trace and Frobenius norm of a residual}

For equally sized row and column sets $I,J$, define the full, zero-padded
interpolatory residual
\[
 \mathcal R(I,J)=A-A[:,J]A[I,J]^{-1}A[I,:].
\]
For empty sets this means $\mathcal R(\varnothing,\varnothing)=A$.
The residual vanishes on the selected rows and columns. By the Schur
complement identity, successive nonzero-pivot updates give this residual.
Thus the selected sets determine the residual, although their selection
orders can have different probabilities. For Cholesky, $I=J$.

The leading term of the pivot residual comes from interpolating the next
column $l_{s+1}=L[:,s+1]$. On a set $S$ of $s<n$ rows, match this column
using the first $s$ columns of $L$. The system
$M_s(S)c_s(S)=l_{s+1}[S]$ gives the coefficients. Let $w_s(S)$ be the
interpolation error and $h_s(S)$ its squared norm.
\begin{equation}\label{rap-eq:w}
 c_s(S)=M_s(S)^{-1}l_{s+1}[S],\qquad
 w_s(S)=l_{s+1}-L[:,1\!:\!s]c_s(S),\qquad
 h_s(S)=\norm{w_s(S)}_2^2.
\end{equation}
In particular, $w_s(S)[S]=0$. For $s=0$ the fit is the zero vector, so
$w_0(\varnothing)=l_1$.

\begin{lemma}[Limiting interpolation residuals]\label{rap-lem:h}
For each fixed $S$ of size $s<n$,
\begin{equation}\label{rap-eq:h-indicator}
 \lim_{t\downarrow0}\frac1{h_s(S)}
 =\rapind_{\{s+1\notin S\}}.
\end{equation}
The limit is $1$ if row $s+1$ lies outside $S$ and $0$ if it lies in $S$.
In both cases $h_s(S)\ge1/4$ for all sufficiently small $t$.
\end{lemma}
\begin{proof}
For small $t$, the smallest singular value of $L$ is at least $1/2$ because
$L\to I$. The vector $w_s(S)$ is $L$ applied to the vector with $-c_s(S)$ in
the first $s$ slots, a $1$ in slot $s+1$, and zeros after that. This gives
\begin{equation}\label{rap-eq:h-lower}
 h_s(S)\ge\frac14\left(1+\norm{c_s(S)}_2^2\right).
\end{equation}
If $s+1\notin S$, the vector $l_{s+1}[S]$ omits the diagonal entry, and all
its entries are $O(t^{s+1})$. The inverse bound in Lemma~\ref{rap-lem:prefix}
then gives $c_s(S)=O(t)$, so $w_s(S)\to e_{s+1}$ and $h_s(S)\to1$. The same
limits hold for the empty set.

If $s+1\in S$, interpolation has to reproduce that diagonal entry, which is
one. The corresponding row of $M_s(S)$ has norm $O(t)$, so the coefficients
must satisfy $\norm{c_s(S)}_2\ge c/t$ for some $c>0$ and small $t$.
Equation~\eqref{rap-eq:h-lower} then gives $h_s(S)\to\infty$.
\end{proof}

\begin{lemma}[Asymptotic trace and squared norm]\label{rap-lem:normalisers}
Let $S$, $I$ and $J$ be subsets of $\{1,\ldots,n\}$, each of cardinality $s<n$. Then
\begin{align}
 T_s(S):=\tr\mathcal R(S,S)&\sim\eps^s h_s(S),\label{rap-eq:trace-asymptotic}\\
 Q_s(I,J):=\norm{\mathcal R(I,J)}_{\mathrm F}^2
 &\sim\eps^{2s}h_s(I)h_s(J).\label{rap-eq:norm-asymptotic}
\end{align}
In particular,
\begin{align}
 \frac{\eps^s}{T_s(S)}&\longrightarrow\rapind_{\{s+1\notin S\}},\label{rap-eq:trace-indicator}\\
 \frac{\eps^{2s}}{Q_s(I,J)}
 &\longrightarrow\rapind_{\{s+1\notin I\}}\rapind_{\{s+1\notin J\}}.
 \label{rap-eq:norm-indicator}
\end{align}
\end{lemma}
\begin{proof}
For $i\notin I$ and $j\notin J$, the bordered determinant identity, with $i$
and $j$ appended to the existing orders, gives
\[
 \mathcal R(I,J)_{ij}
 =\frac{\det A[(I,i),(J,j)]}{\det A[I,J]}.
\]
Applying Lemma~\ref{rap-lem:CB} to the numerator and denominator gives
$\eps^s$ as the ratio of their powers of $\eps$. A second bordered
determinant identity gives
\[
 \det L[(I,i),1\!:\!s+1]=\det M_s(I)\,w_s(I)_i,
\]
and the analogous column identity gives $\det M_s(J)\,w_s(J)_j$.
Consequently,
\begin{equation}\label{rap-eq:entry-asymptotic}
 \mathcal R(I,J)_{ij}
 =\eps^s w_s(I)_i w_s(J)_j(1+o(1)).
\end{equation}
The two components in this formula are nonzero for sufficiently small $t$,
by Lemma~\ref{rap-lem:prefix} applied to the bordered prefix minors. For a fixed
dimension, the relative $o(1)$ is uniform over the finitely many relevant
entries and sets.

For principal residuals, summing the leading diagonal terms
$\eps^s w_s(S)_i^2$ in \eqref{rap-eq:entry-asymptotic} gives
\eqref{rap-eq:trace-asymptotic}. For a general residual, summing the
squares of the leading entries gives \eqref{rap-eq:norm-asymptotic}.
These leading terms are nonnegative, so each sum has relative error
at most the largest relative error of its terms. This tends to zero
even when $w_s$ diverges. The reciprocal limits follow from
Lemma~\ref{rap-lem:h}, whose positive lower bound on $h_s$ prevents
amplification of the relative error.
\end{proof}

\subsection{Cancellation over pivot histories}
For a Cholesky history $(j_1,\ldots,j_r)$ of distinct indices, let
$S_s=\{j_1,\ldots,j_s\}$. For an LU history, define similarly
$I_s=\{i_1,\ldots,i_s\}$ and $J_s=\{j_1,\ldots,j_s\}$, for $0\leq s\leq r$; all three sets are empty when $s=0$.

\begin{lemma}[History identities]\label{rap-lem:paths}
For an $(r+1)\times(r+1)$ positive-definite input with all square minors
nonzero, the expectations after $r$ pivots satisfy
\begin{align}
 \E\tr R_r
 &=\det A\sum_{(j_1,\ldots,j_r)}
     \prod_{s=0}^{r-1}\frac1{T_s(S_s)},\label{rap-eq:chol-path}\\
 \E\norm{B_r}_{\mathrm F}^2
 &=|\det A|^2\sum_{\substack{(i_1,\ldots,i_r)\\(j_1,\ldots,j_r)}}
     \prod_{s=0}^{r-1}\frac1{Q_s(I_s,J_s)}.\label{rap-eq:lu-path}
\end{align}
The sums range over ordered lists of distinct indices; in the second sum
row indices and column indices are separately distinct.
\end{lemma}
\begin{proof}
In Cholesky, the product of the $r$ pivot values is the principal determinant
$\det A[S_r,S_r]$. The probability of the history is this product divided by
$\prod_{s=0}^{r-1}T_s(S_s)$. After $r$ pivots, the Schur complement formula
gives the one remaining diagonal entry of the residual as
\[
 \frac{\det A}{\det A[S_r,S_r]}.
\]
Multiplying by the history probability cancels the determinant factors.

In LU, the product of the pivot values equals the determinant of the selected
$r\times r$ block in the ordered row and column lists. Its squared modulus
is the numerator of the joint row-column history probability. The squared
modulus of the only remaining residual entry is $|\det A|^2/|\det A[I_r,J_r]|^2$,
since the row and column permutations affect only the sign of the determinant.
Multiplying by the history probability cancels $|\det A[I_r,J_r]|^2$.
Summing the weighted residuals gives \eqref{rap-eq:chol-path} and
\eqref{rap-eq:lu-path}.
\end{proof}

\subsection{The limiting contribution of each history}

Since $A$ has order $r+1$ and is positive definite,
$\tau_r(A)=\lambda_n(A)$ and
$\norm{A-A_r^{\mathrm{best}}}_{\mathrm F}^2=\lambda_n(A)^2$.
With $L^{-T}=(L^{-1})^T$, the determinant and inverse of the family satisfy
\begin{equation}\label{rap-eq:det-inverse}
 \det A=\eps^{r(r+1)/2},\qquad
 \eps^r A^{-1}
 =L^{-T}\diag(\eps^r,\eps^{r-1},\ldots,1)L^{-1}
 \longrightarrow e_n e_n^T.
\end{equation}
The largest eigenvalue of $\eps^r A^{-1}$ tends to that of
$e_ne_n^T$, namely one. Since it equals
$\eps^r/\lambda_n(A)$, we obtain
\begin{equation}\label{rap-eq:lambda-scale}
 \frac{\lambda_n(A)}{\eps^r}\longrightarrow1.
\end{equation}
Using \eqref{rap-eq:det-inverse} in the first history identity gives
\begin{equation}\label{rap-eq:normalised-chol-path}
 \frac{\E\tr R_r}{\lambda_n(A)}
 =\frac{\eps^r}{\lambda_n(A)}
 \sum_{(j_1,\ldots,j_r)}\prod_{s=0}^{r-1}
       \frac{\eps^s}{T_s(S_s)}.
\end{equation}
By Lemma~\ref{rap-lem:normalisers}, the product for a fixed history tends to $1$
precisely when
\begin{equation}\label{rap-eq:survive}
 s+1\notin S_s\qquad(0\le s\le r-1),
\end{equation}
and otherwise tends to $0$. The condition at $s=0$ is vacuous.

Under \eqref{rap-eq:survive}, label $q$ first becomes eligible at step $q$
for $2\le q\le r$. Labels $1$ and $n=r+1$ are eligible from the start.
Excluding the labels already used gives the equivalent condition
\begin{equation}\label{rap-eq:twochoices}
 j_u\in\{1,\ldots,u,n\}\setminus\{j_1,\ldots,j_{u-1}\}
 \qquad(1\le u\le r).
\end{equation}
At step $u$ the set $\{1,\ldots,u,n\}$ contains all $u-1$ earlier pivots.
It has $u+1$ elements, so there are two choices at every step. Thus
$2^r$ terms in the finite sum in \eqref{rap-eq:normalised-chol-path}
tend to one and the rest tend to zero. The prefactor tends to one by
\eqref{rap-eq:lambda-scale}.

For LU, the corresponding normalized identity is
\begin{equation}\label{rap-eq:normalised-lu-path}
 \frac{\E\norm{B_r}_{\mathrm F}^2}{\lambda_n(A)^2}
 =\frac{\eps^{2r}}{\lambda_n(A)^2}
 \sum_{\substack{(i_1,\ldots,i_r)\\(j_1,\ldots,j_r)}}
 \prod_{s=0}^{r-1}\frac{\eps^{2s}}{Q_s(I_s,J_s)}.
\end{equation}
By Lemma~\ref{rap-lem:normalisers}, a summand tends to one if both lists pass
the same test \eqref{rap-eq:survive}, and to zero otherwise. We can count
the two lists separately even though the pivot probabilities depend on both.
There are $2^r$ row lists and $2^r$ column lists, giving $4^r$ contributions
of one. Together with \eqref{rap-eq:lambda-scale}, this gives the LU limit
in Theorem~\ref{rap-thm:main}.

\subsection{Sharpness for correlation matrices with positive entries}

\begin{proposition}[Sharpness for correlation matrices]\label{rap-prop:correlation}
For each fixed $r\ge1$, both suprema in Corollary~\ref{rap-cor:sup} remain unchanged
when the inputs are real positive-definite matrices whose entries are positive and whose diagonal entries are all $1$.
\end{proposition}
\begin{proof}
First let $A$ be any rational positive-definite matrix of order $n$ with
positive entries. Choose $c>0$ rational so that
$m_i=A_{ii}/c$ is a positive integer for every $i$. Set $N=\sum_i m_i$.
Partition the $N$ row labels into groups of sizes $m_1,\ldots,m_n$, and
construct $P\in\R^{N\times n}$ by assigning each row in group $i$ the vector
$e_i^T/\sqrt{m_i}$. Then $P^TP=I_n$, and
$
 C=c^{-1}PAP^T
$
is positive semidefinite, entrywise positive, and unit-diagonal. Its nonzero
eigenvalues are those of $A/c$.

To compare pivoting on $C$ and $A$, record which group each selected index
of $C$ belongs to. If the current residual is $c^{-1}PBP^T$, summing the
Cholesky pivot probabilities over group $i$ gives
\[
 \frac{m_i B_{ii}/(cm_i)}{\tr B/c}=\frac{B_{ii}}{\tr B}.
\]
For LU, the aggregate probability of a row in group $i$ and a column in
group $j$ is
\[
 \frac{m_i m_j |B_{ij}|^2/(c^2m_i m_j)}{\norm{B}_{\mathrm F}^2/c^2}
 =\frac{|B_{ij}|^2}{\norm{B}_{\mathrm F}^2}.
\]
Substituting either pivot into its rank-one update gives the new residual
$c^{-1}P B_{\mathrm{new}}P^T$. Thus the probability of any sequence of groups
equals that of the corresponding pivot history on $A$, including those with
zero probability. The trace error and squared Frobenius error scale
by $1/c$ and $1/c^2$ respectively, as do their optimal rank-$r$ errors,
so both expected-error ratios are preserved.

Restore positive definiteness by taking
\[
 C_\delta=\frac{C+\delta I_N}{1+\delta},\qquad \delta>0.
\]
This matrix is positive definite, entrywise positive, and unit-diagonal.
Fix a finite stopping rank $r$ for which the optimal error of $C$ is positive.
Along any pivot history that occurs with positive probability for $C$,
the selected pivots are all nonzero at $\delta=0$. The rank-one updates
and conditional probabilities are continuous there, so going through
the history step by step gives convergence of both its probability and
its final residual as $\delta\downarrow0$. In either expected-error
ratio $F_r$, the denominator converges to the positive optimal error for
$C$ by continuity of the eigenvalues. To bound $F_r(C_\delta)$ from below
we can sum just over the histories that occur for $C$, since every other
history contributes a nonnegative error. Taking the limit in this sum gives
\[
 \liminf_{\delta\downarrow0} F_r(C_\delta)\ge F_r(C).
\]

Choose rational $t$ in Theorem~\ref{rap-thm:main} to approach the desired
sharp constant, apply the replication construction to $A_r(t)$, and then
take $\delta$ sufficiently small. The known universal upper bounds complete
the proof.
\end{proof}

\endgroup

\section{When concave functions spoil low-rank approximations (RA-08)}\label{rac-section}
\begingroup
Suppose we have a positive semidefinite low-rank approximation $\widehat A$
to a positive semidefinite matrix $A$ and want to approximate $f(A)$. It is cheaper to apply $f$ only to
the retained eigenvalues, giving $f(\widehat A)$ when $f(0)=0$. When does a good
approximation to $A$ give a good approximation to $f(A)$? Persson and Kressner
obtained bounds in the Frobenius, nuclear, and spectral norms for nonnegative
operator-monotone functions with $f(0)=0$ \cite{PerssonKressner2023}.

Persson, Meyer, and Musco subsequently proved transfer guarantees independent
of how the approximation is constructed~\cite[Theorems~2.4--2.7]{PerssonMeyerMusco2025}.
Their spectral-norm question for concave functions appears in
Table~1 of \cite{PerssonMeyerMusco2025}.

For $M\succeq0$, write $M_k$ for a truncation to $k$ leading
eigenvectors. Apply $f$ to eigenvalues to define $f(M)$, and use
the same $k$ eigenvectors for $M_k$ and $f(M)_k$.

\begin{problem}{RA-08}
Let $A,\widehat A$ be real symmetric matrices of order $n\ge2$ with
$0\preceq\widehat A\preceq A$, and let $1\le k<n$. For every continuous,
nonnegative, nondecreasing concave $f$ on $[0,\infty)$, every choice of
leading eigenvectors at ties, and every $\eps\ge0$, does
\[
 \norm{A-\widehat A_k}_2\le(1+\eps)\norm{A-A_k}_2
 \quad\Longrightarrow\quad
 \norm{f(A)-f(\widehat A)_k}_2\le(1+\eps)\norm{f(A)-f(A)_k}_2
\]
hold?
\end{problem}

Persson, Meyer, and Musco also found a nuclear-norm counterexample
and asked whether a fixed constant can bound amplification of the
relative excess~\cite[Example~5.1 and Section~5]{PerssonMeyerMusco2025}.
For an actual error $e$ and an optimum $e_*>0$, the relative excess
$(e-e_*)/e_*$ measures the extra error in units of the optimum.

In both examples we use $f(x)=\min(x,1)$ on $[0,\infty)$ and work with
real matrices and sketches. This continuous, nonnegative, nondecreasing
concave function fixes zero and caps the eigenvalues at one.

Given $A\in\R^{n\times n}$ with $A\succeq0$ and a full-column-rank sketch $\Omega\in\R^{n\times\ell}$, $1\le\ell\le n$, the Nystr\"om approximation is
\[
 \mathcal N(A,\Omega)=A\Omega(\Omega^TA\Omega)^{-1}\Omega^TA
\]
when the middle matrix is invertible, as it is in both examples.
In general the Nystr\"om formula uses the Moore--Penrose pseudoinverse
of $\Omega^TA\Omega$~\cite[Section~1.1]{PerssonMeyerMusco2025}. The
approximation depends only on the range of $\Omega$.

The first example starts with a Nystr\"om approximation that is best in spectral norm
for $A$. After applying $f$, its error is larger than the new optimum. For the
nuclear norm, start with an arbitrarily small relative excess over the best error.
The second construction makes the ratio of the new excess to the old one
arbitrarily large. Both examples use exact Nystr\"om approximations.

\subsection{The derivative across the threshold}
For $x,y\ge0$, clipping at one gives $|f(x)-f(y)|\le|x-y|$. If
the matrix has no eigenvalue at one, the Fr\'echet derivative of the
matrix function multiplies each entry of a perturbation in the eigenbasis
by a divided difference of $f$. These numbers are all between zero and one, but reducing
individual entries can still increase the spectral or nuclear norm of the
whole perturbation. To compute the derivative, use separate contours around
the eigenvalues below and above one \cite[Chapter~3]{HighamFunctions2008}.
\begin{lemma}\label{rac-lem:derivative}
Let $B=\diag(b_1,\ldots,b_d)\succeq0$ have no eigenvalue equal to one. For a fixed symmetric perturbation $E\in\R^{d\times d}$ and sufficiently small $t\ge0$ such that $B+tE\succeq0$,
\[
 f(B+tE)-f(B)=tD+O(t^2),\qquad
 D_{ij}=f^{[1]}(b_i,b_j)E_{ij},
\]
where $1\le i,j\le d$, the remainder is measured in the spectral norm as $t\downarrow0$, and $f^{[1]}(x,y)=(f(x)-f(y))/(x-y)$ for $x\ne y$, while $f^{[1]}(x,x)=1$ below one and zero above one.
\end{lemma}
\begin{proof}
Below one, $f$ is the identity, and above one it is constant. Since no
eigenvalue of $B$ equals one, we can extend the two pieces holomorphically
in disjoint neighborhoods of the respective eigenvalues. Take $\Gamma$
to be the union of positively oriented contours in those neighborhoods,
chosen so they still enclose the spectrum of $B+tE$ for small $t$.
The term linear in $t$ in Cauchy's matrix formula gives
\[
 D=\frac1{2\pi i}\int_\Gamma
      f(z)(zI-B)^{-1}E(zI-B)^{-1}\,\mathrm{d}z.
\]
To evaluate the $(i,j)$ entry, factor out $E_{ij}$ and integrate
$f(z)/((z-b_i)(z-b_j))$ with the same prefactor. If $b_i\ne b_j$, partial
fractions give the difference of the two values of $f$ divided by $b_i-b_j$.
If $b_i=b_j$, the double pole gives $f'(b_i)$. The contours are fixed, so the
next term of the resolvent expansion gives a remainder of $O(t^2)$ in spectral norm.
\end{proof}

\subsection{A six-dimensional spectral example}
A $3\times3$ example shows that the derivative of $f$ can increase the spectral
norm of a positive semidefinite perturbation. Two more coordinates allow us to
put that perturbation into a Nystr\"om residual. A sixth coordinate gives
the same best rank-three error before and after applying $f$.

\begin{theorem}\label{rac-thm:spectral}
Let
\[
 B=\diag\left(\frac12,\frac{127}{128},\frac{17}{16}\right),\qquad
 U=\frac19\begin{bmatrix}1&8\\8&1\\-4&4\end{bmatrix},\qquad E=UU^T.
\]
Using block sizes $3,2,1$, set
\[
 F=\frac1{65}\begin{bmatrix}
 64E&8U&0\\8U^T&I_2&0\\0&0&65
 \end{bmatrix},\qquad
 \widehat A=B\oplus0_3,\qquad
 \Omega=\begin{bmatrix}I_3\\-8U^T\\0_{1\times3}\end{bmatrix}.
\]
For $t=1/65536$ and $A=\widehat A+tF$, the target $A$ is positive definite,
\[
 \widehat A=\mathcal N(A,\Omega),\qquad
 0\preceq\widehat A\preceq A,
\]
and, with $k=3$,
\[
 \norm{A-\widehat A}_2=t=\norm{A-A_k}_2.
\]
Nevertheless,
\begin{equation}\label{rac-eq:spectralratio}
 \frac{\norm{f(A)-f(\widehat A)}_2}
      {\norm{f(A)-f(A)_k}_2}
 \ge1+\frac{334583}{15769728}>1.
\end{equation}
Thus no implication with output factor $1+C\eps$ holds for all such Nystr\"om approximations and any finite $C$, even though the input excess here is exactly zero.
\end{theorem}
\begin{proof}
The columns of $U$ satisfy $U^TU=I_2$. Thus the leading $5\times5$ block of
$F$ is the orthogonal projector onto the column space of $[8U^T\ I_2]^T/\sqrt{65}$.
Together with the projector onto the sixth coordinate, this gives $F^2=F$,
$F\succeq0$, $\rank F=3$, and $\norm{F}_2=1$.
Also $F\Omega=0$, so
\[
 A\Omega=\begin{bmatrix}B\\0_{2\times3}\\0_{1\times3}\end{bmatrix},
 \qquad \Omega^TA\Omega=B.
\]
The Nystr\"om formula therefore gives $\widehat A$.

The range of $\widehat A$ occupies the first three coordinates. The two
columns defining the leading projector cover coordinates four and five
after projection onto those coordinates, and the last projector covers
coordinate six. Thus the ranges of $\widehat A$ and $F$ together span $\R^6$.
Since both are positive semidefinite, their kernels intersect only at zero
and $A=\widehat A+tF$ is positive definite.

The rank-three approximation $\widehat A$ has residual $tF$ of
spectral norm $t$, so $\lambda_4(A)\le t$. Positive-semidefinite
order places the first three eigenvalues at or above $1/2$, while
the sixth coordinate is an eigenvector with eigenvalue $t$.
Thus at least four eigenvalues are at or above $t$, and
$\lambda_4(A)=t$. Since $f$ preserves the order and $t<1$, the
fourth transformed eigenvalue, and hence the optimal transformed
spectral error, is again $t$.

Let $\mathcal D$ denote the derivative of $f$ at $\widehat A$ in direction $F$. By Lemma~\ref{rac-lem:derivative}, its upper-left $3\times3$ block is $(64/65)D$, where
\[
 D=\frac1{729}\begin{bmatrix}
 585&144&224\\144&585&-28\\224&-28&0
 \end{bmatrix}.
\]
For $v=(4,3,1)^T/\sqrt{26}$ and $\widetilde v=v\oplus0_3$, exact arithmetic gives
\begin{equation}\label{rac-eq:witness}
 v^TDv=\frac{19705}{18954},\qquad
 \widetilde v^T\mathcal D\widetilde v
   =\frac{64}{65}\frac{19705}{18954}.
\end{equation}
The second Rayleigh quotient in \eqref{rac-eq:witness} exceeds one,
so it suffices to show that the remainder after division by the stated
$t$ is smaller than that excess.

Let $a=17/16$, $d=a-1=1/16$, and $r=9/256$. The circle $\Gamma$ centered at $a$ with radius $r$ isolates the unique eigenvalue above $1$ for both $\widehat A$ and $A$. The gap between the two largest eigenvalues of $\widehat A$ is $2r$, so every point of the circle is at least distance $r$ from its spectrum. The other eigenvalues of $A$ stay below one because $127/128+t<1$.

Set $R(z)=(zI-\widehat A)^{-1}$ and $R_t(z)=(zI-\widehat A-tF)^{-1}$. Since $t<r/2$,
\[
 \norm{R(z)}_2\le r^{-1},\qquad
 \norm{R_t(z)}_2\le2r^{-1}\quad(z\in\Gamma).
\]
Substitute the identity $R_t=R+tRFR+t^2RFRFR_t$ into
\[
 f(A)=A-\frac1{2\pi i}\int_\Gamma(z-1)R_t(z)\,\mathrm{d}z.
\]
The term of order $t$ is $t\mathcal D$. The remaining integral contains
$t^2RFRFR_t$, of norm at most $2t^2/r^3$ since $\norm{F}_2=1$. Multiplying
by $|z-1|\le d+r$ and integrating over the circle of length $2\pi r$ gives
\begin{equation}\label{rac-eq:remainder}
 \norm{f(A)-f(\widehat A)-t\mathcal D}_2
 \le\frac{2(d+r)}{r^2}t^2=\frac{12800}{81}t^2.
\end{equation}
Combining \eqref{rac-eq:witness} and \eqref{rac-eq:remainder},
\[
 \frac{\widetilde v^T(f(A)-f(\widehat A))\widetilde v}{t}
 \ge\frac{64}{65}\frac{19705}{18954}
       -\frac{12800}{81\cdot65536}
 =1+\frac{334583}{15769728}.
\]
The left-hand side is at most $\norm{f(A)-f(\widehat A)}_2/t$. Since $t$
is the best rank-three error for $f(A)$, this proves \eqref{rac-eq:spectralratio}.
\end{proof}

\subsection{Nuclear error and increasing dimension}
\begin{theorem}\label{rac-thm:nuclear}
For every finite $C\ge0$ and every prescribed $\eps_0>0$, there exist a positive-definite matrix $A$, a full-column-rank sketch $\Omega$, and $k=\rank\mathcal N(A,\Omega)$ such that, with $\widehat A=\mathcal N(A,\Omega)$,
\[
 \norm{A-\widehat A}_*=(1+\eps)\norm{A-A_k}_*,
 \qquad 0<\eps<\eps_0,
\]
while
\[
 \norm{f(A)-f(\widehat A)}_*>(1+C\eps)\norm{f(A)-f(A)_k}_*.
\]
The function $f(x)=\min(x,1)$ is fixed independently of $C$ and $\eps_0$.
\end{theorem}

Fix an integer $m\ge1$, set $\rho=1/4$ and $h_i=\rho^i$ for $1\le i\le m$, and let $\one_q$ denote the vector of $q$ ones. Define
\[
 B_m=\diag(1-h_1,\ldots,1-h_m,1+h_1,\ldots,1+h_m),
 \quad u=\frac{\one_{2m}}{\sqrt{2m}},\quad E=uu^T.
\]
Let $J_m$ be the $m\times m$ matrix of ones. By Lemma~\ref{rac-lem:derivative},
the derivative at $B_m$ in direction $E$ is
\begin{equation}\label{rac-eq:Dm}
 D_m=\frac1{2m}\begin{bmatrix}J_m&K_m\\K_m^T&0\end{bmatrix},
 \qquad (K_m)_{ij}=\frac{h_i}{h_i+h_j}
                  =\frac1{1+\rho^{j-i}},
 \quad 1\le i,j\le m.
\end{equation}
As we move away from the diagonal, the entries of $K_m$ approach one above
it and zero below it, with geometrically decaying errors. This makes it
useful to compare $K_m$ with the triangular matrix of ones on and above the
diagonal, whose singular values we can calculate.

\begin{lemma}\label{rac-lem:divergence}
The derivative nuclear norms grow without bound as $m$ increases, with
\[
 \norm{D_m}_*\ge\frac1\pi\log(m+1)-\frac12-\sqrt{\frac2{15}}.
\]
\end{lemma}
\begin{proof}
Removing the diagonal blocks of $D_m$ contracts the nuclear norm: with $S=\diag(I_m,-I_m)$, the triangle inequality and unitary invariance give
\[
 \norm{D_m}_*\ge\norm{\frac{D_m-SD_mS}{2}}_*
 =\frac{\norm{K_m}_*}{m}.
\]
Let $T_m$ have ones on and above its diagonal and zeros below. Entrywise geometric decay gives
\begin{align*}
 \norm{K_m-(T_m-I_m/2)}_{\mathrm F}^2
 =2\sum_{d=1}^{m-1}(m-d)
          \left(\frac{\rho^d}{1+\rho^d}\right)^2
 \le\frac{2m\rho^2}{1-\rho^2}.
\end{align*}
Using $\norm{X}_*\le\sqrt m\norm{X}_{\mathrm F}$ for an $m\times m$ matrix,
\begin{equation}\label{rac-eq:Tbound}
 \frac{\norm{K_m}_*}{m}
 \ge\frac{\norm{T_m}_*}{m}-\frac12
        -\sqrt{\frac{2\rho^2}{1-\rho^2}}.
\end{equation}
To find the singular values of $T_m$, start with $(T_m^T)^{-1}$, the lower
bidiagonal first-difference matrix with diagonal one and subdiagonal minus
one. Its Gram matrix has diagonal $2,\ldots,2,1$ and adjacent off-diagonal
entries minus one. The interior eigenvector recurrence is solved by
$x_i=\sin(i\theta)$. The first and last rows impose $x_0=0$ and
$x_{m+1}=x_m$, respectively, so $\theta=(2j-1)\pi/(2m+1)$. Taking reciprocal
square roots of the Gram eigenvalues $4\sin^2(\theta/2)$ gives
\begin{equation}\label{rac-eq:Tsingular}
 \sigma_j(T_m)=\frac1{2\sin((2j-1)\pi/(4m+2))},
 \qquad 1\le j\le m.
\end{equation}

Since $\sin x\le x$ for $x\ge0$,
\[
 \frac{\norm{T_m}_*}{m}
 \ge\frac{2m+1}{m\pi}\sum_{j=1}^m\frac1{2j-1}
 \ge\frac1\pi\sum_{j=1}^m\frac1j
 \ge\frac1\pi\log(m+1).
\]
Substituting $\rho=1/4$ into \eqref{rac-eq:Tbound} proves the lemma.
\end{proof}

\begin{proof}[Proof of Theorem~\ref{rac-thm:nuclear}]
To embed $D_m$ in a Nystr\"om residual, take $0<\eta<1$, put $d=2m$, and set
\[
 v=\begin{bmatrix}\sqrt\eta\,u\\\sqrt{1-\eta}\end{bmatrix},\qquad
 \widehat A=B_m\oplus[0],\qquad
 A_t=\widehat A+tvv^T,
 \qquad
 \Omega=\begin{bmatrix}I_d\\-\sqrt{\eta/(1-\eta)}\,u^T\end{bmatrix}.
\]
Because $v^T\Omega=0$, we have
\[
 A_t\Omega=\begin{bmatrix}B_m\\0\end{bmatrix},\qquad
 \Omega^TA_t\Omega=B_m,
\]
so $\mathcal N(A_t,\Omega)=\widehat A$. The positive-semidefinite summands have ranges spanning $\R^{d+1}$, so $A_t$ is positive definite for every $t>0$.

Let $a(t)=\lambda_{d+1}(A_t)$. For $0<t<1/4$, the other $d$ eigenvalues are at least $3/4$, while $0<a(t)\le t$. The block Schur complement gives
\begin{align}\label{rac-eq:smalltail}
 a(t)={}&(1-\eta)t
 -t^2\eta(1-\eta)
 u^T(B_m+t\eta uu^T-a(t)I_d)^{-1}u.\notag
\end{align}
The matrix inside the inverse is $\succeq\tfrac12 I_d$, so the inverse has spectral norm at most two and $a(t)=(1-\eta)t+O(t^2)$.
At approximation rank $d$, the optimal nuclear tail is this one eigenvalue, which $f$ leaves unchanged. Both optimal tails are therefore $a(t)$, while the residual $tvv^T$ has nuclear norm $t$.

By Lemma~\ref{rac-lem:derivative}, the derivative of $f(A_t)-f(\widehat A)$ at zero has upper-left block $\eta D_m$ and bottom-right entry $1-\eta$. Pinching to these two diagonal blocks is contractive in nuclear norm: it equals the average of a matrix and its conjugate by $\diag(I_d,-1)$. Therefore
\begin{equation}\label{rac-eq:pinchedlimit}
 \liminf_{t\downarrow0}
 \frac{\norm{f(A_t)-f(\widehat A)}_*}{t}
 \ge\eta\norm{D_m}_*+1-\eta.
\end{equation}
Measure the relative excess above the best error before and after applying $f$ by
\[
 \eps_A(t)=\frac{t}{a(t)}-1,\qquad
 \eps_f(t)=\frac{\norm{f(A_t)-f(\widehat A)}_*}{a(t)}-1.
\]
Using $a(t)=(1-\eta)t+O(t^2)$ and \eqref{rac-eq:pinchedlimit} gives
\[
 \eps_A(t)\longrightarrow\frac\eta{1-\eta}>0,
 \qquad
 \liminf_{t\downarrow0}\frac{\eps_f(t)}{\eps_A(t)}
 \ge\norm{D_m}_*.
\]
Choose $\eta$ so the first limit lies in $(0,\eps_0)$. With $\eta$ fixed,
Lemma~\ref{rac-lem:divergence} lets us take $m$ large enough that
$\norm{D_m}_*>C$. For all sufficiently small $t>0$ we then have
$0<\eps_A(t)<\eps_0$ and $\eps_f(t)/\eps_A(t)>C$.
\end{proof}

\endgroup

\section{Sharp Frobenius transfer beyond operator monotonicity (RA-09)}\label{raf-section}
\begingroup
\newcommand{\rafTA}{T_A}
\newcommand{\rafTf}{T_f}
\newcommand{\rafrange}{\operatorname{range}}
Suppose we already have a rank-$k$ approximation $B$ to $A$, with
both matrices positive semidefinite. We form an approximation $C$
to $f(A)$ by keeping the $k$ eigenvectors used to represent $B$
and applying $f$ to their eigenvalues. We set $C$ to zero on the
remaining directions, where $f(B)$ has eigenvalue $f(0)$. This
keeps the rank at most $k$ and gives $C=f(B)$ when $f(0)=0$.

For a positive semidefinite matrix $M$, let $M_k$ be a truncation to
$k$ leading eigenvectors, including zeros if needed. Apply $f$ to the
eigenvalues by spectral calculus, and form $f(M)_k$ on the same
eigenvectors as $M_k$~\cite[Eq.~(2)]{PerssonMeyerMusco2025}.

\begin{problem}{RA-09}
Let $A,\widehat A$ be real symmetric matrices of order $n\ge2$ with
$0\preceq\widehat A\preceq A$, and let $1\le k<n$. For every continuous,
nonnegative, nondecreasing concave $f$ on $[0,\infty)$, every choice of
leading eigenvectors at ties, and every $\eps\ge0$, does
\begin{equation}\label{raf-eq:canonical-premise}
 \begin{gathered}
 \norm{A}_{\mathrm F}^2-\norm{\widehat A_k}_{\mathrm F}^2
 \le(1+\eps)\norm{A-A_k}_{\mathrm F}^2\\
 \Longrightarrow\quad
 \norm{f(A)-f(\widehat A)_k}_{\mathrm F}^2
 \le(1+\eps)\norm{f(A)-f(A)_k}_{\mathrm F}^2
 \end{gathered}
\end{equation}
hold?
\end{problem}

Using the Frobenius residual $\norm{A-B}_{\mathrm F}$ directly, we show
that an input factor $1+\eps$ in Frobenius norm gives an output factor
of at most $1+\eps$ when $B\preceq A$ and $1+2\eps$ for general
$B\succeq0$, with sharp coefficients of $\eps$ in both cases.

Applying $f$ to the retained eigenvalues of a Nystr\"om approximation
is the idea behind funNystr\"om. Persson and Kressner obtained
probabilistic Frobenius bounds using randomized subspace iteration
for nonnegative operator-monotone $f$ with $f(0)=0$
\cite[Section~3.2]{PerssonKressner2023}.
They also used this approximation to estimate $\tr f(A)$ by computing
$\tr C$ and estimating the trace of the residual $R=f(A)-C$
with Gaussian probes~\cite[Section~5.2]{PerssonKressner2023}.
In an eigenbasis of $R$, a probe $g^TRg$ with $g\sim N(0,I_n)$ is a sum
of independent squared standard Gaussians weighted by the eigenvalues
of $R$. Its mean is $\tr R$ and its variance is $2\norm{R}_{\mathrm F}^2$,
so reducing the Frobenius residual makes the trace estimate less variable.

\subsection{Approximation on the selected eigenspace}
We prove the results for real symmetric matrices; for complex Hermitian matrices, replace transposes by adjoints. For scalar $x$, let $x_+=\max(x,0)$. Throughout, $f$ is continuous, nonnegative and nondecreasing on $[0,\infty)$, and
\begin{equation}\label{raf-eq:subhom}
 x\longmapsto f(x)/x\quad\text{is nonincreasing on }(0,\infty).
\end{equation}
For $0<x<y$, concavity and nonnegativity at zero give
$f(x)\ge (x/y)f(y)+(1-x/y)f(0)\ge(x/y)f(y)$, so
$f(x)/x\ge f(y)/y$ even when $f(0)>0$.

Let $A,B\in\R^{n\times n}$ be positive semidefinite, with
$a_1\ge\cdots\ge a_n\ge0$ the eigenvalues of $A$. Fix $1\le k<n$
and suppose $\rank B\le k$. Choose $k$ orthonormal eigenvectors
$v_j$ of $B$ whose span contains its range, and denote their
eigenvalues by $b_j$. A choice of span matters when $f(0)>0$.
Apply $f$ on this span and extend the result by zero, giving
\[
 B=\sum_{j=1}^k b_jv_jv_j^T,\qquad
 C=\sum_{j=1}^k f(b_j)v_jv_j^T,
\]
The squared Frobenius errors of the best rank-$k$ approximations to
$A$ and $f(A)$ are the sums of squares of the discarded eigenvalues,
\[
 \rafTA=\sum_{i>k}a_i^2,\qquad
 \rafTf=\sum_{i>k}f(a_i)^2.
\]
Let $\tau=a_{k+1}$ be the cutoff for $A$ and, when $\tau>0$,
set $c=f(\tau)/\tau$ to be the slope of the line from the origin
to $(\tau,f(\tau))$.
Since $f(x)/x$ is nonincreasing, this line lies below $f$ up to $\tau$
and above $f$ thereafter.

For ordered approximations $B\preceq A$, Persson, Meyer, and
Musco use $\norm{A}_{\mathrm F}^2-\norm{B}_{\mathrm F}^2$ to
measure the input error in the operator-monotone
case~\cite[Theorem~2.5]{PerssonMeyerMusco2025}. Since this quantity
bounds $\norm{A-B}_{\mathrm F}^2$ from above, a transfer bound using
the squared residual settles the concave case in their Table~1.

\begin{theorem}[Ordered transfer]\label{raf-thm:ordered}
Suppose $0\preceq B\preceq A$ and $\tau>0$. Then
\begin{equation}\label{raf-eq:ordered}
 \norm{f(A)-C}_{\mathrm F}^2-\rafTf
 \le c^2\bigl(\norm{A-B}_{\mathrm F}^2-\rafTA\bigr).
\end{equation}
Consequently, for every $\eta\ge0$,
\[
 \norm{A-B}_{\mathrm F}^2\le(1+\eta)\rafTA
 \quad\Longrightarrow\quad
 \norm{f(A)-C}_{\mathrm F}^2\le(1+\eta)\rafTf.
\]
For $\eps\ge0$, set $\eta=(1+\eps)^2-1$ to obtain the factor $1+\eps$
in Frobenius norm. The multiplier of one on the excess in
\eqref{raf-eq:ordered} is sharp.
\end{theorem}

Example~5.3 of \cite{PerssonMeyerMusco2025} shows that the relative
excess can increase if $B$ is allowed to exceed $A$ in some direction.
Section~5 of that paper asks whether there is a uniform factor for
this increase.

\begin{theorem}[General Frobenius transfer]\label{raf-thm:unordered}
Suppose $A,B\succeq0$ and $\rank B\le k$. If $\tau>0$, then
\begin{equation}\label{raf-eq:unordered}
 \norm{f(A)-C}_{\mathrm F}^2-\rafTf
 \le 2c^2\bigl(\norm{A-B}_{\mathrm F}^2-\rafTA\bigr).
\end{equation}
Thus, a squared relative excess $\eta\ge0$ becomes at most $2\eta$. In norm-relative notation, for every $\eps\ge0$,
\[
 \norm{A-B}_{\mathrm F}\le(1+\eps)\sqrt{\rafTA}
 \quad\Longrightarrow\quad
 \norm{f(A)-C}_{\mathrm F}\le(1+2\eps)\sqrt{\rafTf}.
\]
The universal constant two is optimal in both conventions, even when $f$ is restricted to operator-monotone power functions.
\end{theorem}

\subsection{Reductions for the proofs}
For Problem~\ref{prob:RA-09}, set $B=\widehat A_k$. Then
$0\preceq B\preceq\widehat A\preceq A$ for every permitted choice
within a repeated eigenspace. Expanding the residual gives
$$
 \norm{A-B}_{\mathrm F}^2
 =\norm{A}_{\mathrm F}^2-\norm{B}_{\mathrm F}^2-2\tr\bigl(B(A-B)\bigr)
 \le\norm{A}_{\mathrm F}^2-\norm{B}_{\mathrm F}^2,
$$
since $\tr(B(A-B))=\tr(B^{1/2}(A-B)B^{1/2})\ge0$. So the premise of
\eqref{raf-eq:canonical-premise} bounds $\norm{A-B}_{\mathrm F}^2$ by
$(1+\eps)\rafTA$. For $\tau>0$ the ordered theorem gives the conclusion
with that same $\eps$. The case $\rafTA=0$ is settled at the end.

For $0<x<y$, rearrange $f(x)/x\ge f(y)/y$ to bound the secant slope
by $f(y)/y$. Continuity gives the same bound at $x=0$. If $y\ge\tau>0$,
then $f(y)/y\le c$, so
\begin{equation}\label{raf-eq:lip}
 0\le f(y)-f(x)\le \frac{f(y)}y(y-x)\le c(y-x).
\end{equation}
Applying $f(x)\ge cx$ to each discarded eigenvalue gives
\begin{equation}\label{raf-eq:tail}
 c^2\rafTA\le\rafTf.
\end{equation}
If $f(\tau)=0$, then $c=0$. For $x\le\tau$ monotonicity gives
$f(x)\le f(\tau)=0$, and for $x\ge\tau$ we have $f(x)\le cx=0$.
Nonnegativity now gives $f(A)=C=0$, so the error bound holds with equality.
We can therefore take $c>0$.

\subsection{Proof of ordered transfer}
Choose an orthonormal eigenbasis $u_1,\ldots,u_n$ of $A$. In the
mixed trace $\tr(AB)$, each product $a_i b_j$ carries weight
$|u_i^Tv_j|^2$. We use similar weights for the transformed error.
Since $0\preceq B\preceq A$, we have $A\succeq B\succeq b_jv_jv_j^T$
for each $b_j>0$. Lemma~\ref{raf-lem:harmonic}
shows that $\sum_{a_i>0}|u_i^Tv_j|^2(1-b_j/a_i)\ge0$.
We can therefore use a scalar bound containing $1-b/a$, whose sign
only needs to be controlled after averaging.

\begin{lemma}[Ordered scalar certificate]\label{raf-lem:orderedscalar}
Set $h(a)=(c^2a^2-f(a)^2)_+$ for $a\ge0$. For $a,b>0$,
\begin{align}\label{raf-eq:scalarordered}
 h(a)+2f(b)f(a)-2c^2ba-f(b)^2+c^2b^2
 \ge f(b)(f(b)-cb)_+\left(1-\frac ba\right).
\end{align}
\end{lemma}
\begin{proof}
If $b\ge\tau$, the right side vanishes. The left side is
\[
 c^2(a-b)^2-(f(a)-f(b))^2+(f(a)^2-c^2a^2)_+,
\]
which is nonnegative by \eqref{raf-eq:lip}, since $\max(a,b)\ge\tau$.

For $b\le\tau$, the factor $f(b)-cb$ on the right of
\eqref{raf-eq:scalarordered} can be positive. Scale the arguments and
function values by $b$ and $cb$ to write
\[
 d=\frac{f(b)}{cb}\ge1,\qquad z=\frac ab>0,\qquad
 G=\frac{2f(b)f(a)+h(a)}{c^2b^2}.
\]
For $0<z\le1$, the inequality $f(a)/a\ge f(b)/b$ gives
$f(a)\ge zf(b)$. Since $h(a)\ge0$, this gives $G\ge2d^2z$.
For $1\le z\le d$, monotonicity gives $f(a)\ge f(b)$ and hence
$G\ge2d^2$.

For $z\ge d$, we seek the larger bound $G\ge2dz$. If $a\le\tau$,
this follows from $f(a)\ge ca$ and $h(a)\ge0$. If $a\ge\tau$, set
$\phi=f(a)/(cb)\in[d,z]$. On this interval $G$ is a concave quadratic
in $\phi$, so its minimum occurs at $\phi=d$ or $\phi=z$. Thus
\[
 G=z^2+2d\phi-\phi^2\ge\min\{z^2+d^2,2dz\}=2dz.
\]
Collecting the three bounds gives
\begin{equation}\label{raf-eq:Gbranches}
 G\ge\begin{cases}
 2d^2z,&0<z\le1,\\
 2d^2,&1\le z\le d,\\
 2dz,&z\ge d.
 \end{cases}
\end{equation}

Subtract the right side of \eqref{raf-eq:scalarordered} from the left
and divide by $c^2b^2$. Using \eqref{raf-eq:Gbranches}, we can bound
this difference from below by
\begin{equation}\label{raf-eq:branchcertificates}
 \begin{cases}
 \displaystyle\frac{d-1}{z}
  \bigl(2(d+1)z^2-(2d+1)z+d\bigr),&0<z\le1,\\[5pt]
 \displaystyle\frac{(d-z)(2z+d-1)}{z},&1\le z\le d,\\[5pt]
 \displaystyle\frac{(d-1)(z-d)(2z-1)}{z},&z\ge d.
 \end{cases}
\end{equation}
In the first branch, the quadratic in $z$ is positive because its leading
coefficient is positive and its discriminant is $1-4d(d+1)<0$. Its
prefactor is nonnegative. The other two branches are nonnegative by
their respective bounds on $z$. This proves \eqref{raf-eq:scalarordered}.
\end{proof}

\begin{lemma}[Rank-one bound in the L\"owner order]\label{raf-lem:harmonic}
If $A\succeq bvv^T$ with $b>0$ and $\norm{v}_2=1$, then $v\in\rafrange A$ and
$
 b\,v^TA^+v\le1,
$
where $A^+$ is the Moore--Penrose inverse.
\end{lemma}
\begin{proof}
For $x\in\ker A$, the inequality gives $b|x^Tv|^2\le0$.
Thus $v$ lies in the range of $A$, and $x=A^+v$ solves $Ax=v$.
On this vector the inequality reads $v^Tx\ge b(v^Tx)^2$.
Since $v^Tx>0$, it gives $b\,v^TA^+v\le1$.
\end{proof}

\begin{proof}[Proof of Theorem~\ref{raf-thm:ordered}]
Fix $b_j>0$ and write $v_j$ in an orthonormal eigenbasis $u_i$ of $A$.
Since $A\succeq B\succeq b_jv_jv_j^T$, Lemma~\ref{raf-lem:harmonic}
puts $v_j$ in the range of $A$. Thus its weights $|u_i^Tv_j|^2$ sum
to one over $a_i>0$. Use them to average
Lemma~\ref{raf-lem:orderedscalar}. The factor $1-b_j/a_i$ averages to
$1-b_jv_j^TA^+v_j\ge0$ by the rank-one bound. Its multiplier depends
only on $b_j$ and is nonnegative. The averaged left side is therefore
nonnegative, giving
\begin{align}\label{raf-eq:average}
 f(b_j)^2-c^2b_j^2-2f(b_j)v_j^Tf(A)v_j
       +2c^2b_jv_j^TAv_j
 \le v_j^Th(A)v_j.
\end{align}
If $b_j=0$, then $f(A)\succeq f(0)I$ makes the left side at most
$-f(0)^2$. The right side is nonnegative, so
\eqref{raf-eq:average} holds in this case too.

Sum \eqref{raf-eq:average} over $j$. Since $\sum_jv_jv_j^T\preceq I$
and $h(A)\succeq0$, the right sides sum to at most $\tr h(A)$.
The left side supplies the mixed terms when we expand the two errors.
In the resulting trace bound, $h(a_i)$ leaves only the positive part of
$f(a_i)^2-c^2a_i^2$. This part is zero for $i\le k$ because $a_i\ge\tau$,
and equals $f(a_i)^2-c^2a_i^2$ for $i>k$ because $a_i\le\tau$. Thus
\begin{align*}
 \norm{f(A)-C}_{\mathrm F}^2-c^2\norm{A-B}_{\mathrm F}^2
 &\le\tr\bigl(f(A)^2-c^2A^2\bigr)+\tr h(A)\\
 &=\sum_{i=1}^n\bigl(f(a_i)^2-c^2a_i^2\bigr)_+=\rafTf-c^2\rafTA.
\end{align*}
Rearranging proves \eqref{raf-eq:ordered}, and \eqref{raf-eq:tail}
gives the relative-error statement.

The constant is sharp already for $f(x)=x$, when the two excesses are equal. For example, $A=\diag(2,1)$, $B=\diag(1,0)$ and $k=1$ have positive residual excess, so the multiplier cannot be smaller than one.
\end{proof}

\subsection{Proof of transfer without an order assumption}
For arbitrary $B\succeq0$, we bound the contribution of a pair
$a_i,b_j$ by a function of $a_i$. We need that function to be
nondecreasing. Its trace on any $k$-dimensional subspace is then at most
its trace on the leading eigenspace of $A$. The scalar inequality below
gives such a bound once we rescale the cutoff from one to $\tau$.

\begin{lemma}[Unordered scalar certificate]\label{raf-lem:unorderedscalar}
Let $g$ satisfy the same assumptions as $f$, with $g(1)=1$. For all $a,b\ge0$,
\begin{equation}\label{raf-eq:scalarunordered}
 g(b)^2-2g(b)g(a)-2b^2+4ab
 \le H(a):=\max\{2a^2-g(a)^2,1\}.
\end{equation}
\end{lemma}
\begin{proof}
Suppose $\max(a,b)\ge1$. Subtracting $2a^2-g(a)^2$ from the left
side of \eqref{raf-eq:scalarunordered} gives
\[
 (g(a)-g(b))^2-2(a-b)^2\le0
\]
by \eqref{raf-eq:lip} with $\tau=c=1$. Since $H(a)\ge2a^2-g(a)^2$,
this proves the bound whenever one of the arguments reaches one.

For $a,b\in[0,1]$, put $u=g(a)$ and $v=g(b)$. The ratio bound and
monotonicity give $a\le u\le1$ and $b\le v\le1$. When $a\le b$ and
$b>0$, write $z=a/b$. Then $u\ge zv$, and the left side of
\eqref{raf-eq:scalarunordered} is at most $(1-2z)(v^2-2b^2)$.
Since $b^2\le v^2\le1$ and $0\le b\le1$, we have
$-1\le-b^2\le v^2-2b^2\le1$. The other factor satisfies
$|1-2z|\le1$ for $z\in[0,1]$. Hence the product is at most one.

If $b\le a$, then $u\ge\max(a,v)$. For $v\ge a$, the expression is at most
\[
 -v^2-2b^2+4ab\le-a^2-2b^2+4ab
 =a^2-2(a-b)^2\le1.
\]
When $v\le a$, the left side of \eqref{raf-eq:scalarunordered} is at
most $v^2-2av-2b^2+4ab$. This is convex in $v$, so it suffices to
check the ends of the interval $[b,a]$. At $v=b$ and $v=a$, respectively,
the values are
\[
 2ab-b^2=a^2-(a-b)^2,\qquad
 -a^2-2b^2+4ab=a^2-2(a-b)^2,
\]
both at most one. At $a=b=0$, the left side of
\eqref{raf-eq:scalarunordered} is $-g(0)^2\le1$.
\end{proof}

\begin{proof}[Proof of Theorem~\ref{raf-thm:unordered}]
Lemma~\ref{raf-lem:unorderedscalar} puts the cutoff at one. To use it
at $\tau$, take $g(x)=f(\tau x)/f(\tau)$ and apply the lemma to
$a/\tau$ and $b/\tau$. After multiplying by $f(\tau)^2$, the
right side is
\[
 H_\tau(a)=\max\{2c^2a^2-f(a)^2,f(\tau)^2\}.
\]
We will need $H_\tau$ to be nondecreasing when we take its trace on a
subspace. Up to $\tau$, we have $f(a)\ge ca$, so
$2c^2a^2-f(a)^2\le c^2a^2\le f(\tau)^2$ and $H_\tau$ is constant.
For $y\ge x\ge\tau$, multiply \eqref{raf-eq:lip} by
$f(y)+f(x)\le c(y+x)$ to obtain
\[
 f(y)^2-f(x)^2\le c^2(y^2-x^2).
\]
Thus $2c^2a^2-f(a)^2$ is nondecreasing from its value $f(\tau)^2$
at $a=\tau$.

For each $v_j$, use the squared coordinates $|u_i^Tv_j|^2$ as weights
in the scalar bound. They sum to one, and the weighted right side is
$v_j^TH_\tau(A)v_j$. Adding these over $j$ gives the trace of $H_\tau(A)$
on the subspace with projector $P=\sum_jv_jv_j^T$. The Ky Fan maximum
principle bounds this trace by the sum of the largest $k$ eigenvalues
of $H_\tau(A)$. By monotonicity
these are $H_\tau(a_1),\ldots,H_\tau(a_k)$, so
\begin{align*}
 &\sum_{j=1}^k\bigl(f(b_j)^2-2c^2b_j^2
       -2f(b_j)v_j^Tf(A)v_j+4c^2b_jv_j^TAv_j\bigr)\\
 &\hspace{20mm}\le\tr(PH_\tau(A))
 \le\sum_{i=1}^k\bigl(2c^2a_i^2-f(a_i)^2\bigr).
\end{align*}
Add $\tr(f(A)^2-2c^2A^2)$ to the two outer sides. The left becomes
the difference of the squared errors, and the right becomes
$\rafTf-2c^2\rafTA$, giving \eqref{raf-eq:unordered}. To compare relative excesses, divide
by the optimal tails and use $c^2\rafTA\le\rafTf$ from \eqref{raf-eq:tail}.
An input norm factor $1+\eps$ gives squared relative excess at most
$(1+\eps)^2-1$. The output norm factor then follows from
$\sqrt{1+2((1+\eps)^2-1)}\le1+2\eps$ for $\eps\ge0$.
\end{proof}

\subsection{Sharpness and degenerate cases}
The factor two comes from putting the single nonzero eigenvalue of $B$
in a direction where $A$ is zero. Fix an integer $N\ge1$, $0<b<1$ and
$0<p<1$. With $e_j$ the $j$th coordinate vector in $\R^{N+2}$, take
\[
 A=I_{N+1}\oplus[0],\qquad k=1,\qquad
 B=b\,e_{N+2}e_{N+2}^T,\qquad f(x)=x^p.
\]
Then $\tau=c=1$, $\rafTA=\rafTf=N$, and
\[
 \norm{A-B}_{\mathrm F}^2=N+1+b^2,\qquad
 \norm{f(A)-C}_{\mathrm F}^2=N+1+b^{2p}.
\]
Here $f(A)=A$ and $C=b^pe_{N+2}e_{N+2}^T$. Subtracting the common
optimal tail $N$ leaves the ratio
$(1+b^{2p})/(1+b^2)$. Take $N=L$, $b=e^{-L}$ and $p=L^{-2}$ for
integers $L\ge2$. As $L$ tends to infinity, $b^2$ tends to zero and
$b^{2p}=e^{-2/L}$ tends to one, so the ratio tends to two. As $N$
increases, both squared errors divided by their respective optima
approach one. Taking the square root of $1+x$ multiplies the excess
$x$ by $1/(\sqrt{1+x}+1)$. This factor tends to $1/2$ for both errors,
so the ratio of the two norm excesses also tends to two.

The power functions in this example preserve semidefinite order.
For $0<p<1$, apply spectral calculus to the integral
\[
 x^p=\frac{\sin(\pi p)}\pi
       \int_0^\infty t^{p-1}\,x(x+t)^{-1}\,\mathrm{d}t
\]
to evaluate $X^p$ for a positive semidefinite matrix $X$. For each $t>0$,
the matrix $X(X+tI)^{-1}=I-t(X+tI)^{-1}$ is increasing in $X$ because
inversion reverses the order of $X+tI$. Integrating with positive
weights preserves the order.

If $\tau=0$, the optimal tail for $A$ is zero, and the relative-error
hypothesis gives $B=A$. The chosen $k$ directions contain
$\rafrange A$, so $C$ agrees with $f(A)$ on their span. On the
orthogonal complement the error is $f(0)I$, exactly the optimal tail
on those $n-k$ directions. For $k=n$ the complement is empty.
\endgroup

\section{Sharp constants for adversarial row deletion (IE-21, IE-22)}\label{sec:rows}
Consider a consistent linear system $Ax=b$ with unit-norm rows.
Retaining the equations indexed by $S$ gives the system $A_Sx=b_S$.
If $A_S$ has full column rank, the least-squares estimate from
perturbed data $b_S+e$ is $\widehat x=x+A_S^\dagger e$, where
$A_S^\dagger$ is the pseudoinverse. The sharp stability bound is
\[
 \norm{\widehat x-x}_2\le\frac{\norm{e}_2}{\sigma_{\min}(A_S)}.
\]
A nontrivial kernel of $A_S$ gives multiple solutions to the
retained system. For a prescribed number $k$ of retained equations,
how should we choose the rows of $A$ in advance to maximize
$\min_{|S|=k}\sigma_{\min}(A_S)$?

Strohmer and Vershynin related the expected convergence rate of randomized
Kaczmarz on consistent systems to the smallest singular value and
Frobenius norm of the coefficient matrix
\cite[Theorem~2]{StrohmerVershynin2009}. Haddock, Needell, Rebrova, and
Swartworth used bounds over sufficiently large row subsets to analyze
quantile-based iterations for certain random matrices
\cite{HaddockNeedellRebrovaSwartworth2022}.

For a fixed matrix with normalized rows, Steinerberger gave a
convergence criterion in terms of the least singular value over the
retained subsets and the largest singular value of the whole matrix
\cite[Sections~2.1--2.2]{Steinerberger2023}. The asymptotic
constant he proposed for spherical rows in Section~2.3 is
the basis of both questions below.

\subsection{The sharp constants}
Take $A\in\R^{m\times n}$ with $m,n\ge1$ and unit-norm rows. We
keep $k=\lfloor\theta m\rfloor$ rows for a retained fraction
$\theta\in(0,1)$. This is $q-\beta$ in Steinerberger's notation,
where $q$ is the quantile and $\beta$ is the corruption fraction.
For the submatrix $A_S$ on rows $S$, use the variational
convention $\sigma_{\min}(A_S)=\min_{\norm{x}_2=1}\norm{A_Sx}_2$.
This is zero when $A_S$ has a nontrivial kernel, including wide
and empty submatrices. For each
fixed $x$, we can choose $S$ by keeping the $k$ entries of $Ax$ of
smallest magnitude. Hence
\[
s_\theta(A)=\min_{|S|=k}\sigma_{\min}(A_S)
=\min_{\norm{x}_2=1}\left(\sum_{i=1}^k(|Ax|^2)_{(i)}\right)^{1/2},
\]
where $(i)$ indexes the squared projections in nondecreasing
order. Thus $s_\theta(A)$ plays the role of Steinerberger's
$\sigma_{q-\beta,\min}(A)$ with the row count rounded down and
the sum taken as zero when $k=0$.

For $G\sim N(0,1)$ with density $\varphi$, choose $a_\theta>0$
so that $[-a_\theta,a_\theta]$ has probability $\theta$.
Since $(x\varphi(x))'=\varphi(x)-x^2\varphi(x)$, the integral
of $x^2\varphi(x)$ on $[-a_\theta,a_\theta]$ is the probability
of that interval minus $2a_\theta\varphi(a_\theta)$. Write
\[
h_\theta=\E[G^2\mathbf1_{|G|\le a_\theta}]=\theta-2a_\theta\varphi(a_\theta),\qquad c_\theta=\sqrt{h_\theta}.
\]

\begin{problem}{IE-21}
Let the rows of $A\in\R^{m\times n}$ be independent and uniform on the
unit sphere. Determine whether
\[
 \frac{s_\theta(A)^2}{\norm A_2^2}\xrightarrow{\mathbb P}h_\theta
\]
holds along every sequence with $n\to\infty$ and $m/n\to\infty$ and
give quantitative bounds for the convergence.
\end{problem}

\begin{problem}{IE-22}
Let $M_{m,n}(\theta)$ be the supremum of $\sqrt{n/m}\,s_\theta(A)$ over
all real $m\times n$ matrices with unit-norm rows. Is $c_\theta$ the
smallest constant $c$ such that for every $\eps>0$ there exist integers
$N,R$ with
\[
 M_{m,n}(\theta)\le c+\eps\qquad(n\ge N,\ m/n\ge R)?
\]
\end{problem}

The squared projections in a fixed direction, multiplied by $n$,
approach the law of $G^2$ for uniform spherical rows. The bottom
fraction $\theta$ of these values corresponds to $|G|\le a_\theta$,
so its contribution to the average tends to $h_\theta$.
The minimum over $x$ is why a lower bound for $s_\theta(A)$ needs
this estimate uniformly over the sphere, and any single direction
gives an upper bound.
\begin{theorem}[Random rows]\label{thm:row-random}
If the rows are independent and uniform on $S^{n-1}$, then, along every sequence $n\to\infty$, $m/n\to\infty$, the following limits hold in probability:
\[
\frac nm s_\theta(A)^2\to h_\theta,\qquad \frac nm\norm{A}_2^2\to1,\qquad
\frac{s_\theta(A)^2}{\norm{A}_2^2}\to h_\theta.
\]
\end{theorem}
\begin{theorem}[Universal constant]\label{thm:row-universal}
Uniformly in all integers $m\ge1$ and matrices $A\in\R^{m\times n}$ with unit-norm rows,
\[\sqrt{n/m}\,s_\theta(A)\le c_\theta+o_{n\to\infty}(1).\]
The constant $c_\theta$ is optimal, with the supremum of
$\sqrt{n/m}\,s_\theta(A)$ over matrices with unit-norm rows
converging to $c_\theta$ along every sequence with $n\to\infty$
and $m/n\to\infty$.
\end{theorem}

\subsection{Proof for random rows}
Let $y_{(1)}\le\cdots\le y_{(m)}$ be the observations, all nonnegative,
and put $z_+=\max(z,0)$. The sum of the $k$ smallest observations can
be written as a maximization in one cutoff $t$,
\[
\frac1m\sum_{i=1}^k y_{(i)}=\max_{t\ge0}\left\{\frac km t-\frac1m\sum_i(t-y_i)_+\right\}.
\]
For any $t$, the first $k$ observations contribute at least
$kt-\sum_{i=1}^k y_{(i)}$ to the sum of positive parts. This gives the
upper bound on the maximum in the display. At $t=y_{(k)}$
when $k>0$, the sum of positive parts is exactly
$kt-\sum_{i=1}^k y_{(i)}$, giving equality. For $k=0$,
equality holds at $t=0$.

For a nonnegative random variable $Y$, define
\[
 H_\theta(Y)=\sup_{t\geq0}\{\theta t-\E(t-Y)_+\}.
\]
If $Y$ has a continuous law, choose $b_\theta$ with
$\mathbb P(Y\le b_\theta)=\theta$. The event $\{Y\le b_\theta\}$
plays the role of the $k$ smallest observations. For any $t$,
$(t-Y)_+\ge(t-Y)\mathbf1_{Y\le b_\theta}$, so the objective is
at most $\E[Y\mathbf1_{Y\le b_\theta}]$. At $t=b_\theta$ the
inequality is an equality. Thus $H_\theta(Y)=\E[Y\mathbf1_{Y\le b_\theta}]$.

Below we compare a spherical projection with a Gaussian one by
coupling them. If $Y$ and $Z$ are coupled, the $1$-Lipschitz map
$y\mapsto(t-y)_+$ gives a bound of $\E|Y-Z|$ on the change in
the objective at each $t$. Taking suprema gives
\[|H_\theta(Y)-H_\theta(Z)|\le\E|Y-Z|.\]

Take $u$ uniform on $S^{n-1}$ for $n\ge2$, and set $Y^{(n)}=nu_1^2$
and $h_{\theta,n}=H_\theta(Y^{(n)})$. By rotating the coordinates,
$Y^{(n)}$ describes the squared projection of a random row in any
given unit direction, scaled by $n$. This law is continuous and
$\E Y^{(n)}=1$.

Keep the direction $u$ and give it the length of a Gaussian vector.
That is, take $R^2\sim\chi_n^2$ independently of $u$. Then $Ru$ is
standard Gaussian and its first component $G_1$ satisfies
$G_1^2=(R^2/n)Y^{(n)}$. The absolute difference is then
$Y^{(n)}|1-R^2/n|$. Its expectation factors by independence,
and $\E Y^{(n)}=1$. By Cauchy--Schwarz, the mean of $|1-R^2/n|$
is at most the square root of $\E(1-R^2/n)^2=2/n$. Applying
the bound for $H_\theta$ gives
\[
\E|Y^{(n)}-G_1^2|=\E|1-R^2/n|\le\sqrt{2/n},\qquad
|h_{\theta,n}-h_\theta|\le\sqrt{2/n}.
\]
Set $L=2/(1-\theta)$ and fix $0<\epsilon\le(1-\theta)/2$. For independent
copies $Y_1,\ldots,Y_m$ of $Y^{(n)}$, let $T_m=m^{-1}\sum_{i=1}^kY_{(i)}$.
Let $b$ be the $\theta$-quantile of $Y^{(n)}$. The number of
sample values below $b$ fluctuates around $m\theta$. The
sum of those values fluctuates around $mh_{\theta,n}$.
The cutoff is bounded by $1/(1-\theta)$, since
$b\mathbb P(Y^{(n)}>b)\le\E Y^{(n)}=1$. The individual count
terms are in $[0,1]$ and the sum terms in $[0,b]$, so
Hoeffding's inequality~\cite[Theorem 2]{Hoeffding1963}
gives the following bounds with probability at least
$1-4e^{-2m\epsilon^2}$,
\[
\left|m^{-1}\sum_i\mathbf1_{Y_i\le b}-\theta\right|\le\epsilon,\qquad
\left|m^{-1}\sum_iY_i\mathbf1_{Y_i\le b}-h_{\theta,n}\right|\le b\epsilon.
\]
To turn the sum at $b$ into the sum defining $T_m$, add
or remove up to $m\epsilon+1$ observations, by the count
bound above and the rounding in $k$. All observations
involved are at most $L$ if at least $k$ sample values
are below $L$, since $b\le L$. Markov's inequality gives
$\mathbb P(Y^{(n)}\le L)\ge(1+\theta)/2$, so the mean count
at $L$ exceeds $k$ by at least $m(1-\theta)/2$. Hoeffding's
inequality bounds the chance of a count less than $k$ by
$e^{-m(1-\theta)^2/2}\le e^{-2m\epsilon^2}$. Where the count
bounds hold, correcting the sum and dividing by $m$ adds
at most $L(\epsilon+1/m)$ to the error $b\epsilon$ from the
display. If $k=0$, $T_m=0$ and $h_{\theta,n}\le\theta<1/m$.
It follows that
\[
\mathbb P\bigl(|T_m-h_{\theta,n}|>2L\epsilon+L/m\bigr)\le5e^{-2m\epsilon^2}.
\]

The trimmed sum is the minimum of quadratic forms from
retained row sets. Each of their matrices is bounded in
semidefinite order by $C=(n/m)A^TA$. A bound on $C$ will thus
control the variation of the trimmed sum over the sphere.
In a fixed direction, the quadratic form of $C$ is the
sample mean of $Y^{(n)}$. For $r\ge1$, its summands have moments
\[
\E (Y^{(n)})^r=\frac{n^r(2r-1)!!}{n(n+2)\cdots(n+2r-2)}\le2^r r!,
\]
and hence $\E|Y^{(n)}-1|^r\le4^r r!$. This gives $\E e^{\lambda(Y^{(n)}-1)}\le e^{32\lambda^2}$ for $|\lambda|\le1/8$. Chernoff's inequality applied to the independent copies then yields, for $0<u\le1$,
\[
\mathbb P\left(\left|m^{-1}\sum_iY_i-1\right|>u\right)\le2e^{-mu^2/128}.
\]
For a symmetric matrix $H$, the difference of its quadratic forms
at unit vectors $x,v$ is $(x-v)^TH(x+v)$. If $\norm{x-v}_2\le1/4$,
this changes the form by at most $\norm{H}_2/2$. Thus a $1/4$-net
captures an absolute quadratic form at least half as large as
the maximum on the sphere. Use a net of at most $9^n$ points.
With $H=C-I$, the preceding tail estimate at $t/2$ and a
union bound give
\[\mathbb P(\norm{C-I}_2>t)\le2\,9^n e^{-mt^2/512}\quad(0<t\le1).\]

Write $T(x)=(n/m)\min_{|S|=k}\norm{A_Sx}_2^2$. If
$\norm{C-I}_2\le t$, every retained matrix
$C_S=(n/m)A_S^TA_S$ has norm at most $1+t$, since
$0\preceq C_S\preceq C$. As above, subtracting the quadratic
forms at $x$ and $y$ bounds their difference in absolute value
by $2(1+t)\norm{x-y}_2$, uniformly in $S$. A minimum of functions
with this common Lipschitz constant has the same constant, so
\[|T(x)-T(y)|\le2(1+t)\norm{x-y}_2\]
for unit $x,y$. Applying the fixed-direction estimate on a $\delta$-net of at most $(1+2/\delta)^n$ points and extending by this Lipschitz bound gives
\[
\sup_{\norm{x}_2=1}|T(x)-h_{\theta,n}|\le D:=2L\epsilon+L/m+2(1+t)\delta
\]
except with probability
\[
2\,9^n e^{-mt^2/512}+5(1+2/\delta)^n e^{-2m\epsilon^2},
\]
where $0<t<1$, $0<\delta<1$, $0<\epsilon\le(1-\theta)/2$.

The uniform bound on $T(x)$ also bounds its minimum
$(n/m)s_\theta(A)^2$. Use the Gaussian comparison to replace
$h_{\theta,n}$ by $h_\theta$. For the ratio, the same event gives
$|(n/m)\norm{A}_2^2-1|\le t$, and $h_\theta\le1$. With the same
exceptional probability,
\[
\left|\frac nm s_\theta(A)^2-h_\theta\right|\le D+\sqrt{2/n},\qquad
\left|\frac{s_\theta(A)^2}{\norm{A}_2^2}-h_\theta\right|
\le\frac{D+\sqrt{2/n}+t}{1-t}.
\]
Write $Q=m/n$ and take $t=\epsilon=\delta=32\sqrt{\log Q/Q}$ for
sufficiently large $Q$. The exponential factors in the two probability
bounds become $Q^{-2n}$ and $Q^{-2048n}$. The two nets have at most $9^n$
and $Q^n$ points for large $Q$, so the exceptional probabilities tend to
zero. The error is
\[
O_\theta\left(\sqrt{\frac nm\log\frac mn}+n^{-1/2}+m^{-1}\right).
\]
This proves Theorem~\ref{thm:row-random} for every sequence in its stated regime.

\subsection{The universal bound and sharpness}
For an arbitrary matrix with normalized rows, it suffices to find one
direction with small retained energy. We will pick it from a subspace where $A$ is
small in operator norm. This limits the fluctuations of the Gaussian
threshold estimate and lets us use the same direction for every cutoff.

Choose $1\le r<n$ and put $d=n-r$. In orthonormal coordinates on the
complement of the leading $r$ eigenvectors of $A^TA$, write the restricted
matrix as $B\in\R^{m\times d}$ with rows $b_i^T$. The trace of $A^TA$ is $m$,
so its $(r+1)$st eigenvalue is at most $m/(r+1)$. Projecting the rows can
only reduce their norms. Thus
\[
\norm{B}_2^2\le m/(r+1),\qquad\norm{b_i}_2\le1,\qquad\tr(B^TB)\le m.
\]
Draw $g\sim N(0,I_d)$ and set $Y_i=(b_i^Tg)^2$. To express the sum
of the $k$ smallest observations divided by $m$, we maximize
the following expression over $t\ge0$.
\[\Psi_t(g)=\frac km t-\frac1m\sum_i(t-Y_i)_+.\]
Each $Y_i$ is distributed as $\norm{b_i}_2^2G^2$ for a standard
Gaussian $G$. As $\norm{b_i}_2\le1$, the expectation of $(t-Y_i)_+$
is at least $\E(t-G^2)_+$. Taking expectations term by term and
using $k/m\le\theta$ gives
\[
\E\Psi_t(g)\le\theta t-\E(t-G^2)_+\le h_\theta.
\]

For each fixed $t$, the Gaussian Poincar\'e inequality bounds the
variance of $\Psi_t(g)$ by $\E\norm{\nabla\Psi_t(g)}_2^2$. To bound
this gradient, differentiate $z\mapsto(t-z^2)_+$. Its derivative
has magnitude at most $2\sqrt t$ wherever it exists. The gradient
of $\Psi_t$ is $-B^T/m$ times the vector of these $m$ derivatives, so
\[
\operatorname{Var}(\Psi_t(g))\le\frac{\norm{B}_2^2}{m^2}(4tm)\le\frac{4t}{r+1}.
\]

First confine a maximizing cutoff to a fixed interval. The untrimmed
mean $Z=m^{-1}\norm{Bg}_2^2$ has expectation at most one and variance
$2\tr((B^TB)^2)/m^2\le2/(r+1)$. Chebyshev's inequality gives $Z\le2$
except with probability at most $2/(r+1)$. If $k>0$, the $k$th
smallest observation is a maximizing cutoff, and at least
$m-k+1\ge(1-\theta)m$ observations are that large. On $Z\le2$ the cutoff
is therefore at most $L=2/(1-\theta)$. If $k=0$, use the cutoff zero.

The cutoff objective $\Psi_t(g)$ is piecewise linear in $t$,
with all slopes in $[-1,1]$. Put a grid of at most $L/\delta+2$
points on $[0,L]$, with spacing at most $\delta$ for
$0<\delta<1$. We bound the probability that any grid
value exceeds $h_\theta+\delta$. Chebyshev's inequality
gives at most $4L/((r+1)\delta^2)$ at each point. The total
failure probability is at most
\[\frac{4L(L/\delta+2)}{(r+1)\delta^2}.\]
Between grid points $\Psi_t(g)$ can rise by at most another
$\delta$. Hence on the good event it is bounded by
$h_\theta+2\delta$ throughout $[0,L]$.

The retained sum divided by $m$ scales by $1/\norm{g}_2^2$
when we normalize $g$.
The mean and variance of $\norm{g}_2^2$ are $d$ and $2d$, so
$\norm{g}_2^2\ge(1-\delta)d$ except with probability at most
$2/(d\delta^2)$. We can meet this bound, $Z\le2$, and the grid bound
with the same $g$ if the sum of the three failure probabilities is
less than one, that is,
\[
\frac2{r+1}+\frac{4L(L/\delta+2)}{(r+1)\delta^2}
+\frac2{(n-r)\delta^2}<1.
\]
Use such a $g$ to form a unit vector in the chosen subspace. Its
trimmed squared projections are scaled by
$1/\norm{g}_2^2\le1/((1-\delta)(n-r))$, giving
\[{\frac nm s_\theta(A)^2\le\frac{n}{(n-r)(1-\delta)}(h_\theta+2\delta).}\]
With $r=\lfloor n^{2/3}\rfloor$ and $\delta=n^{-1/6}$, the grid term
in the failure probability is $O_\theta(1/(r\delta^3))$. This is
$O_\theta(n^{-1/6})$, and the other two terms are smaller. Using $r/n=O(n^{-1/3})$ in the
bound gives $(n/m)s_\theta(A)^2\le h_\theta+O_\theta(n^{-1/6})$,
with constants depending only on $\theta$.

Along any sequence with $n\to\infty$ and $m/n\to\infty$,
Theorem~\ref{thm:row-random} gives
matrices of spherical rows with $(n/m)s_\theta(A)^2>h_\theta-\epsilon$
for all sufficiently large $n$, for each $\epsilon>0$. Thus the supremum of
$\sqrt{n/m}\,s_\theta(A)$ over matrices with unit-norm rows has liminf
at least $c_\theta$.

\section{A strict pointwise tensor-train bound (TR-04)}\label{sec:tr04}

The truncated SVD gives a best Frobenius-norm approximation at a
prescribed matrix rank. For a tensor with entries $A(i_1,\ldots,i_d)$,
mode $j$ is the $j$th index position, and its size $n_j$ is the number
of values $i_j$ can take. A tensor train writes each entry as a
matrix product $G_1(i_1)\cdots G_d(i_d)$. For each mode $j$,
we store $n_j$ matrices $G_j(i_j)$ of size $s_{j-1}\times s_j$. The
endpoint dimensions $s_0=s_d=1$ give scalar output, and the total storage
is $\sum_{j=1}^d n_j s_{j-1}s_j$ entries.

A cut after mode $j$ splits the product between factors $j$ and $j+1$.
The \emph{unfolding} at this cut is the matrix with entry
$A(i_1,\ldots,i_d)$ at row $(i_1,\ldots,i_j)$ and column
$(i_{j+1},\ldots,i_d)$. Each tensor entry appears once, so this
reindexing preserves the Frobenius norm~\cite{nmKoldaBader2009}.
Multiplying the factors before the cut gives a row vector of length
$s_j$ for each row, and multiplying the factors after the cut gives a
column vector of length $s_j$ for each column. These vectors factor
the unfolding through the shared dimension $s_j$, so $s_j$ is at least
the rank of that matrix, called the tensor-train rank at cut $j$.
A nonzero tensor $A$ has an exact tensor-train representation in which
the shared dimension at every cut equals the rank of the unfolding at
that cut~\cite[Section~2]{ftOseledets2011}.

TT-SVD builds the train from left to right, applying a truncated SVD
at each cut to the tensor compressed by the earlier steps. The discarded components are
mutually orthogonal in the original tensor space, so their squared
norms add. At each cut, applying the preceding compressions to a
best admissible tensor gives a candidate of the required rank
with squared error at most the original optimum. With $d-1$ cuts,
the usual factor is $d-1$ in squared error, or $\sqrt{d-1}$ in
Frobenius norm~\cite[Theorem~2.2 and Corollary~2.4]{ftOseledets2011}.

If the last singular value kept by the first SVD and the first one
discarded are equal and positive, there is a choice of which left
singular directions to keep at that value. Every choice of the required
number of directions gives the same squared Frobenius error, the sum
of the squared discarded singular values. The remaining SVDs work with
the tensor compressed to the retained subspace, so the later
truncation errors can depend on that choice. We use this freedom to
improve the squared-error bound at fixed ranks. The broader problem asks
for the same improvement when the factors are arranged on a tree,
with one matrix truncation at each edge
\cite[Problem~6.1]{ftAmsel2026}.

\subsection{The approximation problem and the first projection}
For tensors of order $d\ge3$ with mode sizes $n_1,\ldots,n_d\ge2$,
prescribe positive integer rank limits
$r=(r_1,\ldots,r_{d-1})$. Let $S_r$ be the set of tensors in
$\mathbb R^{n_1\times\cdots\times n_d}$ whose unfolding at cut
$j$ has rank at most $r_j$ for each cut. Let $A$ be given as
a dense tensor and write its optimal squared error as
\[
 E_* = \min_{Y\in S_r}\|A-Y\|_{\mathrm F}^2.
\]
The error attains its minimum on the compact set
$\{Y\in S_r:\norm{A-Y}_{\mathrm F}\le\norm A_{\mathrm F}\}$,
since the zero tensor is feasible and the vanishing minors at
each cut make $S_r$ closed.

We count exact arithmetic and SVD operations on the dense input,
with exact comparisons of singular values.

\begin{problem}{TR-04}
In this arithmetic and SVD model, does an algorithm polynomial in
the dense input size and rank parameters always return $X\in S_r$
with
\[
 \norm{A-X}_{\mathrm F}^2<(d-1)E_*\quad\text{if }E_*>0,
 \qquad X=A\quad\text{if }E_*=0?
\]
\end{problem}

\begin{theorem}\label{thm:tr04-main}
In this arithmetic and SVD model, there is a deterministic algorithm,
polynomial in the dense input size and rank parameters, which returns
$X\in S_r$ such that
\[
 \|A-X\|_{\mathrm F}^2<(d-1)E_*\quad\hbox{if }E_*>0,
 \qquad X=A\quad\hbox{if }E_*=0.
\]
It uses at most $n_1$ ordinary TT-SVD completions and does not increase any prescribed rank.
\end{theorem}

At the first cut, write $A_{(1)}\in\mathbb R^{n_1\times(n_2\cdots n_d)}$
for the unfolding. The best rank-$r_1$ matrix error has square
\[
 e^2=\sum_{i>r_1}\sigma_i(A_{(1)})^2,
\]
with zeros appended to the singular values as necessary.
A tensor $Y$ attaining $E_*$ gives a rank-at-most-$r_1$
matrix $Y_{(1)}$ with the same squared error. Thus the first
cut costs $e^2\le E_*$. The standard bound allows $E_*$ at
each cut, so $e^2<E_*$ already gives a strict inequality.
If $e^2=E_*$, we will get a saving by choosing a first
projection that reduces the norm of $A-Y$ before the later cuts.

In a product such as $P_W A$, the orthogonal projector onto
$W\subseteq\mathbb R^{n_1}$ acts on the first index of the tensor.
The same convention applies to other matrices on the left.
\begin{lemma}\label{lem:tr04-completion}
Let $U\in\mathbb R^{n_1\times k}$ have orthonormal columns, where $1\le k\le r_1$, and set $P=UU^T$. Compress $A$ to $U^TA$, group its first two indices, and apply TT-SVD with remaining ranks $r_2,\ldots,r_{d-1}$. If $X_U$ is the lifted result, then $X_U\in S_r$ and, for every $Y\in S_r$,
\begin{equation}\label{eq:tr04-completion}
 \|A-X_U\|_{\mathrm F}^2\le\|(I-P)A\|_{\mathrm F}^2+(d-2)\|P(A-Y)\|_{\mathrm F}^2.
\end{equation}
\end{lemma}
\begin{proof}
Consider the unfolding at a later cut, just before TT-SVD
truncates it. Compress $Y$ in the same way as $A$ up to that point.
Each compression premultiplies the unfolding at this cut, so its
rank stays at most the prescribed one. The compressions after
$U^T$ can only reduce the distance from the compressed $A$ to this
candidate. The resulting distance is at most
$\|U^T(A-Y)\|_{\mathrm F}$, which bounds the cost of the truncation.

Lifting the completed tensor by $U$ preserves its later cut ranks.
The columns of its first unfolding lie in $\operatorname{range} U$,
of dimension at most $r_1$, so $X_U\in S_r$. The $d-2$ later losses
are pairwise orthogonal and retain their norms when lifted by
$U$. The first-cut error is in the orthogonal complement of
$\operatorname{range} U$. Thus the total squared error is bounded
by the first-cut error squared plus $d-2$ copies of
$\|U^T(A-Y)\|_{\mathrm F}^2=\|P(A-Y)\|_{\mathrm F}^2$, as in
\eqref{eq:tr04-completion}.
\end{proof}

\subsection{Choosing what to keep at a tied singular value}
At $A=0$ the algorithm returns zero. For a nonzero input with
$\rank A_{(1)}\le r_1$, take the column space of $A_{(1)}$.
If the rank is larger, let $\tau=\sigma_{r_1}(A_{(1)})>0$.
To attain the best first-cut error, keep the left singular space
$H$ for values greater than $\tau$. The remaining $s$ dimensions
come from the left singular space $E$ at $\tau$, where
\[
 h=\dim H,\qquad t=\dim E,\qquad s=r_1-h,
 \qquad 1\le s\le t.
\]
When $s=t$, the only candidate is $P_H+P_E$. For $s<t$,
fix an orthonormal SVD basis of $E$ indexed by $u_0,\ldots,u_{t-1}$.
Starting at each basis vector, keep that vector and the next $s-1$
in cyclic order. The resulting projectors are given by
\[
 W_j=\operatorname{span}\{u_j,\ldots,u_{j+s-1}\},\qquad
 P_j=P_H+P_{W_j},
\]
for $0\le j<t$, with the subscripts read modulo $t$. Every block
has the same first-cut error, since it retains $s$ directions
of singular value $\tau$. Each direction occurs in $s$ blocks,
so averaging the blocks gives equal weight to the tied directions.

The algorithm compares completions, using \cref{lem:tr04-completion}
once for each first projection and taking the tensor with
least Frobenius error. Use a fixed rule for a tie in this
comparison. There is one completion for a prescribed subspace,
or $t$ for the cyclic choices. Since $t\le n_1$ and each completion
uses unfoldings with at most as many entries as $A$, the
algorithm takes polynomial time in the stated model.

\begin{lemma}\label{lem:tr04-equality}
Suppose $\rank A_{(1)}>r_1$ and a rank-at-most-$r_1$ matrix $Y_{(1)}$ attains the best matrix error $e^2$. Let $Y$ be the tensor with first unfolding $Y_{(1)}$. Then its residual $R=A-Y$ satisfies $P_HR=0$. If $s=t$, it also satisfies $(P_H+P_E)R=0$.
\end{lemma}
\begin{proof}
Let $S$ be the column space of $Y_{(1)}$. The orthogonal decomposition
\[
 \|A_{(1)}-Y_{(1)}\|_{\mathrm F}^2
 =\|(I-P_S)A_{(1)}\|_{\mathrm F}^2+\|P_SA_{(1)}-Y_{(1)}\|_{\mathrm F}^2
\]
gives $Y_{(1)}=P_SA_{(1)}$. Here $\dim S=r_1$: using fewer
dimensions would discard a positive singular value, since
$\sigma_{r_1}(A_{(1)})=\tau>0$.

In a left singular basis $u_i$, put $w_i=\|P_Su_i\|_2^2$. The
projection captures energy $\sum_i w_i\sigma_i^2$, with
$0\le w_i\le1$ and $\sum_i w_i=r_1$. Relative to the baseline
$r_1\tau^2$, the captured energy is $\sum_i w_i(\sigma_i^2-\tau^2)$.
For a positive coefficient the largest contribution uses $w_i=1$.
For a negative coefficient it uses $w_i=0$. Both bounds can be
attained by choosing $H$ and $s$ directions in $E$. Equality
therefore requires
$H\subseteq S\subseteq H\oplus E$. The residual $(I-P_S)A$ projects
to zero on $H$, and also on $E$ if $s=t$.
\end{proof}

\begin{proof}[Proof of \cref{thm:tr04-main}]
Choose a best feasible tensor $Y$ and put $R=A-Y$. The projectors
we try all attain the first-cut error $e^2$. For each of them,
\cref{lem:tr04-completion} gives
\begin{equation}\label{eq:tr04-basic}
 \|A-X_U\|_{\mathrm F}^2\le e^2+(d-2)\|PR\|_{\mathrm F}^2
 \le e^2+(d-2)E_*.
\end{equation}
When $e^2<E_*$, the last expression is already strictly below $(d-1)E_*$.

Suppose $e^2=E_*>0$. The first unfolding has rank greater than
$r_1$, and $Y_{(1)}$ attains its best rank-$r_1$ error. Thus
$P_HR=0$ by \cref{lem:tr04-equality}. When $s=t$, we also have $P_ER=0$.
Then $PR=0$ and \eqref{eq:tr04-basic} gives error at most
$E_*<(d-1)E_*$.

When $s<t$, the residual after a first projection is $P_jR=P_{W_j}R$.
Each tied direction occurs in $s$ of the $t$ tested subspaces, so
the averages of the projectors and squared residuals are
\[
 \frac1t\sum_{j=0}^{t-1}P_{W_j}=\frac{s}{t}P_E,
 \qquad
 \frac1t\sum_{j=0}^{t-1}\|P_jR\|_{\mathrm F}^2
 =\frac{s}{t}\|P_ER\|_{\mathrm F}^2\le\frac{s}{t}E_*<E_*.
\]
One of the tested projections therefore keeps at most $(s/t)E_*$
of the squared residual. By \eqref{eq:tr04-basic}, its squared
error is at most $[1+(d-2)s/t]E_*<(d-1)E_*$, and the selected
completion has error at most this much.

For $E_*=0$, use $Y=A$ in \eqref{eq:tr04-basic}. Both $e^2$ and
the projected residual are zero, which forces $X_U=A$.
\end{proof}

\section{Rank and border rank under tensor powers (TR-27)}\label{sec:tr27}

The number of pure tensors $v_1\otimes\cdots\otimes v_d$ required to represent a given tensor is its tensor rank. This concept is more broadly defined as $X$-rank, where the elements of the variety $X$ correspond to the summands. Let $V$ be a finite dimensional complex vector space and let $X\subset\mathbb P(V)$ be a reduced and irreducible variety which spans $\mathbb P(V)$. The points of $X$ define the ‘directions’ of the summands. For a vector $v\neq 0$, let $p=[v]$ be the line through $v$ and let $R_X(p)$ be the smallest integer $r$ such that $v=w_1+\cdots+w_r$ with $[w_i]\in X$ for all $i$. Ordinary tensor rank is obtained when $X$ consists of the lines spanned by pure tensors. The border rank $\underline R_X(p)$ is the smallest integer $s$ such that $v$ can be approximated by sums of at most $s$ summands. This can be smaller than the exact rank because the individual summands may become large and cancel each other while the sum converges. Geometrically, sums of at most $s$ terms correspond to points in the linear spans of $s$ points of $X$. The Zariski closure of this set is called the $s$th secant variety, denoted $\sigma_s(X)$. Over $\mathbb C$, this closure is the same as the limit closure, and hence $\underline R_X(p)$ is the smallest $s$ such that $p\in\sigma_s(X)$ \cite[Section~2]{ftOnetoVentura2025}.

For $p^{\otimes k}=[v^{\otimes k}]$, the admissible terms are
$x_1\otimes\cdots\otimes x_k$ with each $[x_j]\in X$ chosen
independently. These terms give the Segre image of $X^{\times k}$
in $\mathbb P(V^{\otimes k})$, with a separate copy of $V$ for each
factor. Expanding an $r$-term decomposition of $v$ in all $k$
factors gives the upper bound $R_{X^{\times k}}(p^{\otimes k})\le r^k$.
For ordinary tensor rank, a smaller border rank guarantees that
this bound is strict for some power~\cite{nmChristandlJensenZuiddam2018}.
The same holds for general $X$, since the $k$th roots of rank and
border rank of $p^{\otimes k}$ tend to a common limit at most
$\underline R_X(p)$. Border rank itself can also be strictly
submultiplicative
\cite[Theorem~1.1 and Section~6]{nmChristandlGesmundoJensen2019}.

The following conjecture asks whether the saving must already appear in the
square. This holds for binary forms and ternary cubics when the
summands are powers of linear forms
\cite[Conjecture~1.1 and Section~4]{ftBBGOV2021}.

\begin{problem}{TR-27}
For every reduced irreducible complex projective variety $X$ spanning
$\mathbb P(V)$ and every $p\in\mathbb P(V)$, does
\[
 \underline R_X(p)<R_X(p)
 \quad\Longrightarrow\quad
 R_{X\times X}(p^{\otimes2})<R_X(p)^2
\]
hold, with $X\times X$ embedded by the Segre map
\cite[Conjecture~4.9]{ftOnetoVentura2025}?
\end{problem}

We provide a counterexample that lies on a tangent line to a curve. The
point has border rank two and rank three relative to the
curve, and its square has rank nine. By choosing the curve
we can delay the first saving past any prescribed power.

\begin{theorem}\label{thm:tr27-main}
Given integers $r\ge3$ and $m\ge2$, set $D=r^m+r$. There is a smooth rational nondegenerate curve $X\subset\mathbb P^{D-1}$ and a point $p$ satisfying
\[
 \underline R_X(p)=2,\qquad R_X(p)=r,\qquad
 R_{X^{\times k}}(p^{\otimes k})=r^k\quad(1\le k\le m).
\]
For $r=3$, $m=2$ the curve lies in $\mathbb P^{11}$ and disproves the proposed strict inequality at the square.
\end{theorem}

\subsection{Proof of the theorem}
Projective points are linearly independent if their vector
representatives are linearly independent.
A tangent vector has border rank at most two because it is a
limit of scaled differences of two curve vectors. To give it
rank $r$, project the curve so that the tangent vector becomes a
sum of $r$ curve vectors. We need small sets of projected points to remain
independent. In a proposed short tensor decomposition those
sets then give bases, and the lower bound becomes a count of
nonzero coordinates.

\begin{lemma}\label{lem:tr27-product}
Let $K$ and $k$ be positive integers. Suppose every set of at most $K$ distinct points of $X$ is linearly independent. If $R_X(p)=r$ and $r^k-1\le K$, then $R_{X^{\times k}}(p^{\otimes k})=r^k$.
\end{lemma}
\begin{proof}
Expanding the $k$th tensor power of an $r$-term expression gives the upper bound. For the lower bound, let $v$ represent $p$ and suppose
\[
 v^{\otimes k}=\sum_{\ell=1}^{L}\lambda_\ell
 x_{1\ell}\otimes\cdots\otimes x_{k\ell},\qquad L<r^k.
\]
In each factor $j$, collect a single representative of each
projective point used, absorbing rescalings into $\lambda_\ell$.
The resulting set $A_j$ has at most $L\le K$ vectors, so it is
independent. Apply a functional $\varphi$ with $\varphi(v)=1$
in every factor except $j$. This expresses $v$ in $\operatorname{span}A_j$.
At least $r$ of its coordinates in $A_j$ must be nonzero, since
$R_X(p)=r$. In the product basis, any choice of one of these
coordinates per factor gives a nonzero coordinate of $v^{\otimes k}$.
There are at least $r^k$ of them. The right side therefore needs
at least $r^k$ basis vectors, contrary to $L<r^k$.
\end{proof}

Let $E_D=\mathbb C^{D+1}$ have the ordered basis $e_0,\ldots,e_D$. Define
\[
 c_D(t)=\sum_{j=0}^D t^j e_j\quad(t\in\mathbb C),\qquad
 c_D(\infty)=e_D.
\]
The projective points $[c_D(t)]$, for $t\in\mathbb P^1=\mathbb C\cup\{\infty\}$,
form the rational normal curve $C_D\subset\mathbb P(E_D)$. Any set of at most
$D+1$ distinct finite parameters gives independent vectors by the Vandermonde
determinant. If one parameter is infinity, the first coordinates show that
the other vectors are independent, and the last coordinate separates $e_D$.

\begin{lemma}\label{lem:tr27-tangent}
For $D\ge2$, $R_{C_D}([e_1])\ge D$.
\end{lemma}
\begin{proof}
Suppose $e_1$ lies in the span of the $c_D(a)$ for a set
$S\subset\mathbb P^1$ of $\ell\le D-1$ distinct parameters. A linear
functional with coefficients $f_j$ acts on $c_D(a)$ at a finite $a$
by evaluating the polynomial $f(T)=\sum_{j=0}^D f_jT^j$. It acts on $e_1$
by taking $f'(0)$. To separate $e_1$ from the proposed span, take
\[
 f(T)=T\prod_{a\in S\setminus\{0,\infty\}}(T-a)
\]
This vanishes at the finite parameters in $S$ and has a simple zero
at the origin. Its degree is at most $D$, and at most $D-1$ when
$\infty\in S$, so in that case its coefficient of $T^D$ also vanishes.
The functional $L_f(e_j)=f_j$ is therefore zero on the proposed span
and nonzero on $e_1$.
\end{proof}

To make $e_1$ a sum of $r$ curve vectors after projection, choose
distinct nonzero parameters $a_1,\ldots,a_r$ and set
\[
 b=\sum_{i=1}^r c_D(a_i),\qquad z=e_1-b,\qquad
 V=E_D/\mathbb Cz,
\]
with quotient map $\pi:E_D\to V$. The zeroth coordinate of $z$ is
$-r$, so the quotient removes one dimension. If we add the $r$
terms of $b$ to an expression for $z$, we obtain an expression
for $e_1$. The latter needs at least $D$ terms by
\cref{lem:tr27-tangent}, so
\begin{equation}\label{eq:tr27-centre}
 R_{C_D}([z])\ge D-r=r^m.
\end{equation}
We need the independence in \cref{lem:tr27-product} with
$K=r^m-1$. Take any $s\le K$ distinct points of $C_D$.
Their representatives are independent by the Vandermonde
calculation above. If the span contained $z$, it would express
$z$ using only $s$ points, contrary to \eqref{eq:tr27-centre}.
Thus $\pi$ is injective on this span and the $s$ images remain
independent. Applying this to one or two points shows that the
map on $C_D$ is defined everywhere and injective. Its reduced
image $X$ is closed and irreducible, and spans $\mathbb P(V)$
because $C_D$ spans $\mathbb P(E_D)$ and $\pi$ is onto.

Since the quotient identifies $e_1$ with $b$, the vector $v=\pi(e_1)$ satisfies
\begin{equation}\label{eq:tr27-r-terms}
 v=\sum_{i=1}^r\pi(c_D(a_i)).
\end{equation}
Take the union of the points in \eqref{eq:tr27-r-terms} and any
set of at most $r-1$ points of $X$. This gives at most $2r-1\le K$
points, whose representatives are independent. In the resulting
basis of their span, $v$ has $r$ nonzero coordinates by
\eqref{eq:tr27-r-terms}. This rules out a shorter expression
and gives $R_X([v])=r$.

For the border rank, approximate $v$ by the difference
quotient of two projected curve vectors,
\[
 \frac{\pi(c_D(t))-\pi(c_D(0))}{t}
 =v+\sum_{j=2}^D t^{j-1}\pi(e_j)\longrightarrow v
 \quad(t\to0).
\]
For $t\ne0$, the difference quotient represents a point of
$X$-rank at most two. Every homogeneous polynomial that vanishes
on all such points vanishes on the quotients and hence, by
continuity, at $v$. This proves $\underline R_X([v])\le2$.
Since $[v]$ has rank $r>1$
and $X$ is closed, its border rank is exactly two. The ranks
through power $m$ now follow from \cref{lem:tr27-product}.

For the image to be smooth, we also need to preserve tangent
directions. The obstruction is a center on a tangent line,
just as a center on a secant line would identify points.
Both kinds of lines lie in the border-rank-at-most-two locus.
A Hankel matrix gives the needed lower bound for the border
rank of the center $[z]$.

Write $z_j$ for the $e_j$-coordinate of $z$ and $\delta_{j1}$ for
the indicator of $j=1$. Since $D\ge2r+2$, these coordinates give
the Hankel matrix
\[
 H=(z_{i+j})_{0\le i,j\le r+1},\qquad
 z_j=\delta_{j1}-\sum_{\nu=1}^r a_\nu^j.
\]
For a curve vector with finite parameter $a$, this construction
gives entries $a^{i+j}$ and hence a rank-one matrix. The vector at
infinity gives at most one nonzero entry. It follows that a sum of
$s$ curve terms gives a Hankel matrix of rank at most $s$. The
minors expressing this condition still vanish on the Zariski
closure. Thus $\rank H$ is a lower bound for the border rank of $z$.

To compute the rank of $H$, factor it using the vectors
$u(a)=(1,a,\ldots,a^{r+1})^T$. Set
\[
 B=[u(0)\ \ u'(0)\ \ u(a_1)\ \cdots\ \ u(a_r)],\qquad
 J=\begin{pmatrix}0&1\\1&0\end{pmatrix}\oplus(-I_r).
\]
The $e_1$ term of $z$ contributes
$u(0)u'(0)^T+u'(0)u(0)^T$ to $H$, while each $-c_D(a_i)$
contributes $-u(a_i)u(a_i)^T$. Thus $H=BJB^T$. A vector in the
left kernel of $B$ gives a polynomial of degree at most $r+1$
with a double root at zero and $r$ distinct roots at the $a_i$.
It must vanish identically, so $B$ is invertible. Since $J$ is
invertible as well, $\rank H=r+2$ and $\underline R_{C_D}([z])\ge r+2$.

Since $\underline R_{C_D}([z])\ge r+2>2$, the projection center lies
off every secant and tangent line, including the ones at infinity.
Projection separates distinct points and preserves their tangent
directions. The resulting proper injective immersion of the smooth
projective curve is an embedding with smooth rational image.

\subsection{An explicit counterexample}
Take $r=3$, $m=2$, and $D=12$, choosing $a_i=i$ for $1\le i\le3$. Set
\[
 S_j=1+2^j+3^j,\qquad J_0=\{0,2,3,\ldots,12\}.
\]
We can write the quotient map in coordinates by eliminating
coordinate $1$ along $\mathbb C z$. With the target coordinates
ordered by $J_0$, take $P:\mathbb C^{13}\to\mathbb C^{12}$ given by
\[
 (P(w))_j=5w_j-S_jw_1\quad(j\in J_0).
\]
Indeed $z_1=-5$ and $z_j=-S_j$ for $j\in J_0$, so $Pz=0$. The
columns indexed by $J_0$ form $5I_{12}$, making $P$ onto with
kernel exactly $\mathbb C z$.

Define $g(t)=(5t^j-S_jt)_{j\in J_0}$ and $v=(-S_j)_{j\in J_0}$. Then
\[
 g'(0)=v,
 \qquad g(1)+g(2)+g(3)=v,
 \qquad v^{\otimes2}=\sum_{a=1}^3\sum_{b=1}^3 g(a)\otimes g(b).
\]
Here $K=r^m-1=8$, so any set of at most eight distinct projected points
is independent. Since $[v]$ has rank three, \cref{lem:tr27-product} shows
that all nine terms for its square are necessary.

\subsection{Eventual rank savings}
For each fixed curve, a later power does give a saving. We use
the coefficient of a polynomial that is a sum of few curve
terms at every evaluation~\cite[Section 1]{ftBBGOV2021}.
Put $v=\pi(e_1)$, $p=[v]$, and set
\[
 q(t)=\pi(c_D(t)-c_D(0))=tv+t^2\pi(e_2)+\cdots+t^D\pi(e_D).
\]
The coefficient of $t^k$ in $q(t)^{\otimes k}$ is $v^{\otimes k}$.
Its possible degrees run from $k$ to $kD$. With $L=k(D-1)+1$,
only $k$ in that range is congruent to $k$ modulo $L$. To extract
its coefficient, multiply by $\zeta^{-k}$ and average over the
$L$th roots of unity to obtain
\[
 v^{\otimes k}=\frac1L\sum_{\zeta^L=1}\zeta^{-k}q(\zeta)^{\otimes k}.
\]
Each $q(\zeta)$ is a difference of two curve vectors, so its
$k$th power uses at most $2^k$ terms. This gives
\[
 R_{X^{\times k}}(p^{\otimes k})\le[k(D-1)+1]2^k.
\]
The ratio of this bound to $r^k$ is $[k(D-1)+1](2/r)^k$, which
tends to zero for each fixed $D$ and $r>2$. For the curve in
$\mathbb P^{11}$, at $k=13$ the bound is $144\cdot2^{13}<3^{13}$.
The first saving is finite for each fixed curve, and
\cref{thm:tr27-main} shows that these first-saving powers
are unbounded across the examples.

\section{Differentiation need not preserve nonnegative realizability (IS-03)}\label{sec:is03}

A list of numbers is nonnegatively realizable if it is the
spectrum, with algebraic multiplicities, of an entrywise nonnegative
matrix. For such a matrix $A\in\R^{n\times n}$ with characteristic
polynomial $p_A(z)=\det(zI-A)$, the critical points are the roots of
the monic polynomial $p_A'/n$ of degree $n-1$.

\begin{problem}{IS-03}
For every $n\ge5$ and every real entrywise nonnegative matrix $A$ of
order $n$, does there exist a real entrywise nonnegative matrix $B$ of
order $n-1$ with
\[
 \det(zI-B)=\frac1n p_A'(z)
\]
as polynomials in $z$? This is Johnson's derivative-realization
conjecture in the remaining dimensions~\cite[Conjecture~1.2]{Hoover2018}.
\end{problem}

Cronin and Laffey~\cite{CroninLaffey2014} proved the conjecture
for lists of length up to four and for trace-zero lists of length
up to six, as discussed in \cite[Section~1]{Hoover2018}. Further cases
include nonnegative companion matrices~\cite[Theorem 3.6]{Hoover2018}
and nonnegative matrices unitarily diagonalized by a normalized
complex Hadamard matrix~\cite[Corollary 5.13]{Hoover2018}. We provide a counterexample in \cref{thm:is03}.

For a nonnegative realization $B\in\mathbb R^{m\times m}$, the
powers $B^k$ are entrywise nonnegative. Their traces give the
nonnegative power-sum condition on the eigenvalues,
\begin{equation}\label{eq:is03-moment}
 \sum_{j=1}^{m}\mu_j^k=\tr(B^k)\geq0,\qquad k\geq1,
\end{equation}
where $\mu_1,\ldots,\mu_m$ are the eigenvalues counted with algebraic multiplicity.
Monov~\cite[Conjecture 3.1]{Monov2008} conjectured that the roots
of $p_A'$ always pass the trace test when $A$ is entrywise
nonnegative. He also asked if a larger nonnegative matrix
could realize the same list of roots with extra zeros. The
size of that matrix is immaterial in \eqref{eq:is03-moment},
since the added roots contribute $0^k=0$ for $k\ge1$. To find
a negative power sum, we use directed-cycle permutation
matrices $C_k$, with characteristic polynomials $z^k-1$.

\begin{theorem}\label{thm:is03}
The entrywise-nonnegative matrix
\[
 A=\diag\left(\tfrac12,C_2,C_4\right)\in\R^{7\times7}
\]
has the property that the roots of $p_A'/7$ are not the spectrum
of any entrywise nonnegative matrix of order six. Adjoining any number
of zeros still leaves a list with no nonnegative realization.
\end{theorem}

\begin{proof}
Multiplying the characteristic polynomials of the three blocks and differentiating gives
\begin{align}
 p_A(z)&=(z-\tfrac12)(z^2-1)(z^4-1)=z^7-\tfrac12z^6-z^5+\tfrac12z^4-z^3+\tfrac12z^2+z-\tfrac12,\label{eq:is03-p}\\
 q(z)&=\frac{p_A'(z)}7
 =z^6-\frac37z^5-\frac57z^4+\frac27z^3-\frac37z^2+\frac17z+\frac17.\label{eq:is03-q}
\end{align}
Newton's identities compute the power sums directly from the
coefficients of $q$. In general, if the monic polynomial is
$z^d+c_1z^{d-1}+\cdots+c_d$, write $s_k$ for the $k$th power
sum of its roots, counted with multiplicity. Then
\begin{equation}\label{eq:is03-newton}
\begin{cases}
 s_k+\sum_{j=1}^{k-1}c_js_{k-j}+kc_k=0,&1\leq k\leq d,\\
 s_k+\sum_{j=1}^{d}c_js_{k-j}=0,&k>d.
\end{cases}
\end{equation}
For $q$ in \eqref{eq:is03-q}, the first six values are
\[
\begin{array}{lll}
 s_1=3/7,&s_2=79/49,&s_3=48/343,\\
 s_4=6731/2401,&s_5=5213/16807,&s_6=219766/117649.
\end{array}
\]
Since $q$ has degree six, the $k>d$ case of \eqref{eq:is03-newton} gives
\[
 7s_7=3s_6+5s_5-2s_4+3s_3-s_2-s_1.
\]
Hence
\begin{equation}\label{eq:is03-negative}
 s_7=-\frac{8593}{823543}<0.
\end{equation}
The roots of $q$, even with extra zeros, therefore violate
the trace condition \eqref{eq:is03-moment}.
\end{proof}

\section{A null vector for the Sinkhorn coefficient identity (NM-04)}\label{sec:nm04}

For an $m\times n$ matrix $A=(a_{ij})$ with strictly positive entries,
its Sinkhorn limit $\operatorname{Sink}(A)$ is the matrix obtained by
positive diagonal scaling with row sums one and column sums $m/n$.
With these prescribed sums the scaled matrix is
unique~\cite{Sinkhorn1967}. Alternating row and column scaling
to those sums converges to it~\cite[Theorem 3.1 and Section 4]{ftIdel2016}.

The top-left entry $x=\operatorname{Sink}(A)_{11}$ can be irrational
even when $A$ has integer entries. We seek a polynomial that vanishes
at $x$ and whose coefficients can be computed from $A$. For a $2\times2$
input, the scaled matrix has $x$ on the diagonal and $1-x$ off the diagonal.
Diagonal scaling preserves the ratio of the diagonal and off-diagonal
products, so $a_{12}a_{21}x^2=a_{11}a_{22}(1-x)^2$. Solving for
$x\in(0,1)$ gives Nathanson's square-root formula~\cite[Theorem~3]{nmNathanson2020}.

In the symmetric $3\times3$ case, Ekhad and
Zeilberger~\cite[Theorem 1]{nmEkhadZeilberger2019} obtained a quartic
equation for $x$. Rowland and Wu~\cite[Theorems 3 and 6]{ftRowlandWu2025}
gave a degree-six equation for general $3\times3$ input and proposed
a coefficient formula for arbitrary sizes in \cite[Conjecture 2]{ftRowlandWu2025}.
We write the polynomial in $z$ as $\det K(A,z)$, and the conjecture
asserts that it vanishes at $z=x$.
We begin with a matrix $A$ that already has the prescribed row and
column sums, so $x=a_{11}$. These sums give a null vector of $K(A,x)$.
To return to the original matrix, we use the fact that scaling $A$
only multiplies the columns of $K(A,z)$ by positive factors.

\subsection{The exact coefficient formula}
Let $A=(a_{ij})\in\mathbb R_{>0}^{m\times n}$, with $m,n\ge1$,
and let $x=\operatorname{Sink}(A)_{11}$.
The coefficient formula uses pairs of trailing row and
column sets,
\[
 \mathcal D=\{(R,C):R\subseteq\{2,\ldots,m\},\ C\subseteq\{2,\ldots,n\},\ |R|=|C|\},
\]
and two quantities at each pair,
\[
 \Delta_{R,C}=\det A_{\{1\}\cup R,\{1\}\cup C},\qquad
 \Gamma_{R,C}=a_{11}\det A_{R,C}.
\]
For example, at $R=\{i\}$ and $C=\{j\}$ the two entries are
$a_{11}a_{ij}-a_{1j}a_{i1}$ and $a_{11}a_{ij}$. Each product
uses rows $1,i$ and columns $1,j$ once. In general the products
use the rows $\{1\}\cup R$ and columns $\{1\}\cup C$ once each,
so the two expressions acquire the same factor under scaling.

Keep the original row and column order in each minor and take empty
determinants and products to be one. To specify
the cofactor signs, write $p_s(U)$ for the position of $s$ in
the ordered set $U$, with the first position numbered one.
Write $R+s=R\cup\{s\}$ and $R-s=R\setminus\{s\}$, and likewise
for columns.

Index both axes of $H$ and both vectors $\Delta=(\Delta_{R,C})$,
$\Gamma=(\Gamma_{R,C})$ by $\mathcal D$ in the same order. On an
index subset $\mathcal S$, the principal submatrix $H_{\mathcal S}$
is Rowland and Wu's signed adjacency matrix. Its diagonal entry
depends only on the size of the minor,
\[
 H_{(R,C),(R,C)}=|R|(m+n)-mn
\]
The nonzero off-diagonal entries in row $(R,C)$ connect the minor to those
obtained by adding or deleting a row and a column together,
or replacing a single row or column. In an addition or
replacement, choose the new label from the trailing indices
outside the set it enters. This gives
\begin{align*}
 H_{(R,C),(R+s,C+t)}&=(-1)^{p_s(R+s)+p_t(C+t)}m,\\
 H_{(R,C),(R-s,C-t)}&=-(-1)^{p_s(R)+p_t(C)}n,\\
 H_{(R,C),(R,C-s+t)}&=(-1)^{p_s(C)+p_t(C-s+t)}m,\\
 H_{(R,C),(R-s+t,C)}&=(-1)^{p_s(R)+p_t(R-s+t)}n.
\end{align*}

\begin{problem}{NM-04}
For every $m,n\ge1$ and $A\in\R_{>0}^{m\times n}$, does the top-left
entry $x$ of the Sinkhorn limit satisfy the Rowland--Wu identity
\begin{equation}\label{eq:nm04-target}
 \sum_{\mathcal S\subseteq\mathcal D}\det(m^{-1}H_{\mathcal S})
 \prod_{(R,C)\in\mathcal S}\Delta_{R,C}
 \prod_{(R,C)\in\mathcal D\setminus\mathcal S}\Gamma_{R,C}\,x^{|\mathcal S|}=0?
\end{equation}
\end{problem}

\begin{theorem}[Rowland--Wu coefficient identity]\label{thm:nm04-main}
Equation~\eqref{eq:nm04-target} holds for every $A\in\mathbb R_{>0}^{m\times n}$,
even when some minors vanish.
\end{theorem}

To express the sum in \eqref{eq:nm04-target} as a determinant, set
\begin{equation}\label{eq:nm04-K}
 K(A,z)=\diag(\Gamma)+\frac z m H\diag(\Delta).
\end{equation}
By multilinearity, $\det K(A,z)$ has one term for each choice of
columns from the two matrices in \eqref{eq:nm04-K}. Let $\mathcal S$
mark the columns taken from the term containing $z$.
Outside $\mathcal S$ the only possible entry in each column
is $\Gamma_{R,C}$ on the diagonal. Expanding along such a column
takes out that factor and removes the row with the same index,
with positive cofactor sign. The remaining
matrix is $m^{-1}H_{\mathcal S}$ with column $(R,C)$ multiplied
by $z\Delta_{R,C}$. Its determinant and the removed factors give
exactly the $\mathcal S$-summand of \eqref{eq:nm04-target} with
$z$ in place of $x$. It remains to show $K(A,x)$ is singular.

\subsection{Four ways to differentiate a minor}
The four operations on index sets in $H$ also arise by
differentiating minors in four rank-one directions.
Let $T$ be a real $p\times q$ matrix with ordered row and column labels, and let $f=(f_{R,C})$ denote its array of square minors, where $f_{R,C}=\det T_{R,C}$ and $|R|=|C|=k$. In the application, the labels are the trailing row and column indices of $A$. Define the operators
\begin{align*}
 (\mathcal Lf)_{R,C}
 &=\sum_{s\in R,\ t\in C}(-1)^{p_s(R)+p_t(C)}f_{R-s,C-t},\\
 (\mathcal Uf)_{R,C}
 &=\sum_{s\notin R,\ t\notin C}(-1)^{p_s(R+s)+p_t(C+t)}f_{R+s,C+t},\\
 (\mathcal Cf)_{R,C}
 &=\sum_{s\in C,\ t\notin C}(-1)^{p_s(C)+p_t(C-s+t)}f_{R,C-s+t},\\
 (\mathcal Rf)_{R,C}
 &=\sum_{s\in R,\ t\notin R}(-1)^{p_s(R)+p_t(R-s+t)}f_{R-s+t,C}.
\end{align*}
Missing indices range over the full row or column label set, and empty sums vanish. Write $Df[Z]$ for the array of directional derivatives at $T$ in direction $Z\in\mathbb R^{p\times q}$. Here and below, $\mathbf1$ is a column vector of ones of the required dimension, including dimension zero. Set $s(T)=\mathbf1_p^TT\mathbf1_q$.

\begin{lemma}\label{lem:nm04-derivatives}
For every $T$, even when $T_{R,C}$ is singular,
\begin{align*}
 Df[\mathbf1\mathbf1^T]&=\mathcal Lf,\\
 Df[(T\mathbf1)\mathbf1^T]&=kf+\mathcal Cf,\\
 Df[\mathbf1(\mathbf1^TT)]&=kf+\mathcal Rf,\\
 Df[(T\mathbf1)(\mathbf1^TT)]&=s(T)f-\mathcal Uf.
\end{align*}
Each identity is understood at the coordinate $(R,C)$, with $k=|R|=|C|$.
\end{lemma}
\begin{proof}
Fix $(R,C)$ with $k>0$. The first derivative of the determinant
is a sum over choices of one column to perturb. In the all-ones
direction, expansion down that column gives precisely the cofactors
summed in $(\mathcal Lf)_{R,C}$.

For the second identity, the perturbation in each column is
$(T\mathbf1)_R=\sum_t T_{R,t}$. Consider replacing column $s\in C$.
The summand $t=s$ contributes $f_{R,C}$, and the other $t\in C$
give determinants with a repeated column. Thus the columns in
$C$ contribute $kf_{R,C}$ in total. For $t\notin C$, the new column
set is $C-s+t$. Moving $t$ into its ordered position gives the
sign in $(\mathcal Cf)_{R,C}$. Apply the same calculation to $T^T$
for the third identity.

For the fourth, set $V=T_{R,C}$ and use a bordered determinant. The classical adjugate $\adj(V)$ satisfies $V\adj(V)=\det(V)I_k$. The polynomial identity
\[
 T_{i,C}\adj(V)T_{R,j}
 =T_{ij}\det V-
 \det\begin{pmatrix}V&T_{R,j}\\T_{i,C}&T_{ij}\end{pmatrix}
\]
holds for all $V$. Sum over the row and column labels $i,j$ of $T$.
The left side is $(\mathbf1_p^TT)_C\adj(V)(T\mathbf1_q)_R$, the
cofactor formula for the derivative in the fourth direction. The
first term on the right sums to $s(T)f_{R,C}$. If $i\in R$ or $j\in C$,
the bordered determinant has a repeated row or column and vanishes.
For the remaining pairs, moving the last row and column to their
ordered positions contributes the sign $(-1)^{p_i(R+i)+p_j(C+j)}$.
These determinants give $(\mathcal Uf)_{R,C}$.

For $k=0$, the minor is constant one and all derivatives vanish. The fourth identity is $0=s(T)-\sum_{i,j}T_{ij}$; the first three have zero on both sides.
\end{proof}

\subsection{A null vector from the prescribed row and column sums}
Suppose $A$ already has row sums one and column sums $m/n$.
The minors through its top-left entry $x$ can be read from
the Schur complement of that entry. Write
\[
 A=\begin{pmatrix}x&r^T\\c&B\end{pmatrix},\qquad
 \rho=m/n,\qquad T=B-cr^T/x.
\]
The row and column sums give
\begin{equation}\label{eq:nm04-sums}
 T\mathbf1=\mathbf1-c/x,\qquad
 \mathbf1^TT=\rho(\mathbf1^T-r^T/x),\qquad
 s(T)=m-\rho/x.
\end{equation}
Define $u=\mathbf1-T\mathbf1$ and
$v^T=\mathbf1^T-\rho^{-1}\mathbf1^TT$.
By \eqref{eq:nm04-sums}, $u=c/x$ and $v=r/x$, so $B=T+xuv^T$.
The bordered minor formula gives $\Delta_{R,C}=xf_{R,C}$.
For the trailing minor $b_{R,C}=\det B_{R,C}$, each column
of $xuv^T$ in the determinant expansion is a multiple of
$u_R$ after restriction to $R$. Thus terms with two or more
of these columns vanish. Using no update column gives
$f_{R,C}$, and the single-column terms give the derivative
in direction $uv^T$ multiplied by $x$. With $b=(b_{R,C})$,
\[
 b=f+xDf[uv^T].
\]
The update direction expands into the four directions from \cref{lem:nm04-derivatives}:
\[
 uv^T=\mathbf1\mathbf1^T
 -(T\mathbf1)\mathbf1^T
 -\rho^{-1}\mathbf1(\mathbf1^TT)
 +\rho^{-1}(T\mathbf1)(\mathbf1^TT).
\]
Substituting \cref{lem:nm04-derivatives} and \eqref{eq:nm04-sums} into the rank-one update gives, at order $k$,
\[
 b=x\left[(n-k-k/\rho)f+\mathcal Lf-\mathcal Cf
              -\rho^{-1}\mathcal Rf-\rho^{-1}\mathcal Uf\right].
\]
Putting $\Delta=xf$ and collecting terms gives
\begin{equation}\label{eq:nm04-master}
 mb+[k(m+n)-mn]\Delta+n\mathcal U\Delta-m\mathcal L\Delta
       +m\mathcal C\Delta+n\mathcal R\Delta=0.
\end{equation}

To match \eqref{eq:nm04-master} with $H$, weight the columns of
$H\diag(\Delta)$ and divide each row by its own weight. An addition
goes from order $k$ to $k+1$ and has coefficient $m$ in $H$ and $n$
in \eqref{eq:nm04-master}, so the ratio of the weights must be $n/m$.
Take
\[
 w_{R,C}=\rho^{-|R|}=(n/m)^{|R|}.
\]
The reverse ratio $m/n$ converts deletion coefficients from $-n$ to
$-m$. An exchange preserves the order of the minor and hence its
coefficient. Applying these weights to \eqref{eq:nm04-master} gives
\[
 m\diag(b)w+H\diag(\Delta)w=0.
\]
Since $\Gamma=xb$, multiplication by $x/m$ proves
\begin{equation}\label{eq:nm04-null}
 K(A,x)w=0.
\end{equation}
When $m=1$ or $n=1$, the trailing minors are empty. The same identities
apply with the empty determinant and vector conventions.

\begin{proof}[Completion of \cref{thm:nm04-main}]
Write $A'=\operatorname{Sink}(A)=D_\lambda A D_\mu$.
The positive diagonal factors $D_\lambda$ and $D_\mu$ have diagonal
entries $\lambda_i$ and $\mu_j$, respectively. Both $\Gamma_{R,C}$
and $\Delta_{R,C}$ are multiplied by
\[
 \kappa_{R,C}=\lambda_1\mu_1\prod_{i\in R}\lambda_i\prod_{j\in C}\mu_j>0.
\]
With $\kappa=(\kappa_{R,C})$, this gives $K(A',z)=K(A,z)\diag(\kappa)$ for every $z$. Applying \eqref{eq:nm04-null} at the top-left entry $x$ of $A'$ gives $K(A,x)\diag(\kappa)w=0$. Since $\diag(\kappa)w$ is nonzero, $\det K(A,x)=0$, proving \eqref{eq:nm04-target}.
\end{proof}

\section{Disconnected spaces of minimal PSD factorizations (PF-02)}\label{sec:pf02}

In a positive semidefinite (PSD) factorization, the entry in row $i$
and column $j$ is the trace product of a row factor $A_i$ and a column
factor $B_j$. For a polytope's slack matrix, these factors encode a
semidefinite lift, a projection of an affine slice of a PSD cone.
The theorem of Gouveia, Parrilo and Thomas~\cite{nmGouveiaParriloThomas2013}
relates the size of the factors to the size of the lift. For the
theory of PSD rank, see Fawzi et al.~\cite{ftFawziEtAl2015}.

\subsection{The factorization space and an explicit matrix}One matrix can have several PSD factorizations of the smallest
size. A change of basis by congruence preserves the trace products,
so we identify the factorizations in each orbit under this action.
Let $\mathbb S^k$ denote the space of real symmetric $k\times k$ matrices, and let $\mathbb S_+^k=\{X\in\mathbb S^k:X\succeq0\}$ be its PSD cone. For an entrywise-nonnegative matrix $M\in\mathbb R^{p\times q}$, its real PSD rank is the least positive integer $k$ for which
\[
 M_{ij}=\tr(A_iB_j),\qquad A_i,B_j\in\mathbb S_+^k.
\]
Write $A=(A_1,\ldots,A_p)$ and $B=(B_1,\ldots,B_q)$. Let
$\mathcal F_k(M)$ be the set of pairs $(A,B)$ satisfying these
factorization equations, with the topology inherited from Euclidean
space. An invertible $S\in\operatorname{GL}(k,\mathbb R)$ acts by
the following congruences, where $S^{-T}=(S^{-1})^T$.
\begin{equation}\label{eq:pf02-action}
 A_i\longmapsto S^TA_iS,\qquad
 B_j\longmapsto S^{-1}B_jS^{-T},\qquad S\in \operatorname{GL}(k,\mathbb R).
\end{equation}
The transformed product is $S^TA_iB_jS^{-T}$, whose trace is still
$\tr(A_iB_j)$. We identify pairs related by this action and give
$\mathcal F_k(M)/\operatorname{GL}(k,\mathbb R)$ the quotient topology.

\begin{problem}{PF-02}
Let $k\ge3$ and suppose the real entrywise nonnegative matrix $M$ has
real PSD rank $k$ and ordinary rank $k(k+1)/2$. Is
$\mathcal F_k(M)/\operatorname{GL}(k,\R)$ connected for every such
$M$~\cite[Problem~9.4]{ftFawziEtAl2015}?
\end{problem}

For real factors of size two and ordinary rank three, Fawzi et al.\ proved connectedness by identifying the orbits with ellipses nested between two polyhedra~\cite[Proposition 7.9]{ftFawziEtAl2015}.
Connectedness fails for every factor size at least three. Our first example
is a strictly positive integer matrix of ordinary rank six and PSD rank
three. In every minimal factorization, its six row factors form an ordered
basis of the symmetric $3\times3$ matrices. A congruence preserves the
orientation of this basis, and we find two factorizations with opposite orientations.

For the example, perturb $2I_3$ along the six symmetric coordinate
directions. Denote the matrix units by $E_{ij}$ and take
\[
 A_i=2I_3+2E_{ii}\quad(1\le i\le3),\qquad
 A_4=2I_3+E_{12}+E_{21},
\]
\[
 A_5=2I_3+E_{13}+E_{31},\qquad
 A_6=2I_3+E_{23}+E_{32},\qquad B_i=A_i.
\]
Their trace products $M_{ij}=\tr(A_iA_j)$ give the matrix
\begin{equation}\label{eq:pf02-M}
 M=\begin{pmatrix}
 24&20&20&16&16&16\\
 20&24&20&16&16&16\\
 20&20&24&16&16&16\\
 16&16&16&14&12&12\\
 16&16&16&12&14&12\\
 16&16&16&12&12&14
 \end{pmatrix}.
\end{equation}
To obtain the second factorization, conjugate only $A_4$ and
$B_4$ by $\diag(-1,1,1)$. This changes the signs of their $(1,2)$
and $(2,1)$ entries and preserves positive definiteness and
$\tr(A_4^2)$. All the other factors have zeros in those two
positions, so their trace products with the changed factor
are unaffected. The new factors therefore give the same $M$.

\begin{theorem}\label{thm:pf02-basic}
The matrix \eqref{eq:pf02-M} has ordinary rank six and real PSD rank three, and $\mathcal F_3(M)/\operatorname{GL}(3,\mathbb R)$ is disconnected.
\end{theorem}

\subsection{An orientation defined on the entire fiber}
Fix coordinates on $\mathbb S^3$ in the order
\[
 (X_{11},X_{22},X_{33},X_{12},X_{13},X_{23}).
\]
Let the coordinate row of $A_i$ be row $i$ of $U_A$. The
upper-right block vanishes, since the first three factors
are diagonal. Their diagonal coordinates give $2I_3+2ee^T$
for $e=(1,1,1)^T$. The last three factors give the block $I_3$
on the off-diagonal coordinates, so $U_A$ is invertible with
positive determinant. In the trace pairing, each off-diagonal
coordinate counts twice, so $M=U_AGU_A^T$ for
$G=\diag(1,1,1,2,2,2)$, and $\rank M=6$. The reflected
factorization negates one coordinate column of $U_A$ and
reverses the sign of its determinant.

In any size-$k$ real PSD factorization, $M$ pairs the row and column
factors through the space $\mathbb S^k$. Hence
\begin{equation}\label{eq:pf02-rank-bound}
 \rank M\le\dim\mathbb S^k=k(k+1)/2.
\end{equation}
For the rank-six matrix \eqref{eq:pf02-M}, this makes size three minimal.

\begin{lemma}\label{lem:pf02-congruence}
For $S\in\mathbb R^{k\times k}$, the determinant of the linear map $C_S:X\mapsto S^TXS$ on $\mathbb S^k$ equals $(\det S)^{k+1}$.
\end{lemma}
\begin{proof}
First take $S=\diag(s_1,\ldots,s_k)$. The diagonal matrix units are multiplied by $s_i^2$, and the symmetric off-diagonal units by $s_is_j$. In their product, each $s_i$ occurs with exponent $2+(k-1)=k+1$, so the determinant is $\prod_i s_i^{k+1}$. Changing basis gives the same formula for a real diagonalizable $S$. Both sides are polynomials in the entries of $S$. Equality on the nonempty open set of matrices with distinct real eigenvalues proves the identity for every real $S$.
\end{proof}

\begin{proof}[Proof of \cref{thm:pf02-basic}]
In any size-three factorization, a dependence among the six row
factors would give the same dependence among the rows of $M$.
Thus the row factors always form a basis, and
\[
 \varepsilon(A,B)=\operatorname{sgn}\det U_A\in\{-1,+1\}
\]
is continuous on the factorization space. A congruence multiplies
the determinant by $(\det S)^4>0$ by \cref{lem:pf02-congruence},
so the sign is also continuous on the quotient. The two explicit
factorizations take opposite values, separating the quotient into
two nonempty open-and-closed sets.
\end{proof}

\subsection{Extending the example to larger factor sizes}
When $k$ is even, a congruence can reverse the orientation of $\mathbb S^k$.
Instead we embed the size-three example in a larger one. The construction
lets us recover its three-dimensional subspace from every minimal factorization,
so we can still use the orientation argument on that subspace.

\begin{theorem}\label{thm:pf02-all}
For every $k\ge3$, there is a nonnegative integer matrix of ordinary rank
$k(k+1)/2$ and real PSD rank $k$ whose space of minimum-size PSD factorizations
is disconnected modulo congruence. There is also a strictly positive
rational matrix with these properties for each such $k$.
\end{theorem}
\begin{proof}[Nonnegative integer construction]
We have already dealt with $k=3$. Write $k=3+r$, where $r\ge1$, and set $d=r(r+1)/2$ and $N=k(k+1)/2$. Take the following PSD matrices $F_1,\ldots,F_N$ in the order listed:
\begin{enumerate}
\item the six size-three factors above, embedded in the leading block;
\item the $d$ basis factors $2I_r+2E_{ii}$, for $1\le i\le r$, and $2I_r+E_{ij}+E_{ji}$, for $1\le i<j\le r$, embedded in the trailing block;
\item the $3r$ full-size factors $2I_k+E_{ij}+E_{ji}$ with $1\le i\le3<j\le k$.
\end{enumerate}
Read the off-diagonal rectangular block of a vanishing linear
combination of the $F_i$. Each position $(i,j)$ with $i\le3<j$
picks out the coefficient of just $2I_k+E_{ij}+E_{ji}$, so all
coefficients from the third group vanish. On the two diagonal
blocks the remaining factors are bases, forcing their
coefficients to vanish as well. Consequently the trace Gram
matrix $(M_0)_{ij}=\tr(F_iF_j)$ has rank $N$. Its entries are
nonnegative integers and its PSD rank is $k$ by \eqref{eq:pf02-rank-bound}.

Let $\mathcal T$ and $\mathcal D_0$ be the indices of the first and second groups. The block $M_0[\mathcal T,\mathcal D_0]$ is zero because these groups occupy disjoint diagonal blocks. This is the block that lets us recover the support. Take any size-$k$ factorization $(\widetilde A,\widetilde B)$ of $M_0$ and set
\[
 Q=\sum_{j\in\mathcal D_0}\widetilde B_j.
\]
For PSD matrices, $\tr(AB)=\|A^{1/2}B^{1/2}\|_{\mathrm F}^2$, so zero trace product means orthogonal ranges. The six leading row factors therefore have ranges in $U=\ker Q$.

The block $M_0[\mathcal T,\mathcal T]$ has rank six, so the leading
factors are independent and $\dim U\ge3$. The column factors of the
rank-$d$ block $M_0[\mathcal D_0,\mathcal D_0]$ are PSD summands of
$Q$, so their ranges lie in $\operatorname{im}Q$. Writing $s=\rank Q$,
we have $d\le s(s+1)/2$, which gives $s\ge r$ and $\dim U\le3$.

So $U$ has dimension three, and the six factors restricted to it form an ordered basis of $\mathbb S(U)$, the symmetric operators on $U$. To assign the six factors a sign, choose an ordered orthonormal basis of $U$ and use the same coordinates as in the size-three example. By \cref{lem:pf02-congruence}, changing this basis acts on the symmetric matrices with determinant a positive fourth power. The sign is thus independent of the choice. Since $Q$ has constant rank $r$, the orthonormal basis can be chosen continuously near any factorization. The coordinate determinant never vanishes, so the sign is continuous.

Under \eqref{eq:pf02-action}, $Q$ becomes $S^{-1}QS^{-T}$ and its kernel becomes $S^TU$. The induced map on the symmetric matrices is a size-three congruence, so it preserves the sign.

Reflecting the $(1,2)$ coordinate changes just one factor in the first group. As before, the reflected factor is PSD and the trace products are unchanged. The two factorizations have opposite signs on the leading support. The sign descends to the quotient and takes both values, so the quotient is disconnected.
\end{proof}

\begin{proof}[Strictly positive rational construction]
The example for $k=3$ is already strictly positive. For $k>3$, add a small multiple of the identity to each factor:
\[
 F_i(\eta)=F_i+\eta I_k,\qquad
 M(\eta)_{ij}=\tr(F_i(\eta)F_j(\eta)).
\]
The factors are positive definite for $\eta>0$ and remain independent for small $\eta$. Thus $M(\eta)$ is strictly positive, with ordinary rank $N$ and PSD rank $k$. The reflection fixes $I_k$ and gives a second factorization for each $\eta$.

Adding $\eta I_k$ removes the zero block used to recover the support.
We now use the subspace for the three smallest eigenvalues of a sum of
column factors.
To control this space over all size-$k$ factorizations, first normalize
the row factors. Their sum $S_A=\sum_i A_i$ is positive definite. A
vector in its kernel would be in every $\ker A_i$, so all row factors
would be supported on a proper subspace, contradicting rank $N$. Apply
the congruence
\[
 A_i'=S_A^{-1/2}A_iS_A^{-1/2},\qquad
 B_j'=S_A^{1/2}B_jS_A^{1/2}
\]
to make the row factors sum to $I_k$. They are PSD and
bounded above by $I_k$. Each $B_j'$ is PSD with trace equal
to column sum $j$ of $M(\eta)$. Thus all normalized factors
are bounded when $\eta$ lies in a small closed interval
starting at zero. The factorization
equations and PSD conditions are closed, giving compactness
of the normalized factors together with $\eta$.

At $\eta=0$, the support is the three-dimensional kernel of
$Q'=\sum_{j\in\mathcal D_0}B_j'$, separated by a spectral gap
from its positive eigenvalues. Let $U'$ be the subspace of the three lowest
eigenvalues when that gap is open, and let $P_{U'}$ be its
orthogonal projector. At every normalized factorization
at zero, the gap is positive and the compressions
$P_{U'}A_i'|_{U'}$, $i\in\mathcal T$, are independent.
Both are open conditions, by continuity of the eigenvalues
and of the projector onto the separated subspace. A sequence
of failures at $\eta\to0$ would have a convergent subsequence
of normalized factors with limit at zero. Near that limit the
gap and independence hold, a contradiction.

The determinant of the six ordered compressions has a continuous
sign, as in the integer construction. For two factorizations related
by $S$, the congruence between their normalized factors is given by
the orthogonal matrix
\[
 O=S_A^{1/2}S(S^TS_AS)^{-1/2}.
\]
Under $X\mapsto O^TXO$, the spectral subspace becomes $O^TU'$,
and the compressions change by a size-three congruence. This
preserves the sign of their coordinate determinant. Since
normalization is continuous, the sign gives a continuous
function on the quotient.

The original family and its reflection have opposite
orientations at $\eta=0$. Their signs persist for small
positive $\eta$, by continuity and the nonvanishing of both
determinants. Choose a rational $\eta$ in that range.
\end{proof}

\section{Stable patterns on random cubic graphs (MD-06)}\label{sec:md}
Let $G=(V,E)$ be a finite simple graph. For phases $\theta\in\mathbb T^V$, with $\mathbb T=\R/(2\pi\mathbb Z)$, the homogeneous attractive Kuramoto energy is
\begin{equation}\label{md-energy}
 \mathcal E_G(\theta)=\sum_{\{u,v\}\in E}\bigl(1-\cos(\theta_u-\theta_v)\bigr).
\end{equation}
The phases evolve by $\dot\theta=-\nabla\mathcal E_G(\theta)$
and are synchronized when every pair agrees modulo $2\pi$.
In the sense of \cite[Definition~2.1]{Randomstrasse2025}, a graph is
globally synchronizing if every local minimum of \eqref{md-energy}
on the phase torus is synchronized. The energy is constant under
common phase shifts along the all-ones vector $\one\in\R^V$.
We therefore test strict local minimality modulo these shifts
by restricting to $\one^\perp$ with the Euclidean norm and inner product.

Abdalla, Bandeira, Kassabov, Souza, Strogatz and Townsend~\cite{Abdalla2026} proved global synchronization for sufficiently strong expanders, including typical regular graphs of sufficiently large degree. DeVille and Ermentrout~\cite{DeVille2016} found stable nonsynchronized patterns in cubic graphs, and their simulations pointed to the existence of such a pattern with high probability.

\begin{problem}{MD-06}
For each even $n\ge4$, let $G_n$ be uniformly distributed over the
labeled simple $3$-regular graphs on $\{1,\ldots,n\}$. As $n$ tends
to infinity through even integers, does the probability that every
local minimum of $\mathcal E_{G_n}$ on the phase torus is synchronized
tend to one~\cite[Conjecture~3]{Randomstrasse2025}?
\end{problem}

The following theorem shows that this does not hold in general.

\begin{theorem}\label{md-main}
Let $G_n$ be uniformly distributed over the simple $3$-regular graphs on $\{1,\ldots,n\}$, with $n\geq4$ even. Then
\[
 \lim_{\substack{n\to\infty\\ n\text{ even}}}
 \mathbb P\bigl(\text{every local minimum of }\mathcal E_{G_n}
                       \text{ is synchronized}\bigr)=0.
\]
Indeed, with probability tending to one, there is a nonsynchronized critical point $\theta_*$ at which every edge cosine exceeds $1/32$ and
\[
 z^T\nabla^2\mathcal E_{G_n}(\theta_*)z
       \geq\frac1{320}\norm{z}_2^2\qquad(z\in\one^\perp).
\]
\end{theorem}

At each cycle vertex we balance the force from its tree with one of the two cycle edges. To continue the balance down the tree, each of the two child edges carries half the force on the parent edge, with the opposite sign. So only the leaves contribute to the gradient. We take the trees deep enough for this residual to be small, then use the Laplacian gap to correct the phases.

\subsection{Correcting an approximate equilibrium}
For real lifts $\phi\in\R^V$ of the phases and vectors $z\in\R^V$,
differentiating \eqref{md-energy} gives the following formulas, where
$w\sim v$ denotes adjacency.
\begin{align}
 (\nabla\mathcal E_G(\phi))_v
   &=\sum_{w\sim v}\sin(\phi_v-\phi_w),\label{md-gradient}\\
 z^T\nabla^2\mathcal E_G(\phi)z
   &=\sum_{\{u,v\}\in E}\cos(\phi_u-\phi_v)(z_u-z_v)^2.
   \label{md-hessian}
\end{align}
Write $L_G$ for the combinatorial Laplacian, the degree matrix
minus the adjacency matrix. For $|V|\ge2$ its spectral gap
$\lambda_{\mathrm{gap}}(L_G)$ is its second-smallest eigenvalue,
counting multiplicity. The Laplacian has quadratic form $\sum_{\{u,v\}\in E}(z_u-z_v)^2$.
Thus if the cosines in \eqref{md-hessian} are all at least $c>0$,
the Hessian is bounded below by $cL_G$, which on $\one^\perp$ is
bounded below by $c\lambda_{\mathrm{gap}}(L_G)I$.

\begin{lemma}\label{md-correction}
Suppose that $G$ is connected with at least two vertices, $\lambda_{\mathrm{gap}}(L_G)\geq\gamma>0$, and $\phi\in\R^V$ satisfies
\[
 \cos(\phi_u-\phi_v)\geq c>0\quad(\{u,v\}\in E),
 \qquad \norm{\nabla\mathcal E_G(\phi)}_2<\frac{c^2\gamma}{4\sqrt2}.
\]
There is $h_*\in\one^\perp$, with $\norm{h_*}_2<c/(2\sqrt2)$, such that $\phi+h_*$ is an exact local minimum. At this point every edge cosine is at least $c/2$, and the Hessian on $\one^\perp$ is bounded below by $(c\gamma/2)I$.
\end{lemma}
\begin{proof}
Take $\rho=c/(2\sqrt2)$ and let
$B=\{h\in\one^\perp:\norm{h}_2\le\rho\}$. The energy
$\mathcal E_G(\phi+h)$ attains a minimum on this compact ball.
For $h\in B$, an edge difference changes by $h_u-h_v$, of
magnitude at most $\sqrt2\norm{h}_2\le c/2$. Because cosine
is $1$-Lipschitz, every weight in \eqref{md-hessian} stays
at least $c/2$. Together with the gap bound this gives
\[
 z^T\nabla^2\mathcal E_G(\phi+h)z
 \geq\frac c2 z^TL_Gz\geq\mu\norm{z}_2^2,
 \qquad \mu=\frac{c\gamma}{2},\quad z\in\one^\perp.
\]
For $\norm{h}_2=\rho$, we calculate the derivative in direction
$h$ at $\phi+h$. The Hessian bound applies along the entire
segment from $\phi$ to $\phi+h$, and the gradient hypothesis
at $\phi$ gives $\norm{\nabla\mathcal E_G(\phi)}_2<\mu\rho$.
Integrating along the segment yields
\begin{align*}
 \langle h,\nabla\mathcal E_G(\phi+h)\rangle
 &=\langle h,\nabla\mathcal E_G(\phi)\rangle
   +\int_0^1h^T\nabla^2\mathcal E_G(\phi+th)h\,\mathrm{d}t\\
 &\geq-\rho\norm{\nabla\mathcal E_G(\phi)}_2+\mu\rho^2>0.
\end{align*}
Thus a short step in direction $-h$ from any boundary point
lowers the energy, placing the minimizer $h_*$ in the interior
of $B$. The Hessian bound gives a strict minimum in the directions
in $\one^\perp$. Common shifts leave the energy constant, so
$\phi+h_*$ is a local minimum on the torus.
\end{proof}

\subsection{Phases on a cycle and its trees}
Call the radius-$R$ neighborhood of a cycle \emph{clean} if its induced subgraph has no other cycle. In a cubic graph, each cycle vertex then has one outward child, and vertices at depths $1,\ldots,R-1$ have two. For a cycle of length $\ell$, there are $\ell2^{j-1}$ vertices at depth $j$, for $1\leq j\leq R$.

Take $R\geq4$ and a clean cycle of length $\ell$ divisible by four.
We want one of the two cycle edges at each root to balance its tree
edge. Assign signs in pairs so that each root has one cycle neighbor
of each sign,
\[
 (+,+,-,-,+,+,-,-,\ldots).
\]
In a tree with root sign $\varepsilon$, put $\phi_v=\varepsilon a_j$
at depth $0\leq j\leq R$, with $a_R=0$, and set all other phases to
zero. For $0<t<1$, set $a_j-a_{j+1}=\arcsin(t/2^j)$. This halves the
sine at each branch and keeps the edge cosines positive. At a positive
root the equal-sign edge contributes zero. The other cycle edge
has sine $\sin(2a_0)$, which must be $-t$. Taking
$2a_0=2\pi-\arcsin t$ gives that sine with a positive cosine.
Since $a_0$ is the sum of the differences down the tree, the root
equation is $F_R(t)=2\pi$, where
\[
 F_R(t)=3\arcsin t+2\sum_{j=1}^{R-1}\arcsin(t/2^j),\qquad 0\leq t\leq1.
\]

We need a bound on the root of $F_R(t)=2\pi$ that holds for every
$R\ge4$. Extra levels only add positive terms to $F_R$, so we
test $R=4$ at $t_*=\sqrt{255}/16$. For the first term, use
$\arcsin t_*=\pi/2-\arcsin(1/16)$ and
$\arcsin(1/16)<1/\sqrt{255}<1/15$. For the remaining terms, use
$\arcsin x\geq x+x^3/6$. With $t_*>99/100$, $t_*^3>97/100$,
and $151/96>11/7>\pi/2$, this gives
\[
 F_4(t_*)>\frac{3\pi}{2}-\frac15+\frac{7t_*}{4}+\frac{t_*^3}{24}
          >\frac{3\pi}{2}+\frac{151}{96}>2\pi.
\]
As $F_R$ increases strictly from $F_R(0)=0$, its unique solution
satisfies $0<t_R\leq t_4<t_*$. In particular the cosine
$\sqrt{1-t_R^2}$ is bounded below by $c_0:=\sqrt{1-t_4^2}>1/16$,
uniformly in $R$.

Use $t_R$ to set the phase differences on a tree, and sum them from
the leaves back to the root to define the profile
\begin{equation}\label{md-profile}
 \delta_j=\arcsin(t_R/2^j),\qquad
 a_j=\sum_{s=j}^{R-1}\delta_s\quad(0\leq j<R),\qquad a_R=0.
\end{equation}
The root equation gives
\begin{equation}\label{md-rootbalance}
 2a_0+\delta_0=2\pi,\qquad
 \sin(2a_0)=-t_R,\qquad \cos(2a_0)=\sqrt{1-t_R^2}.
\end{equation}
With the interior forces balanced, the residual comes from the
$\ell 2^{R-1}$ leaves at depth $R$, each with force
$-\varepsilon t_R/2^{R-1}$ from its parent. Thus
\begin{equation}\label{md-residual}
 \norm{\nabla\mathcal E_G(\phi)}_2^2
   =\ell2^{R-1}(t_R/2^{R-1})^2
   =\ell t_R^2 2^{1-R}.
\end{equation}
The tree-edge cosines $\cos\delta_j$ and the opposite-sign cycle-edge cosine $\cos(2a_0)$ are at least $c_0$; all other edge cosines are one.

\begin{proposition}\label{md-obstruction}
Fix $\gamma>0$ and a positive multiple of four $\ell$. For an integer $R$, if
\begin{equation}\label{md-radius}
 R\geq4,\qquad 2^{R-1}>\frac{32\ell}{c_0^4\gamma^2},
\end{equation}
then every connected cubic graph with $\lambda_{\mathrm{gap}}(L_G)\geq\gamma$ and a clean radius-$R$ neighborhood of an $\ell$-cycle has a nonsynchronized local minimum. Every edge cosine at this minimum is at least $c_0/2$, and its Hessian on $\one^\perp$ is at least $(c_0\gamma/2)I$.
\end{proposition}
\begin{proof}
The residual in \eqref{md-residual} has squared norm at most
$\ell2^{1-R}$. By \eqref{md-radius} this is smaller than
$c_0^4\gamma^2/32$, so \cref{md-correction} applies.

For an opposite-sign cycle edge, \eqref{md-rootbalance}
puts its phase difference at distance $\delta_0$ from
$2\pi\mathbb Z$. Convexity of $\arcsin$ gives
$\delta_j\le2^{-j}\delta_0$, so $2\pi=F_R(t_R)<5\delta_0$.
The correction changes the difference by at most $c_0/2\le1/2$.
Since $\delta_0>2\pi/5$, the two vertices still have different phases.
\end{proof}

\subsection{The gap and the clean cycle}
For a cubic graph the Laplacian gap is three minus the
second-largest adjacency eigenvalue. Bordenave's form of
Friedman's theorem~\cite[Theorem~1]{Bordenave2020} shows that
for every $\delta>0$, this eigenvalue is at most $2\sqrt2+\delta$
with high probability.
Since $3-2\sqrt2>1/10$, this gives
\begin{equation}\label{md-gap}
 \mathbb P\bigl(\lambda_{\mathrm{gap}}(L_{G_n})\geq1/10\bigr)\longrightarrow1,
\end{equation}
so $G_n$ is connected with high probability. For each fixed
$\ell\ge3$, let $C_\ell(G_n)$ count the $\ell$-cycles.
Johnson's Poisson approximation~\cite[Theorem~11]{Johnson2015}
gives a limit of mean $\mu_\ell=2^\ell/(2\ell)$ and in particular
\begin{equation}\label{md-poisson}
 \mathbb P(C_\ell(G_n)=0)\longrightarrow e^{-\mu_\ell}.
\end{equation}

A nonclean neighborhood contains the original $\ell$-cycle and
a second cycle. It is connected, so it has at least two edges in
excess of a spanning tree. The first layer around the cycle has
at most $\ell$ vertices and each layer has at most twice as many
as the last. This puts at most $M=\ell2^R$ vertices in the
neighborhood. Prune the leaves one at a time to obtain a subgraph
$H$ of minimum degree two. Each step preserves the difference
between the numbers of edges and vertices, so
\[
 v(H)\leq M,\qquad e(H)-v(H)\geq1.
\]
For any specified labeled copy of $H$, Johnson's subgraph bound~\cite[Proposition~1(a)]{Johnson2015} gives
\[
 \mathbb P(H\subseteq G_n)\leq c_1 2^{e(H)}n^{-e(H)}
\]
for an absolute constant $c_1$ and large $n$. Summing over
the at most $n^{v(H)}$ placements of each type gives a bound
of $O(n^{v(H)-e(H)})=O(n^{-1})$. The number of types is finite
for fixed $\ell$ and $R$, so
\begin{equation}\label{md-clean}
 \mathbb P(\text{some $\ell$-cycle has a nonclean radius-$R$ neighborhood})
     =O_{\ell,R}(n^{-1}).
\end{equation}

\begin{proof}[Proof of \cref{md-main}]
Fix a positive multiple of four $\ell$ and choose $R$ from
\eqref{md-radius} with $\gamma=1/10$. If the gap is at least $1/10$
and there is a clean $\ell$-cycle, \cref{md-obstruction} gives the
required minimum. Since $c_0>1/16$, its edge cosines exceed $1/32$
and its Hessian on $\one^\perp$ is bounded below by $(1/320)I$.
For the construction to fail, the gap must be too small, the graph
must have no $\ell$-cycle, or some $\ell$-cycle must be nonclean.
For fixed $\ell,R$, the first and third probabilities tend to zero
by \eqref{md-gap} and \eqref{md-clean}. Using \eqref{md-poisson}
for the probability of no $\ell$-cycle, the union bound gives
\[
 \limsup_{\substack{n\to\infty\\ n\text{ even}}}
 \mathbb P(\text{no such minimum exists})
 \leq \exp\!\left(-\frac{2^\ell}{2\ell}\right).
\]
The limsup on the left is independent of $\ell$. Taking $\ell$
to infinity through the multiples of four makes the upper bound
tend to zero.
\end{proof}

\section{Interval determinants and solution hulls (\mbox{IV-02}, \mbox{IV-04})}\label{sec:intervals}
We can compute the determinant of a given tridiagonal matrix by
a recurrence. Let each entry vary independently over a given
closed interval, defining an interval matrix
$\mathcal A\subset\R^{n\times n}$. All supplied interval endpoints
are binary-encoded rationals.
The range is a closed interval because the determinant is continuous
on the compact, connected set $\mathcal A$.

\begin{problem}{IV-02}
Given a tridiagonal interval matrix $\mathcal A$ of order $n\ge2$,
which may contain singular matrices, can a deterministic algorithm
compute the two exact rational endpoints
\[
 d_-=\min_{A\in\mathcal A}\det A,\qquad d_+=\max_{A\in\mathcal A}\det A
\]
in time polynomial in the total binary input length?
\end{problem}

For $A\in\mathcal A$ and $b$ in an interval
vector $\mathbf b\subset\R^n$, their possible solutions form
the united solution set
\[
 \Sigma(\mathcal A,\mathbf b)=\{x\in\R^n:Ax=b\text{ for some }A\in\mathcal A,\ b\in\mathbf b\},
\]
whose smallest axis-aligned enclosing box is its interval hull.
Its bounds in coordinate $i$ are the infimum and supremum of
$x_i$ over $\Sigma(\mathcal A,\mathbf b)$.

\begin{problem}{IV-04}
Given a tridiagonal interval system $(\mathcal A,\mathbf b)$ of
order $n\ge1$ whose matrix family may contain singular members,
can a deterministic algorithm compute the exact interval hull
of $\Sigma(\mathcal A,\mathbf b)$
in time polynomial in the total binary input length? It must
report an empty solution set when appropriate and otherwise
return each finite endpoint as an exact rational and each infinite
endpoint explicitly.
\end{problem}

Both computations are NP-hard even for \emph{regular} families,
where every $A\in\mathcal A$ is nonsingular.

For tridiagonal families, Bar-On, Codenotti and Leoncini~\cite{BarOnCodenottiLeoncini1996} gave a linear-time test for regularity. Hor\'a\v cek, Hlad\'ik and Mat\v ejka~\cite[Section~5.4, Proposition~5.6]{IntervalDet2018} computed the determinant range in polynomial time for tridiagonal interval $H$-matrices and asked about general tridiagonal families. For solution hulls, Hor\'a\v cek, Hlad\'ik and \v Cern\'y~\cite[Section~1.4.3, Theorem~9]{IntervalComplexity2017} treated regular bidiagonal systems and left the tridiagonal case open.

We reduce from PARTITION~\cite{Karp1972}, the NP-complete decision
problem that asks whether a list of positive integer weights
$w_1,\ldots,w_m$ can be split into two groups with equal total weight.
In the equivalent sign formulation, the answer is yes exactly when
there are signs $\sigma_i\in\{-1,1\}$ with $\sum_{i=1}^m\sigma_iw_i=0$.
The reduction associates one planar rotation with each weight and
uses a choice of interval endpoint between rotations to keep or
reverse the accumulated angle.

\subsection{Rational rotations as continuants}
Let $T$ be tridiagonal of order $N$, with diagonal entries $a_j$,
superdiagonal ones and subdiagonal entries $\beta_j$ for $j\ge2$.
Let $p_j$ be the continuant, the determinant of the leading
$j\times j$ block. Expansion along its last row gives
\[
 p_{-1}=0,\quad p_0=1,\quad p_j=a_jp_{j-1}-\beta_jp_{j-2}.
\]
The first step uses $p_{-1}=0$, so $\beta_1$ can be arbitrary.
To get the next pair $(p_j,p_{j-1})^T$, multiply $(p_{j-1},p_{j-2})^T$
by
\[
 M(a,\beta)=\begin{pmatrix}a&-\beta\\1&0\end{pmatrix}.
\]
Starting from $(p_0,p_{-1})^T=(1,0)^T$, we read the determinant
$p_N$ from the first coordinate of
\[
 M(a_N,\beta_N)\cdots M(a_1,\beta_1)(1,0)^T.
\]

Use an angle $\theta_i=2\arctan(t_i)$ for weight $w_i$, where $t_i$
is a small rational multiple of $w_i$. Then the sine and cosine
are rational and the angle is close to $2t_i$. Specifically,
\[
 W=\sum_iw_i,\quad\delta=\frac1{10W^2},\quad t_i=\delta w_i,\quad
 c_i=\frac{1-t_i^2}{1+t_i^2},\quad s_i=\frac{2t_i}{1+t_i^2}.
\]
Here $c_i^2+s_i^2=1$ and $c_i>0$. The rotation matrix
for $\theta_i$ factors into two recurrence steps,
\[
 R_i=\begin{pmatrix}c_i&-s_i\\s_i&c_i\end{pmatrix}
 =M(-s_i/c_i,-1/c_i)M(s_i,-c_i)
\]
Between two rotations we can insert two steps with product
\[
 M(0,-1)M(0,q)=\diag(1,-q),\qquad q\in[-1,1].
\]
Choosing $q=-1$ makes the pair the identity and choosing
$q=1$ makes it a reflection. Use a separate interval entry
$q_i$ between rotations $i$ and $i+1$. Reading the steps
in order gives the $N=4m-2$ rows of $T(q)$,
\[
\begin{array}{c|cc|l}
 j&a_j&\beta_j&\text{range}\\\hline
 4i-3&s_i&-c_i&i=1,\ldots,m\\
 4i-2&-s_i/c_i&-1/c_i&i=1,\ldots,m\\
 4i-1&0&q_i&i=1,\ldots,m-1,\quad q_i\in[-1,1]\\
 4i&0&-1&i=1,\ldots,m-1.
\end{array}
\]
At a vertex of the box, the vector stays on the unit circle.
Let $\phi_i$ be its angle after applying rotation $i$.
The next pair sends $\phi_i$ to $\epsilon_i\phi_i$ for
$\epsilon_i=-q_i\in\{-1,1\}$, and then rotation $i+1$ adds
$\theta_{i+1}$. Thus
\[
 \phi_1=\theta_1,\qquad \phi_{i+1}=\theta_{i+1}+\epsilon_i\phi_i.
\]
The coefficient of $\theta_i$ in $\phi_m$ is the product
$\epsilon_i\cdots\epsilon_{m-1}$, with coefficient one for
$\theta_m$. Since $\det T(q)$ is the first coordinate of
the final vector, its values at the vertices are
\begin{equation}\label{eq:iv-cos}
 \cos\Big(\sum_{i=1}^m\sigma_i\theta_i\Big),
 \qquad \sigma_i\in\{-1,1\},\quad\sigma_m=1.
\end{equation}
Conversely, every pattern in \eqref{eq:iv-cos} comes from the
choices $\epsilon_i=\sigma_i\sigma_{i+1}$. It suffices to look for
a partition with $\sigma_m=1$ because reversing all the signs
preserves a zero sum.

\subsection{Regularity and separation}
First consider a vertex of the parameter box. The angle in
\eqref{eq:iv-cos} is at most $\sum_i\theta_i\le2\delta W=1/(5W)\le1/5$
in absolute value, so its cosine is positive. For arbitrary $q$,
hold all variables but one $q_i$ fixed. The determinant depends
affinely on $q_i$ since it occurs in just one entry of $T(q)$.
Repeating linear interpolation over all coordinates writes
$\det T(q)$ as a convex combination of the vertex values.
Thus $T(q)$ is nonsingular throughout the box and the extrema
are attained at vertices. The largest determinant comes from
the smallest absolute angle in \eqref{eq:iv-cos}.

Compare the angle $\sum_i\sigma_i\theta_i$ in \eqref{eq:iv-cos}
with $2\delta\sum_i\sigma_i w_i$. This scaled sum is either zero
or has magnitude at least $2\delta$, since the $w_i$ are integers.
To bound the difference, use that the derivative of
$2t-2\arctan t$ lies between $0$ and $2t^2$. Integration gives
\[
 0\le2t-2\arctan t\le\frac23t^3\qquad(t\ge0)
\]
and summing over $t_i=\delta w_i$ bounds the angle error by
\[
 E:=\sum_i|\theta_i-2\delta w_i|
 \le\frac23\delta^3W^3=\frac{\delta}{150W}\le\frac\delta{10}.
\]
If a partition exists, use its signs with $\sigma_m=1$ in
\eqref{eq:iv-cos}. The angle has absolute value at most $E$,
so $\cos x\ge1-x^2/2$ at this vertex gives
\[
 \max_q\det T(q)\ge1-\frac{\delta^2}{200}.
\]
If the weights admit no partition, subtracting the angle error
$E$ from the smallest possible nonzero scaled sum leaves a lower
bound of $2\delta-E\ge19\delta/10$ on every absolute vertex angle.
The bound $\cos x\le1-x^2/3$ on $[-1/5,1/5]$ then gives
\[
 \max_q\det T(q)\le1-\frac{361}{300}\delta^2<1-\delta^2.
\]
To distinguish the two cases, compare the maximum determinant
with the threshold
\[
 \tau=1-\frac{\delta^2}{2}
\]
between the displayed bounds. The rationals specifying the
$4m-2$ rows of $T(q)$ and the threshold $\tau$ have numerators
and denominators bounded by fixed powers of $1+W$. They are
formed by arithmetic on $O(\log(1+W))$-bit integers, so the
reduction takes polynomial time.

\begin{theorem}\label{thm:iv-det}
Deciding whether the maximum determinant of a rational tridiagonal interval matrix is at least a rational threshold is NP-complete. NP-hardness holds even when every matrix in the interval family is nonsingular and every superdiagonal entry is fixed at 1. In particular exact determinant-range computation is NP-hard.
\end{theorem}
\begin{proof}
The weights have a zero signed sum exactly when
$\max_q\det T(q)\ge\tau$, which gives the reduction from
PARTITION~\cite{Karp1972}. For membership in NP, the certificate
specifies a vertex of the input box by choosing an endpoint
of each interval. A vertex suffices because the determinant
is affine in each uncertain entry separately. Exact rational
elimination computes the determinant of the vertex matrix
and checks the threshold in polynomial time.
\end{proof}

\subsection{Transferring hardness to a solution hull}
\begin{theorem}\label{thm:iv-hull}
Computing the exact solution hull of a rational tridiagonal interval system
is NP-hard even for instances with nonempty bounded solution sets in
which every matrix is nonsingular, the right-hand side is fixed at
$-e_N$, and every superdiagonal entry equals one.
\end{theorem}
\begin{proof}
Use the same $T(q)$ and take $b=-e_N$, where $e_N$ is the last
coordinate vector. This makes $x_1=-(T^{-1})_{1N}$. The cofactor
for $(T^{-1})_{1N}$ uses rows $1,\ldots,N-1$ and columns
$2,\ldots,N$ of $T$. In that submatrix the diagonal entries
are the original superdiagonal ones, and entries above the
diagonal vanish by tridiagonality. Its determinant is one,
so the inverse entry is
\[
 (T^{-1})_{1N}=\frac{(-1)^{N+1}}{\det T}.
\]
Since $N=4m-2$ is even, multiplying this entry by $-1$ gives
\[
 x_1=(T^{-1}b)_1=\frac1{\det T}>0.
\]
The map $q\mapsto T(q)^{-1}b$ is continuous on the compact box
because every $T(q)$ is invertible. In particular its range is
nonempty and compact. The formula for $x_1$ shows that the
smallest first coordinate is attained at the largest determinant,
\[
 \underline{x}_1=\frac1{\max_q\det T(q)}.
\]
We can thus decide whether the maximum determinant is at least
$\tau$ by checking whether the lower endpoint is at most $1/\tau$.
\end{proof}

\section{An ordering obstruction for the \texorpdfstring{$q$}{q}-permanent (MI-19)}
\label{sec:mi19}

For $A=(a_{ij})\in\C^{n\times n}$ and $q\in\R$, the
$q$-permanent weights each product of entries by one factor
of $q$ per inversion of its permutation. With $S_n$ the set
of permutations of $\{1,\ldots,n\}$, write
\[
 \operatorname{per}_q(A)=\sum_{\sigma\in S_n}q^{\ell(\sigma)}\prod_{i=1}^n a_{i,\sigma(i)},\qquad
 \ell(\sigma)=\#\{(i,j):1\leq i<j\leq n,\ \sigma(i)>\sigma(j)\}.
\]
The cases $q=-1$ and $q=1$ give the determinant and permanent.
Use $q^0=1$ also at $q=0$. A simultaneous relabeling of rows
and columns conjugates the permutations, which can change
their inversion counts. This matters for the comparison Bapat
and Lal~\cite{BapatLal1994} proposed with the sum over permutations
that preserve a fixed $S\subseteq\{1,\ldots,n\}$ setwise,
\[
 R_{q,S}(A)=\sum_{\substack{\sigma\in S_n\\\sigma(S)=S}}
 q^{\ell(\sigma)}\prod_{i=1}^n a_{i,\sigma(i)}.
\]
The inversion count in $R_{q,S}(A)$ still uses the full index
ordering. Both sums are real for Hermitian $A$ and real $q$ because
$\ell(\sigma^{-1})=\ell(\sigma)$ and the terms for $\sigma$ and
$\sigma^{-1}$ are conjugates.

\begin{problem}{MI-19}
Let $A\in\C^{n\times n}$ be Hermitian positive semidefinite, with
$n\ge2$. For every $q\in[0,1]$ and every nonempty proper subset
$S\subset\{1,\ldots,n\}$, is
\[
 \operatorname{per}_q(A)\ge R_{q,S}(A)?
\]
Da Fonseca states the conjecture \cite[Conjecture~3, equation (5.1)]{Fonseca2018}.
\end{problem}

\begin{theorem}
The subset inequality is false for a real positive semidefinite matrix of order four and rank two. A counterexample is
\[
 A=\begin{pmatrix}
 170&1&167&170\\
 1&1&-2&1\\
 167&-2&173&167\\
 170&1&167&170
 \end{pmatrix},\qquad q=\frac78,\qquad S=\{2\}.
\]
\end{theorem}
\begin{proof}
The matrix is PSD of rank two because it is the Gram matrix of
the columns of the rank-two matrix $X$ below,
\[
 A=X^TX,\qquad X=\begin{pmatrix}13&0&13&13\\1&1&-2&1\end{pmatrix}.
\]
In $\operatorname{per}_q(A)-R_{q,\{2\}}(A)$ the permutations fixing
$2$ cancel, leaving $4!-3!=18$ terms. The only negative matrix
entries are $a_{23}$ and $a_{32}$. A term using just one of these
has a negative product of entries. Collecting terms by inversion
count gives
\begin{align*}
 \operatorname{per}_q(A)-R_{q,\{2\}}(A)
 &=q(q+1)\bigl(115600q^4-229160q^3+315869q^2\\
 &\hspace{8em}-344241q+145010\bigr).
\end{align*}
At $q=7/8$ this equals
$
 -\frac{3235575}{16384}<0.
$
\end{proof}

\begin{remark}[Positive definite inputs]
With $q=7/8$ fixed, the difference in the proof evaluated at
$A+\varepsilon I_4$ is a polynomial in $\varepsilon$. It is
negative at zero and therefore at sufficiently small positive
$\varepsilon$, for which the matrix is positive definite.
\end{remark}

\section{Geometric means and a norm inequality (MI-21)}
\label{sec:mi21}

Take positive semidefinite pairs $A_i,B_i$ with $A_iB_i=B_iA_i$
for each $i$. Audenaert~\cite{Audenaert2015} showed that
$\norm{\sum_i A_iB_i}\le\norm{(\sum_i A_i)(\sum_i B_i)}$ for
every unitarily invariant norm. Freewan and Hayajneh asked
for a comparison of this kind using matrix geometric means.

For Hermitian positive definite $A,B\in\C^{n\times n}$ and
$0\leq t\leq1$, the weighted matrix geometric mean uses spectral
calculus for the real powers in
\[
 A\#_tB=A^{1/2}(A^{-1/2}BA^{-1/2})^tA^{1/2},\qquad A\#B=A\#_{1/2}B.
\]
A unitarily invariant norm
satisfies $\norm{UXV}=\norm{X}$ for all unitary $U,V$.

\begin{problem}{MI-21}
For $m$ pairs of Hermitian positive definite matrices $(A_i,B_i)$ of
order $n$, set $A=\sum_i A_i$ and $B=\sum_i B_i$. Does
Freewan and Hayajneh's inequality~\cite[Conjecture~3]{FH2024}
\begin{equation}\label{eq:mi21-target}
 \left\|\sum_{i=1}^m(A_i^s\#_tB_i^s)^r\right\|
 \leq
 \left\|\left(A^{(1-t)srp/2}B^{tsrp}A^{(1-t)srp/2}\right)^{1/p}\right\|
\end{equation}
hold for all $m,n\ge1$, all such pairs, every unitarily invariant
norm, $0\leq t\leq1$, and $p,r,s>0$ with $sr\geq1$?
\end{problem}

Freewan and Hayajneh proved the inequality for $t=1/2$,
$s\ge2$, $r\ge1$ and $rp\ge1$~\cite[Theorem~3.6]{FreewanHayajnehNormMeans2023}.
The conjecture also allows $s=1/2$ with $r=2$. At those
parameters it can fail in dimension two.

\begin{theorem}\label{thm:mi21}
Inequality \eqref{eq:mi21-target} fails in operator norm for two $2\times2$
summands at $s=t=1/2$ and $r=2$. The same counterexample works for every
$p>0$ and has rational positive definite inputs satisfying $A=B=I_2$.
\end{theorem}

\begin{proof}
We make the right-hand side $1$ for every $p>0$ by arranging $A=B=I_2$. Take
\[
 C=\diag\left(\frac{12}{37},\frac{21}{29}\right),\qquad
 E=\diag\left(\frac{35}{37},\frac{20}{29}\right),\qquad
 S=\frac1{17}\begin{pmatrix}15&8\\8&-15\end{pmatrix}.
\]
The Pythagorean triples give $C^2+E^2=I_2$, and $S$ is an orthogonal reflection, $S=S^T=S^{-1}$. Use the two squares for the first sum and reflect and interchange them for the second:
\[
 A_1=C^2,\quad A_2=E^2,\qquad
 B_1=SE^2S,\quad B_2=SC^2S.
\]
Both sums are $I_2$, and all four inputs are rational and positive definite.

Put $D=SES$ and $G=C\#D$. At $s=t=1/2$ the first mean is $G$, and the second is $E\#(SCS)=SGS$ by symmetry and orthogonal covariance. Since $r=2$, we have $sr=1$ and the matrix on the left is
$
 L=G^2+SG^2S.
$

To evaluate $G^2$ in $L$, first compute the square root of
$K=C^{-1/2}DC^{-1/2}$ in the formula for $G$. Write
\[
 \delta=\sqrt{\det D/\det C},\qquad h=\tr(C^{-1}D)+2\delta.
\]
Since $\det K=\delta^2$ and $\tr K=h-2\delta$, Cayley--Hamilton gives
$(K+\delta I_2)^2=hK$. The matrix $K+\delta I_2$ is positive definite,
so $K^{1/2}=(K+\delta I_2)/\sqrt h$. Congruence by $C^{1/2}$ gives
\begin{equation}\label{eq:mi21-2by2}
 C\#D=\frac{D+\delta C}{\sqrt h},\qquad
 (C\#D)^2=\frac{(D+\delta C)^2}{h}.
\end{equation}

Here $\delta=5/3$ and $h=61697295/8682716$. Substituting into \eqref{eq:mi21-2by2} and forming $G^2+SG^2S$ gives
\begin{equation}\label{eq:mi21-L}
 L=\frac1{12158163}
 \begin{pmatrix}8216600&-985600\\-985600&11912600\end{pmatrix}.
\end{equation}
Since $L=G^2+SG^2S$, it commutes with the reflection $S$. In particular,
it preserves the line spanned by $w=(1,-4)^T$, on which $Sw=-w$.
Multiplication along this line gives
\[
 Lw=\frac{1351000}{1350907}w,
 \qquad\frac{1351000}{1350907}=1+\frac{93}{1350907}>1.
\]
Thus $\opnorm{L}>1=\opnorm{I_2}$, contradicting \eqref{eq:mi21-target}.
\end{proof}

\section{Weighted log-majorization and certified matrix roots (MI-22)}
\label{sec:mi22}

For positive vectors $x,y\in\R^n$ in nonincreasing order,
log-majorization $x\prec_{\log}y$ means $\prod_{j=1}^kx_j\le\prod_{j=1}^ky_j$
at every $1\le k<n$, with equality at $k=n$.
Let $s(M)$ be the singular-value vector in this order and use the
weighted geometric mean $\#_t$ from the preceding section.

\begin{problem}{MI-22}
For every $n\ge1$, every Hermitian positive definite pair
$A,B\in\C^{n\times n}$, and every $0\le t\le1$, does the
log-majorization
\begin{equation}\label{eq:mi22-target}
 s\bigl(A^t(A\#_t B)B^{1-t}\bigr)\prec_{\log}s(AB)
\end{equation}
hold, as conjectured by Ghabries, Abbas, Mourad and Assi
\cite[Conjecture~1.1]{GAMA2021}?
\end{problem}

Equality at $k=n$ follows from
$\det(A\#_t B)=(\det A)^{1-t}(\det B)^t$.

At $t=1/2$, \eqref{eq:mi22-target} is Zou's
inequality~\cite[Theorem~2.10]{Zou2017}. Lemos and Soares
also proved this case and asked about general
$t$~\cite[Corollary~3.4 and the question following it]{LemosSoares2018}.
Ghabries et al.\ strengthen the $t=1/2$ inequality in Theorem~3.1
of their paper. Their Theorem~3.2 compares eigenvalue moduli for
a different pair of products. Shi, Wang and Wei~\cite{ShiWangWei2026}
used Furuta-type inequalities to obtain log-majorizations for the
metric geometric mean and weak log-majorizations for the spectral
geometric mean.

If $A$ and $B$ commute, then $A\#_tB=A^{1-t}B^t$ and the product on
the left is $AB$. In the example below, the failure shows up already at
$k=1$, where the comparison is just between the operator norms of the
two products.

\begin{theorem}\label{thm:mi22}
A counterexample to \eqref{eq:mi22-target} is
\[
 A=\diag\left(256,\frac1{256},1\right),\qquad
 B=\begin{pmatrix}17&-4&0\\-4&16385&-8192\\0&-8192&4096\end{pmatrix},
 \qquad t=\frac18.
\]
In fact,
\[
 \opnorm{A^{1/8}(A\#_{1/8}B)B^{7/8}}>10900,
 \qquad \opnorm{AB}<10200.
\]
\end{theorem}

For the lower bound in Theorem~\ref{thm:mi22}, we estimate the absolute
value of the $(1,2)$ entry of $A^{1/8}(A\#_{1/8}B)B^{7/8}$. Computing
this product involves two principal matrix eighth roots. For a
rational approximation $R$ to the positive definite eighth root
of a rational positive definite matrix $C$, we can compute $R^8-C$
exactly. The following estimate turns that residual
into a bound on $R-C^{1/8}$.

\begin{lemma}\label{lem:root-lipschitz}
If $P,Q\in\C^{n\times n}$ are Hermitian with $P,Q\succeq mI$ for some $m>0$, and $0<\alpha<1$, then
\[
 \opnorm{P^\alpha-Q^\alpha}
 \leq \alpha m^{\alpha-1}\opnorm{P-Q}.
\]
\end{lemma}
\begin{proof}
The integral representation of $x^\alpha$ and the resolvent identity give
\[
 P^\alpha-Q^\alpha
 =\frac{\sin(\pi\alpha)}{\pi}\int_0^\infty
 u^\alpha(uI+Q)^{-1}(P-Q)(uI+P)^{-1}\,\mathrm{d}u.
\]
The norm of the integrand is bounded by the integrable function $u^\alpha(u+m)^{-2}\opnorm{P-Q}$. Differentiating the scalar integral formula at $m$ gives
\[
 \frac{\sin(\pi\alpha)}{\pi}\int_0^\infty\frac{u^\alpha}{(u+m)^2}\,\mathrm{d}u
 =\alpha m^{\alpha-1},
\]
which proves the estimate.
\end{proof}

\begin{proof}[Proof of Theorem \ref{thm:mi22}]
Here $A\succ0$, and the leading principal minors $17$, $278529$, and $4096$ give $B\succ0$. For the right-hand side,
\[
 AB=\begin{pmatrix}4352&-1024&0\\-1/64&16385/256&-32\\0&-8192&4096\end{pmatrix},
\]
so its Frobenius norm gives the required upper bound:
\[
 \norm{AB}_{\mathrm F}^2=\frac{6807858741265}{65536}<10200^2.
\]

For the left-hand side, put
\[
 C=A^{-1/2}BA^{-1/2}
 =\begin{pmatrix}
 17/256&-4&0\\-4&4194560&-131072\\0&-131072&4096
 \end{pmatrix}.
\]
With $X=C^{1/8}$ and $Y=B^{1/8}$, the product to be bounded below is
\[
 L=A^{5/8}XA^{1/2}Y^7,
 \qquad A^{5/8}=\diag(32,1/32,1),\quad A^{1/2}=\diag(16,1/16,1).
\]
High-precision calculations supplied the rational approximations $R$ to $X$ and $S$ to $Y$:
\begingroup\small
\[
 R=10^{-15}\begin{pmatrix}
 563431071954661&-4774979464975&-152620714401404\\
 -4774979464975&6722248446399974&-185449303604146\\
 -152620714401404&-185449303604146&793308473748417
 \end{pmatrix},
\]
\[
 S=10^{-15}\begin{pmatrix}
 1417565567728212&-40084698659651&-79375309604622\\
 -40084698659651&2883518412755210&-1150490885652445\\
 -79375309604622&-1150490885652445&1157680400969996
 \end{pmatrix}.
\]
\endgroup
To bound the root errors, we need lower eigenvalue bounds for $C,R^8$
and $B,S^8$. Upper bounds for $R,S$ control how those errors affect the
product. Check the spectral inequalities by applying Sylvester's
criterion to the matrix differences, and compute the squared
Frobenius norms by summing squares of the residual entries. In exact
rational arithmetic this gives
\begin{align}
 B,C&\succ2^{-10}I,\label{eq:mi22-cert1}\\
 \frac9{20}I\prec R\prec8I,&\qquad \frac9{20}I\prec S\prec4I,\label{eq:mi22-cert2}\\
 \norm{R^8-C}_{\mathrm F}^2&<10^{-16},\qquad \norm{S^8-B}_{\mathrm F}^2<10^{-16}.\label{eq:mi22-cert3}
\end{align}

The lower bounds on $R,S$ give $R^8,S^8\succ2^{-10}I$, since $(9/20)^8>2^{-10}$.
The positive definite eighth roots of these powers are $R$ and $S$,
respectively. The residuals have operator norm less than $10^{-8}$.
Apply Lemma~\ref{lem:root-lipschitz} to the pairs $(C,R^8)$ and $(B,S^8)$
with $\alpha=1/8$ and $m=2^{-10}$. The multiplier in both bounds is
$\alpha m^{\alpha-1}=2^{23/4}<64$, so
\[
 \opnorm{X-R}<\varepsilon,\qquad \opnorm{Y-S}<\varepsilon,
 \qquad \varepsilon=64\cdot10^{-8}.
\]

Since $\tr B<4^8$, we have $\opnorm{Y}<4$. We also have $\opnorm{S}<4$ by \eqref{eq:mi22-cert2}, so telescoping gives
\[
 \opnorm{Y^7-S^7}\leq7\cdot4^6\varepsilon.
\]
Set $\widetilde L=A^{5/8}RA^{1/2}S^7$. The change in any entry is
bounded by $\opnorm{L-\widetilde L}$. Replace $X$ by $R$ first and
then $Y^7$ by $S^7$, using $\opnorm{R}<8$, to get
\begin{align*}
 \opnorm{L-\widetilde L}
 &\leq32\cdot16\,\varepsilon\bigl(4^7+8\cdot7\cdot4^6\bigr)\\
 &=\frac{6291456}{78125}<100.
\end{align*}
Exact multiplication gives
\begin{equation}\label{eq:mi22-entry}
 \widetilde L_{12}>11000.
\end{equation}
The bounds $\opnorm{L}\geq|L_{12}|>11000-100=10900$ and
$\opnorm{AB}\leq\norm{AB}_{\mathrm F}<10200$ reverse the $k=1$
comparison in \eqref{eq:mi22-target}.
\end{proof}

\small

\begingroup
\interlinepenalty=10000
\bibliographystyle{siamplain}
\bibliography{references}
\endgroup
\end{document}